\documentclass{article}
\usepackage{amsmath}
\usepackage{amssymb}
\begin{document}
%%%%%%%%%%%%%%%%%%%%%%%%%%%%%%%%%%%%%%%%%%%%%%%%%%%%%%%%%%%%%%%%%%%%%%%
%%%%%%%%%%%%%%%%%%%%%%%%%     Macros      %%%%%%%%%%%%%%%%%%%%%%%%%%%%%
%%%%%%%%%%%%%%%%%%%%%%%%%%%%%%%%%%%%%%%%%%%%%%%%%%%%%%%%%%%%%%%%%%%%%%%
\def\e#1\e{\begin{equation}#1\end{equation}}
\def\ea#1\ea{\begin{align}#1\end{align}}
\def\eq#1{{\rm(\ref{#1})}}
\newtheorem{thm}{Theorem}[section]
\newtheorem{lem}[thm]{Lemma}
\newtheorem{prop}[thm]{Proposition}
\newtheorem{cor}[thm]{Corollary}
\newenvironment{dfn}{\medskip\refstepcounter{thm}
\noindent{\bf Definition \thesection.\arabic{thm}\ }}{\medskip}
\newenvironment{proof}[1][,]{\medskip\ifcat,#1
\noindent{\it Proof.\ }\else\noindent{\it Proof of #1.\ }\fi}
{\hfill$\square$\medskip}
\def\dim{\mathop{\rm dim}}
\def\Re{\mathop{\rm Re}}
\def\Im{\mathop{\rm Im}}
\def\Ker{\mathop{\rm Ker}}
\def\diam{\mathop{\rm diam}}
\def\vol{\mathop{\rm vol}}
\def\supp{\mathop{\rm supp}}
\def\id{\mathop{\rm id}}
\def\SU{\mathop{\rm SU}}
\def\ge{\geqslant} 
\def\le{\leqslant} 
\def\N{\mathbin{\mathbb N}}
\def\R{\mathbin{\mathbb R}}
\def\Z{\mathbin{\mathbb Z}}
\def\C{\mathbin{\mathbb C}}
\def\g{\mathfrak{g}} 
\def\su{\mathfrak{su}} 
\def\al{\alpha}
\def\be{\beta}
\def\ga{\gamma}
\def\de{\delta}
\def\ep{\epsilon}
\def\ve{\varepsilon}
\def\la{\lambda}
\def\ka{\kappa}
\def\th{\theta}
\def\ze{\zeta}
\def\up{\upsilon}
\def\vp{\varphi}
\def\si{\sigma}
\def\om{\omega}
\def\De{\Delta}
\def\La{\Lambda}
\def\Om{\Omega}
\def\Ga{\Gamma}
\def\Si{\Sigma}
\def\Th{\Theta}
\def\Up{\Upsilon}
\def\d{{\rm d}}
\def\pd{\partial}
\def\ts{\textstyle}
\def\sst{\scriptscriptstyle}
\def\w{\wedge}
\def\lt{\ltimes}
\def\sm{\setminus}
\def\op{\oplus}
\def\ot{\otimes}
\def\bigot{\bigotimes}
\def\iy{\infty}
\def\ra{\rightarrow}
\def\hookra{\hookrightarrow}
\def\longra{\longrightarrow}
\def\t{\times}
\def\na{\nabla}
\def\ha{{\textstyle\frac{1}{2}}}
\def\ti{\tilde}
\def\ov{\overline}
\def\sC{{\smash{\sst C}}}
\def\sCi{{\smash{\sst C_i}}}
\def\sL{{\smash{\sst L}}}
\def\sLi{{\smash{\sst L_i}}}
\def\sSi{{\smash{\sst\Si}}}
\def\sSii{{\smash{\sst\Si_i}}}
\def\sX{{\smash{\sst X}}}
\def\sNt{{\smash{\sst N^t}}}
\def\sW{{\smash{\sst W}}}
\def\sWt{{\smash{\sst W^t}}}
\def\sXp{{\smash{\sst X'}}}
\def\A{{\cal A}}
\def\B{{\cal B}}
\def\D{{\cal D}}
\def\H{{\cal H}}
\def\I{{\cal I}}
\def\K{{\cal K}}
\def\O{{\cal O}}
\def\ms#1{\vert#1\vert^2}
\def\bms#1{\bigl\vert#1\bigr\vert^2}
\def\md#1{\vert #1 \vert}
\def\bmd#1{\big\vert #1 \big\vert}
\def\cnm#1#2{\Vert #1 \Vert_{C^{#2}}} %% C^k or C^{k,\al} norms
\def\lnm#1#2{\Vert #1 \Vert_{L^{#2}}} %% L^q norms
\def\snm#1#2#3{\Vert #1 \Vert_{L^{#2}_{#3}}} %% L^q_k norms
\def\bcnm#1#2{\bigl\Vert #1 \bigr\Vert_{C^{#2}}} %% big C^k etc norms.
\def\blnm#1#2{\bigl\Vert #1 \bigr\Vert_{L^{#2}}} %% big L^q norms
\def\bsnm#1#2#3{\bigl\Vert #1 \bigr\Vert_{L^{#2}_{#3}}} %% big L^q_k norms
\def\an#1{\langle#1\rangle}
\def\ban#1{\bigl\langle#1\bigr\rangle}
%%%%%%%%%%%%%%%%%%%%%%%%%%%%%%%%%%%%%%%%%%%%%%%%%%%%%%%%%%%%%%%%%%%%%%%
%%%%%%%%%%%%%%%%%%%%%     Text of paper    %%%%%%%%%%%%%%%%%%%%%%%%%%%%
%%%%%%%%%%%%%%%%%%%%%%%%%%%%%%%%%%%%%%%%%%%%%%%%%%%%%%%%%%%%%%%%%%%%%%%
\title{Special Lagrangian submanifolds with isolated\\
conical singularities. III. Desingularization, \\
the unobstructed case}
\author{Dominic Joyce \\ Lincoln College, Oxford}
\date{}
\maketitle

\section{Introduction}
\label{cu1}

{\it Special Lagrangian $m$-folds (SL\/ $m$-folds)} are a
distinguished class of real $m$-dimensional minimal submanifolds
which may be defined in $\C^m$, or in {\it Calabi--Yau $m$-folds},
or more generally in {\it almost Calabi--Yau $m$-folds} (compact
K\"ahler $m$-folds with trivial canonical bundle). We write an
almost Calabi--Yau $m$-fold as $M$ or $(M,J,\om,\Om)$, where
the manifold $M$ has complex structure $J$, K\"ahler form
$\om$ and holomorphic volume form~$\Om$.

This is the third in a series of five papers
\cite{Joyc4,Joyc5,Joyc6,Joyc7} studying SL $m$-folds with
{\it isolated conical singularities}. That is, we consider an
SL $m$-fold $X$ in an almost Calabi--Yau $m$-fold $M$ for $m>2$
with singularities at $x_1,\ldots,x_n$ in $M$, such that for
some special Lagrangian cones $C_i$ in $T_{\smash{x_i}}M\cong\C^m$
with $C_i\sm\{0\}$ nonsingular, $X$ approaches $C_i$ near $x_i$ in
an asymptotic $C^1$ sense. Readers are advised to begin with
the final paper \cite{Joyc7}, which surveys the series, and
applies the results to prove some conjectures.

The first paper \cite{Joyc4} laid the foundations for the series,
and studied the {\it regularity} of SL $m$-folds with conical
singularities near their singular points. The second paper
\cite{Joyc5} discussed the {\it deformation theory} of compact
SL $m$-folds $X$ with conical singularities in an almost
Calabi--Yau $m$-fold~$M$.

This paper and the sequel \cite{Joyc6} study
{\it desingularizations} of compact SL $m$-folds $X$ with
conical singularities. That is, we construct a family of
compact, {\it nonsingular} SL $m$-folds $\smash{\ti N^t}$
in $M$ for $t\in(0,\ep]$ such that $\smash{\ti N^t}\ra X$
as $t\ra 0$, in the sense of currents.

Having a good understanding of the singularities of special
Lagrangian submanifolds will be essential in clarifying the
Strominger--Yau--Zaslow conjecture on the Mirror Symmetry
of Calabi--Yau 3-folds \cite{SYZ}, and also in resolving
conjectures made by the author \cite{Joyc1} on defining
new invariants of Calabi--Yau 3-folds by counting special
Lagrangian homology 3-spheres with weights. The series aims
to develop such an understanding for simple singularities
of SL $m$-folds.

Here is the basic idea of the paper. Let $X$ be a compact SL
$m$-fold with conical singularities $x_1,\ldots,x_n$ in an
almost Calabi--Yau $m$-fold $(M,J,\om,\Om)$. Choose an
isomorphism $\up_i:\C^m\ra T_{x_i}M$ for $i=1,\ldots,n$.
Then there is a unique {\it SL cone} $C_i$ in $\C^m$ with $X$
asymptotic to $\up_i(C_i)$ at~$x_i$.

Let $L_i$ be an {\it Asymptotically Conical SL\/ $m$-fold\/}
({\it AC SL\/ $m$-fold\/}) in $\C^m$, asymptotic to $C_i$ at
infinity. As $C_i$ is a cone it is invariant under dilations,
so $tC_i=C_i$ for all $t>0$. Thus $tL_i=\{t\,{\bf x}:{\bf x}\in
L_i\}$ is also an AC SL $m$-fold asymptotic to $C_i$ for
$t>0$. We explicitly construct a 1-parameter family of compact,
nonsingular {\it Lagrangian} $m$-folds $N^t$ in $(M,\om)$ for
$t\in(0,\de)$ by gluing $tL_i$ into $X$ at $x_i$, using a
partition of unity.

When $t$ is small, $N^t$ is close to being special Lagrangian (its
phase is nearly constant), but also close to being singular (it has
large curvature and small injectivity radius). We prove using
analysis that for small $\ep\in(0,\de)$ we can deform $N^t$ to a
{\it special\/} Lagrangian $m$-fold $\smash{\ti N^t}$ in $M$ for
all $t\in(0,\ep]$, using a small Hamiltonian deformation. The proof
involves a delicate balancing act, showing that the advantage of
being close to special Lagrangian outweighs the disadvantage of
being nearly singular.

Here are some of the issues involved in doing this in full
generality:
\begin{itemize}
\item[(i)] To ensure $N^t$ and $\smash{\ti N^t}$ are connected,
we suppose $X$ is connected. But $X'=X\sm\{x_1,\ldots,x_n\}$, the
nonsingular part of $X$, may not be connected. If it is not then
the Laplacian $\De$ on $N^t$ has small positive eigenvalues, of
size $O(t^{m-2})$. These cause analytic problems in
constructing~$\smash{\ti N^t}$.
\item[(ii)] Let $\Si_i=C_i\cap{\cal S}^{2m-1}$. Then $\Si_i$ is a
compact $(m\!-\!1)$-manifold, and $L_i$ effectively has boundary
$\Si_i$ at infinity. There are natural {\it cohomological invariants}
$Y(L_i)\in H^1(\Si_i,\R)$ and $Z(L_i)\in H^{m-1}(\Si_i,\R)$. It
turns out that there are {\it topological obstructions} to the
existence of $N^t$ or $\smash{\ti N^t}$, involving the $Y(L_i)$
and~$Z(L_i)$.
\item[(iii)] Let $\bigl\{(M,J^s,\om^s,\Om^s):s\in{\cal F}\bigr\}$
be a smooth family of almost Calabi--Yau $m$-folds for $0\in{\cal
F}\subset\R^d$ with $(M,J^0,\om^0,\Om^0)=(M,J,\om,\Om)$. Then we
can consider special Lagrangian desingularizations
$\smash{\ti N^{s,t}}$ of $X$ not just in $(M,J,\om,\Om)$
but in $(M,J^s,\om^s,\Om^s)$ for small $s\in{\cal F}$. To do
this introduces new analytic problems, and new topological
obstructions involving the cohomology classes $[\om^s]$
and~$[\Im\Om^s]$.
\end{itemize}

Rather than tackling these questions all at once, we prove
our first main result in \S\ref{cu6} assuming that $X'$ is
connected, that $Y(L_i)=0$ and $L_i$ converges quickly to
$C_i$ at infinity in $\C^m$, and working in a single almost
Calabi--Yau $m$-fold $(M,J,\om,\Om)$ rather than a family.
This simplifies (i)--(iii) above.

Section \ref{cu7} extends this to the case when $X'$ is not
connected, as in (i), but still supposing $Y(L_i)=0$ and
$L_i$ converges quickly to $C_i$. The sequel \cite{Joyc6}
deals with issues (ii) and (iii), allowing $Y(L_i)\ne 0$ and
$L_i$ to converge more slowly to $C_i$, and working in a
family of almost Calabi--Yau $m$-folds~$(M,J^s,\om^s,\Om^s)$.

We begin in \S\ref{cu2} with an introduction to special
Lagrangian geometry. Sections \ref{cu3} and \ref{cu4}
discuss SL $m$-folds with conical singularities and
Asymptotically Conical SL $m$-folds respectively,
recalling results we will need from~\cite{Joyc4}.

Given a compact Lagrangian $m$-fold $N$ in an almost
Calabi--Yau $m$-fold $(M,J,\om,\Om)$ which is close to being
special Lagrangian, \S\ref{cu5} uses analysis to construct an
SL $m$-fold $\ti N$ as a small Hamiltonian deformation of $N$.
This existence result, Theorem \ref{cu5thm2} below, can probably
be used elsewhere. In each of \S\ref{cu6} and \S\ref{cu7} we
construct a family of Lagrangian $m$-folds $N^t$ in $(M,J,\om,\Om)$,
and apply Theorem \ref{cu5thm2} to show that $N^t$ can be deformed
to an SL $m$-fold $\smash{\ti N^t}$ for small~$t$.

For simplicity we generally take all submanifolds to be
{\it embedded}. However, all our results generalize immediately
to {\it immersed\/} submanifolds, with only cosmetic changes.

We conclude by discussing similar work by other authors.
Salur \cite{Salu1,Salu2} considers a nonsingular, connected,
immersed SL 3-fold $N$ in a Calabi--Yau 3-fold with a
codimension two self-intersection along an ${\cal S}^1$, and
constructs new SL 3-folds by smoothing along the~${\cal S}^1$.

Butscher \cite{Buts} studies SL $m$-folds $N$ in $\C^m$
with boundary in a symplectic submanifold $W^{2m-2}
\subset\C^m$. Given two such SL $m$-folds $N_1,N_2$
intersecting transversely at $x$ and satisfying an
angle criterion, he constructs a 1-parameter family
of connect sum SL $m$-folds $N_1\#_xN_2$ in $\C^m$,
with boundary, by gluing in an explicit AC SL $m$-fold
$L$ in $\C^m$ due to Lawlor \cite{Lawl}, diffeomorphic
to ${\cal S}^{m-1}\t\R$ and asymptotic to the union of
two SL planes $\R^m$ in~$\C^m$.

Closest to the present paper is Lee \cite{Lee}. She
considers a compact, connected, immersed SL $m$-fold
$N$ in a Calabi--Yau $m$-fold $M$ with transverse
self-intersection points $x_1,\ldots,x_n$ satisfying
an angle criterion. She shows that $N$ can be
desingularized by gluing in one of Lawlor's AC SL
$m$-folds $L_i$ at $x_i$ for $i=1,\ldots,n$, to get
a family of compact, embedded SL $m$-folds in $M$.
Her result follows from Theorem \ref{cu6thm5} below.
\medskip

\noindent{\it Acknowledgements.} I would like to thank
Stephen Marshall, Sema Salur and Adrian Butscher for useful
conversations. I was supported by an EPSRC Advanced Research
Fellowship whilst writing this paper.

\section{Special Lagrangian geometry}
\label{cu2}

We introduce special Lagrangian submanifolds (SL $m$-folds)
in two different geometric contexts. First we define SL
$m$-folds in $\C^m$. Then we discuss SL $m$-folds in
{\it almost Calabi--Yau $m$-folds}, compact K\"ahler
manifolds equipped with a holomorphic volume form, which
generalize Calabi--Yau manifolds. Some references for
this section are Harvey and Lawson \cite{HaLa} and
the author \cite{Joyc3}. We begin by defining {\it
calibrations} and {\it calibrated submanifolds},
following~\cite{HaLa}.

\begin{dfn} Let $(M,g)$ be a Riemannian manifold. An {\it oriented
tangent $k$-plane} $V$ on $M$ is a vector subspace $V$ of
some tangent space $T_xM$ to $M$ with $\dim V=k$, equipped
with an orientation. If $V$ is an oriented tangent $k$-plane
on $M$ then $g\vert_V$ is a Euclidean metric on $V$, so 
combining $g\vert_V$ with the orientation on $V$ gives a 
natural {\it volume form} $\vol_V$ on $V$, which is a 
$k$-form on~$V$.

Now let $\vp$ be a closed $k$-form on $M$. We say that
$\vp$ is a {\it calibration} on $M$ if for every oriented
$k$-plane $V$ on $M$ we have $\vp\vert_V\le \vol_V$. Here
$\vp\vert_V=\al\cdot\vol_V$ for some $\al\in\R$, and 
$\vp\vert_V\le\vol_V$ if $\al\le 1$. Let $N$ be an 
oriented submanifold of $M$ with dimension $k$. Then 
each tangent space $T_xN$ for $x\in N$ is an oriented
tangent $k$-plane. We say that $N$ is a {\it calibrated 
submanifold\/} if $\vp\vert_{T_xN}=\vol_{T_xN}$ for all~$x\in N$.
\label{cu2def1}
\end{dfn}

It is easy to show that calibrated submanifolds are automatically
{\it minimal submanifolds} \cite[Th.~II.4.2]{HaLa}. Here is the 
definition of special Lagrangian submanifolds in $\C^m$, taken
from~\cite[\S III]{HaLa}.

\begin{dfn} Let $\C^m$ have complex coordinates $(z_1,\dots,z_m)$, 
and define a metric $g'$, a real 2-form $\om'$ and a complex $m$-form 
$\Om'$ on $\C^m$ by
\e
\begin{split}
g'=\ms{\d z_1}+\cdots+\ms{\d z_m},\quad
\om'&=\ts\frac{i}{2}(\d z_1\w\d\bar z_1+\cdots+\d z_m\w\d\bar z_m),\\
\text{and}\quad\Om'&=\d z_1\w\cdots\w\d z_m.
\end{split}
\label{cu2eq1}
\e
Then $\Re\Om'$ and $\Im\Om'$ are real $m$-forms on $\C^m$. Let $L$
be an oriented real submanifold of $\C^m$ of real dimension $m$. We
say that $L$ is a {\it special Lagrangian submanifold\/} of $\C^m$,
or {\it SL\/ $m$-fold}\/ for short, if $L$ is calibrated with respect
to $\Re\Om'$, in the sense of Definition~\ref{cu2def1}.
\label{cu2def2}
\end{dfn}

Harvey and Lawson \cite[Cor.~III.1.11]{HaLa} give the following
alternative characterization of special Lagrangian submanifolds:

\begin{prop} Let\/ $L$ be a real $m$-dimensional submanifold 
of\/ $\C^m$. Then $L$ admits an orientation making it into an
SL submanifold of\/ $\C^m$ if and only if\/ $\om'\vert_L\equiv 0$ 
and\/~$\Im\Om'\vert_L\equiv 0$.
\label{cu2prop1}
\end{prop}

An $m$-dimensional submanifold $L$ in $\C^m$ is called {\it Lagrangian} 
if $\om'\vert_L\equiv 0$. Thus special Lagrangian submanifolds are 
Lagrangian submanifolds satisfying the extra condition that 
$\Im\Om'\vert_L\equiv 0$, which is how they get their name. We
shall define special Lagrangian submanifolds not just in
Calabi--Yau manifolds, but in the much larger
class of {\it almost Calabi--Yau manifolds}.

\begin{dfn} Let $m\ge 2$. An {\it almost Calabi--Yau $m$-fold\/}
is a quadruple $(M,J,\om,\Om)$ such that $(M,J)$ is a compact
$m$-dimensional complex manifold, $\om$ is the K\"ahler form
of a K\"ahler metric $g$ on $M$, and $\Om$ is a non-vanishing
holomorphic $(m,0)$-form on~$M$.

We call $(M,J,\om,\Om)$ a {\it Calabi--Yau $m$-fold\/} if in
addition $\om$ and $\Om$ satisfy
\e
\om^m/m!=(-1)^{m(m-1)/2}(i/2)^m\Om\w\bar\Om.
\label{cu2eq2}
\e
Then for each $x\in M$ there exists an isomorphism $T_xM\cong\C^m$
that identifies $g_x,\om_x$ and $\Om_x$ with the flat versions
$g',\om',\Om'$ on $\C^m$ in \eq{cu2eq1}. Furthermore, $g$ is
Ricci-flat and its holonomy group is a subgroup of~$\SU(m)$.
\label{cu2def3}
\end{dfn}

This is not the usual definition of a Calabi--Yau manifold, but
is essentially equivalent to it.

\begin{dfn} Let $(M,J,\om,\Om)$ be an almost Calabi--Yau $m$-fold,
and $N$ a real $m$-dimensional submanifold of $M$. We call $N$ a
{\it special Lagrangian submanifold}, or {\it SL $m$-fold\/} for
short, if $\om\vert_N\equiv\Im\Om\vert_N\equiv 0$. It easily
follows that $\Re\Om\vert_N$ is a nonvanishing $m$-form on $N$.
Thus $N$ is orientable, with a unique orientation in which
$\Re\Om\vert_N$ is positive.
\label{cu2def4}
\end{dfn}

Again, this is not the usual definition of SL $m$-fold, but is
essentially equivalent to it. Suppose $(M,J,\om,\Om)$ is an
almost Calabi--Yau $m$-fold, with metric $g$. Let
$\psi:M\ra(0,\iy)$ be the unique smooth function such that
\e
\psi^{2m}\om^m/m!=(-1)^{m(m-1)/2}(i/2)^m\Om\w\bar\Om,
\label{cu2eq3}
\e
and define $\ti g$ to be the conformally equivalent metric $\psi^2g$
on $M$. Then $\Re\Om$ is a {\it calibration} on the Riemannian manifold
$(M,\ti g)$, and SL $m$-folds $N$ in $(M,J,\om,\Om)$ are calibrated
with respect to it, so that they are minimal with respect to~$\ti g$.

If $M$ is a Calabi--Yau $m$-fold then $\psi\equiv 1$ by \eq{cu2eq2},
so $\ti g=g$, and an $m$-submanifold $N$ in $M$ is special Lagrangian
if and only if it is calibrated w.r.t.\ $\Re\Om$ on $(M,g)$, as in
Definition \ref{cu2def2}. This recovers the usual definition of
special Lagrangian $m$-folds in Calabi--Yau $m$-folds.

\section{SL $m$-folds with conical singularities}
\label{cu3}

The preceding papers \cite{Joyc4,Joyc5} studied SL $m$-folds
$X$ with {\it conical singularities} in an almost Calabi--Yau
$m$-fold $(M,J,\om,\Om)$. We now recall the definitions and
results from \cite{Joyc4} that we will need later. For brevity
we keep explanations to a minimum, and readers are referred to
\cite{Joyc4} for further details.

\subsection{Preliminaries on special Lagrangian cones}
\label{cu31}

Following \cite[\S 2.1]{Joyc4} we give definitions
and results on {\it special Lagrangian cones}.

\begin{dfn} A (singular) SL $m$-fold $C$ in $\C^m$ is called a
{\it cone} if $C=tC$ for all $t>0$, where $tC=\{t\,{\bf x}:{\bf x}
\in C\}$. Let $C$ be a closed SL cone in $\C^m$ with an isolated
singularity at 0. Then $\Si=C\cap{\cal S}^{2m-1}$ is a compact,
nonsingular $(m\!-\!1)$-submanifold of ${\cal S}^{2m-1}$, not
necessarily connected. Let $g_\sSi$ be the restriction
of $g'$ to $\Si$, where $g'$ is as in~\eq{cu2eq1}.

Set $C'=C\sm\{0\}$. Define $\iota:\Si\t(0,\iy)\ra\C^m$ by
$\iota(\si,r)=r\si$. Then $\iota$ has image $C'$. By an abuse
of notation, {\it identify} $C'$ with $\Si\t(0,\iy)$ using
$\iota$. The {\it cone metric} on $C'\cong\Si\t(0,\iy)$
is~$g'=\iota^*(g')=\d r^2+r^2g_\sSi$.

For $\al\in\R$, we say that a function $u:C'\ra\R$ is
{\it homogeneous of order} $\al$ if $u\circ t\equiv t^\al u$ for
all $t>0$. Equivalently, $u$ is homogeneous of order $\al$ if
$u(\si,r)\equiv r^\al v(\si)$ for some function~$v:\Si\ra\R$.
\label{cu3def1}
\end{dfn}

In \cite[Lem.~2.3]{Joyc4} we study {\it homogeneous harmonic
functions} on~$C'$.

\begin{lem} In the situation of Definition \ref{cu3def1},
let\/ $u(\si,r)\equiv r^\al v(\si)$ be a homogeneous function
of order $\al$ on $C'=\Si\t(0,\iy)$, for $v\in C^2(\Si)$. Then
\begin{equation*}
\De u(\si,r)=r^{\al-2}\bigl(\De_\sSi v-\al(\al+m-2)v\bigr),
\end{equation*}
where $\De$, $\De_\sSi$ are the Laplacians on $(C',g')$
and\/ $(\Si,g_\sSi)$. Hence, $u$ is harmonic on $C'$
if and only if\/ $v$ is an eigenfunction of\/ $\De_\sSi$
with eigenvalue~$\al(\al+m-2)$.
\label{cu3lem}
\end{lem}

Following \cite[Def.~2.5]{Joyc4}, we define:

\begin{dfn} In the situation of Definition \ref{cu3def1},
suppose $m>2$ and define
\e
\D_\sSi=\bigl\{\al\in\R:\text{$\al(\al+m-2)$ is
an eigenvalue of $\De_\sSi$}\bigr\}.
\label{cu3eq1}
\e
Then $\D_\sSi$ is a countable, discrete subset of
$\R$. By Lemma \ref{cu3lem}, an equivalent definition is that
$\D_\sSi$ is the set of $\al\in\R$ for which there
exists a nonzero homogeneous harmonic function $u$ of order
$\al$ on~$C'$.

Define $m_\sSi:\D_\sSi\ra\N$ by taking
$m_\sSi(\al)$ to be the multiplicity of the eigenvalue
$\al(\al+m-2)$ of $\De_\sSi$, or equivalently the
dimension of the vector space of homogeneous harmonic
functions $u$ of order $\al$ on $C'$. Define
$N_\sSi:\R\ra\Z$ by
\begin{equation*}
N_\sSi(\de)=
-\sum_{\!\!\!\!\al\in\D_\sSi\cap(\de,0)\!\!\!\!}m_\sSi(\al)
\;\>\text{if $\de<0$, and}\;\>
N_\sSi(\de)=
\sum_{\!\!\!\!\al\in\D_\sSi\cap[0,\de]\!\!\!\!}m_\sSi(\al)
\;\>\text{if $\de\ge 0$.}
\end{equation*}
Then $N_\sSi$ is monotone increasing and upper semicontinuous,
and is discontinuous exactly on $\D_\sSi$, increasing by
$m_\sSi(\al)$ at each $\al\in\D_\sSi$. As the
eigenvalues of $\De_\sSi$ are nonnegative, we see that
$\D_\sSi\cap(2-m,0)=\emptyset$ and $N_\sSi\equiv 0$
on~$(2-m,0)$.
\label{cu3def2}
\end{dfn}

\subsection{The definition of SL $m$-folds with conical singularities}
\label{cu32}

Now we can define {\it conical singularities} of SL $m$-folds,
following~\cite[Def.~3.6]{Joyc4}.

\begin{dfn} Let $(M,J,\om,\Om)$ be an almost Calabi--Yau $m$-fold
for $m>2$, and define $\psi:M\ra(0,\iy)$ as in \eq{cu2eq3}. Suppose
$X$ is a compact singular SL $m$-fold in $M$ with singularities at
distinct points $x_1,\ldots,x_n\in X$, and no other singularities.

Fix isomorphisms $\up_i:\C^m\ra T_{x_i}M$ for $i=1,\ldots,n$
such that $\up_i^*(\om)=\om'$ and $\up_i^*(\Om)=\psi(x_i)^m\Om'$,
where $\om',\Om'$ are as in \eq{cu2eq1}. Let $C_1,\ldots,C_n$ be SL
cones in $\C^m$ with isolated singularities at 0. For $i=1,\ldots,n$
let $\Si_i=C_i\cap{\cal S}^{2m-1}$, and let $\mu_i\in(2,3)$ with
$(2,\mu_i]\cap\D_\sSii=\emptyset$, where $\D_\sSii$ is defined in
\eq{cu3eq1}. Then we say that $X$ has a {\it conical singularity}
at $x_i$, with {\it rate} $\mu_i$ and {\it cone} $C_i$ for
$i=1,\ldots,n$, if the following holds.

By Darboux' Theorem \cite[Th.~3.15]{McSa} there exist embeddings
$\Up_i:B_R\ra M$ for $i=1,\ldots,n$ satisfying $\Up_i(0)=x_i$,
$\d\Up_i\vert_0=\up_i$ and $\Up_i^*(\om)=\om'$, where $B_R$ is
the open ball of radius $R$ about 0 in $\C^m$ for some small $R>0$.
Define $\iota_i:\Si_i\t(0,R)\ra B_R$ by $\iota_i(\si,r)=r\si$
for~$i=1,\ldots,n$.

Define $X'=X\sm\{x_1,\ldots,x_n\}$. Then there should exist a
compact subset $K\subset X'$ such that $X'\sm K$ is a union of
open sets $S_1,\ldots,S_n$ with $S_i\subset\Up_i(B_R)$, whose
closures $\bar S_1,\ldots,\bar S_n$ are disjoint in $X$. For
$i=1,\ldots,n$ and some $R'\in(0,R]$ there should exist a smooth
$\phi_i:\Si_i\t(0,R')\ra B_R$ such that $\Up_i\circ\phi_i:\Si_i
\t(0,R')\ra M$ is a diffeomorphism $\Si_i\t(0,R')\ra S_i$, and
\e
\bmd{\na^k(\phi_i-\iota_i)}=O(r^{\mu_i-1-k})
\quad\text{as $r\ra 0$ for $k=0,1$.}
\label{cu3eq2}
\e
Here $\na,\md{\,.\,}$ are computed using the cone
metric $\iota_i^*(g')$ on~$\Si_i\t(0,R')$.
\label{cu3def3}
\end{dfn}

The reasoning behind this definition was discussed in
\cite[\S 3.3]{Joyc4}. We suppose $m>2$ for two reasons.
Firstly, the only SL cones in $\C^2$ are finite unions
of SL planes $\R^2$ in $\C^2$ intersecting only at 0.
Thus any SL 2-fold with conical singularities is actually
{\it nonsingular} as an immersed 2-fold. Secondly, $m=2$
is a special case in the analysis of \cite[\S 2]{Joyc4},
and it is simpler to exclude it. Therefore we will suppose
$m>2$ throughout the paper.

\subsection{Lagrangian Neighbourhood Theorems and regularity}
\label{cu33}

We recall some symplectic geometry, which can be found in McDuff
and Salamon \cite{McSa}. Let $N$ be a real $m$-manifold. Then its
tangent bundle $T^*N$ has a {\it canonical symplectic form}
$\hat\om$, defined as follows. Let $(x_1,\ldots,x_m)$ be local
coordinates on $N$. Extend them to local coordinates
$(x_1,\ldots,x_m,y_1,\ldots,y_m)$ on $T^*N$ such that
$(x_1,\ldots,y_m)$ represents the 1-form $y_1\d x_1+\cdots+y_m
\d x_m$ in $T_{(x_1,\ldots,x_m)}^*N$. Then~$\hat\om=\d x_1\w\d y_1+
\cdots+\d x_m\w\d y_m$.

Identify $N$ with the zero section in $T^*N$. Then $N$ is a
{\it Lagrangian submanifold\/} of $T^*N$. The {\it Lagrangian
Neighbourhood Theorem} \cite[Th.~3.33]{McSa} shows that any
compact Lagrangian submanifold $N$ in a symplectic manifold
looks locally like the zero section in~$T^*N$.

\begin{thm} Let\/ $(M,\om)$ be a symplectic manifold and\/
$N\subset M$ a compact Lagrangian submanifold. Then there
exists an open tubular neighbourhood\/ $U$ of the zero
section $N$ in $T^*N$, and an embedding $\Phi:U\ra M$ with\/
$\Phi\vert_N=\id:N\ra N$ and\/ $\Phi^*(\om)=\hat\om$, where
$\hat\om$ is the canonical symplectic structure on~$T^*N$.
\label{cu3thm1}
\end{thm}

In \cite[\S 4]{Joyc4} we extend Theorem \ref{cu3thm1}
to situations involving conical singularities, first to
{\it SL cones},~\cite[Th.~4.3]{Joyc4}.

\begin{thm} Let\/ $C$ be an SL cone in $\C^m$ with isolated
singularity at\/ $0$, and set\/ $\Si=C\cap{\cal S}^{2m-1}$.
Define $\iota:\Si\t(0,\iy)\ra\C^m$ by $\iota(\si,r)=r\si$,
with image $C\sm\{0\}$. For $\si\in\Si$, $\tau\in T_\si^*\Si$,
$r\in(0,\iy)$ and\/ $u\in\R$, let\/ $(\si,r,\tau,u)$ represent
the point\/ $\tau+u\,\d r$ in $T^*_{\smash{(\si,r)}}\bigl(\Si\!
\t\!(0,\iy)\bigr)$. Identify $\Si\t(0,\iy)$ with the zero section
$\tau\!=\!u\!=\!0$ in $T^*\bigl(\Si\t(0,\iy)\bigr)$. Define an
action of\/ $(0,\iy)$ on $T^*\bigl(\Si\!\t\!(0,\iy)\bigr)$ by
\e
t:(\si,r,\tau,u)\longmapsto (\si,tr,t^2\tau,tu)
\quad\text{for $t\in(0,\iy),$}
\label{cu3eq3}
\e
so that\/ $t^*(\hat\om)\!=\!t^2\hat\om$, for $\hat\om$ the
canonical symplectic structure on~$T^*\bigl(\Si\!\t\!(0,\iy)\bigr)$.

Then there exists an open neighbourhood\/ $U_\sC$ of\/
$\Si\t(0,\iy)$ in $T^*\bigl(\Si\t(0,\iy)\bigr)$ invariant under
\eq{cu3eq3} given by
\begin{equation*}
U_\sC=\bigl\{(\si,r,\tau,u)\in T^*\bigl(\Si\t(0,\iy)\bigr):
\bmd{(\tau,u)}<2\ze r\bigr\}\quad\text{for some $\ze>0,$}
\end{equation*}
where $\md{\,.\,}$ is calculated using the cone metric $\iota^*(g')$
on $\Si\t(0,\iy)$, and an embedding $\Phi_\sC:U_\sC\ra\C^m$
with\/ $\Phi_\sC\vert_{\Si\t(0,\iy)}=\iota$, $\Phi_{\sst
C}^*(\om')=\hat\om$ and\/ $\Phi_\sC\circ t=t\,\Phi_\sC$
for all\/ $t>0$, where $t$ acts on $U_\sC$ as in \eq{cu3eq3}
and on $\C^m$ by multiplication.
\label{cu3thm2}
\end{thm}

In \cite[Th.~4.4]{Joyc4} we construct a particular choice of
$\phi_i$ in Definition~\ref{cu3def3}.

\begin{thm} Let\/ $(M,J,\om,\Om)$, $\psi,X,n,x_i,\up_i,C_i,\Si_i,
\mu_i,R,\Up_i$ and\/ $\iota_i$ be as in Definition \ref{cu3def3}.
Theorem \ref{cu3thm2} gives $\ze>0$, neighbourhoods $U_\sCi$
of\/ $\Si_i\t(0,\iy)$ in $T^*\bigl(\Si_i\t(0,\iy)\bigr)$ and
embeddings $\Phi_\sCi:U_\sCi\ra\C^m$ for~$i=1,\ldots,n$.

Then for sufficiently small\/ $R'\in(0,R]$ there exist unique
closed\/ $1$-forms $\eta_i$ on $\Si_i\t(0,R')$ for $i=1,\ldots,n$
written $\eta_i(\si,r)=\eta_i^1(\si,r)+\eta_i^2(\si,r)\d r$ for
$\eta_i^1(\si,r)\in T_\si^*\Si_i$ and\/ $\eta_i^2(\si,r)\in\R$,
and satisfying $\md{\eta_i(\si,r)}<\ze r$ and\/ $\md{\na^k\eta_i}
=O(r^{\mu_i-1-k})$ as $r\ra 0$ for $k=0,1,$ computing
$\na,\md{\,.\,}$ using the cone metric $\iota_i^*(g')$, such that
the following holds.

Define $\phi_i:\Si_i\t(0,R')\ra B_R$ by $\phi_i(\si,r)=\Phi_{\sst
C_i}\bigl(\si,r,\eta_i^1(\si,r),\eta_i^2(\si,r)\bigr)$. Then
$\Up_i\circ\phi_i:\Si_i\t(0,R')\ra M$ is a diffeomorphism
$\Si_i\t(0,R')\ra S_i$ for open sets $S_1,\ldots,S_n$ in $X'$
with\/ $\bar S_1,\ldots,\bar S_n$ disjoint, and\/ $K=X'\sm(S_1
\cup\cdots\cup S_n)$ is compact. Also $\phi_i$ satisfies
\eq{cu3eq2}, so that\/ $R',\phi_i,S_i,K$ satisfy
Definition~\ref{cu3def3}.
\label{cu3thm3}
\end{thm}

In \cite[\S 5]{Joyc4} we study the asymptotic behaviour of the maps
$\phi_i$ of Theorem \ref{cu3thm3}, using the elliptic regularity of
the special Lagrangian condition. Combining \cite[Th.~5.1]{Joyc4},
\cite[Lem.~4.5]{Joyc4} and \cite[Th.~5.5]{Joyc4} proves:

\begin{thm} In the situation of Theorem \ref{cu3thm3} we have
$\eta_i=\d A_i$ for $i=1,\ldots,n$, where $A_i:\Si_i\t(0,R')\ra\R$
is given by $A_i(\si,r)=\int_0^r\eta_i^2(\si,s)\d s$. Suppose
$\mu_i'\in(2,3)$ with\/ $(2,\mu_i']\cap\D_\sSii=\emptyset$ for
$i=1,\ldots,n$. Then
\e
\begin{gathered}
\bmd{\na^k(\phi_i-\iota_i)}=O(r^{\mu_i'-1-k}),\quad
\bmd{\na^k\eta_i}=O(r^{\mu_i'-1-k})\quad\text{and}\\
\bmd{\na^kA_i}=O(r^{\mu_i'-k})
\quad\text{as $r\ra 0$ for all\/ $k\ge 0$ and\/ $i=1,\ldots,n$.}
\end{gathered}
\label{cu3eq4}
\e

Hence $X$ has conical singularities at $x_i$ with cone $C_i$
and rate $\mu_i'$, for all possible rates $\mu_i'$ allowed by
Definition \ref{cu3def3}. Therefore, the definition of
conical singularities is essentially independent of the
choice of rate~$\mu_i$.
\label{cu3thm4}
\end{thm}

Finally we extend Theorem \ref{cu3thm1} to SL $m$-folds with
conical singularities \cite[Th.~4.6]{Joyc4}, in a way
compatible with Theorems \ref{cu3thm2} and~\ref{cu3thm3}.

\begin{thm} Suppose $(M,J,\om,\Om)$ is an almost Calabi--Yau
$m$-fold and\/ $X$ a compact SL\/ $m$-fold in $M$ with conical
singularities at\/ $x_1,\ldots,x_n$. Let the notation $\psi,\up_i,
C_i,\Si_i,\mu_i,R,\Up_i$ and\/ $\iota_i$ be as in Definition
\ref{cu3def3}, and let\/ $\ze,U_\sCi,\allowbreak
\Phi_\sCi,\allowbreak R',\allowbreak \eta_i,\allowbreak
\eta_i^1,\eta_i^2,\phi_i,S_i$ and\/ $K$ be as in Theorem~\ref{cu3thm3}.

Then making $R'$ smaller if necessary, there exists an open tubular
neighbourhood\/ $U_\sXp\subset T^*X'$ of the zero section
$X'$ in $T^*X'$, such that under $\d(\Up_i\circ\phi_i):T^*\bigl(
\Si_i\t(0,R')\bigr)\ra T^*X'$ for $i=1,\ldots,n$ we have
\e
\bigl(\d(\Up_i\circ\phi_i)\bigr)^*(U_\sXp)=\bigl\{(\si,r,\tau,u)
\in T^*\bigl(\Si_i\t(0,R')\bigr):\bmd{(\tau,u)}<\ze r\bigr\},
\label{cu3eq5}
\e
and there exists an embedding $\Phi_\sXp:U_\sXp\ra M$ with\/
$\Phi_\sXp\vert_{X'}=\id:X'\ra X'$ and\/ $\Phi_\sXp^*(\om)=\hat\om$,
where $\hat\om$ is the canonical symplectic structure on $T^*X'$,
such that
\e
\Phi_\sXp\circ\d(\Up_i\circ\phi_i)(\si,r,\tau,u)\equiv\Up_i\circ
\Phi_\sCi\bigl(\si,r,\tau+\eta_i^1(\si,r),u+\eta_i^2(\si,r)\bigr)
\label{cu3eq6}
\e
for all\/ $i=1,\ldots,n$ and\/ $(\si,r,\tau,u)\in T^*\bigl(\Si_i\t(0,R')
\bigr)$ with\/ $\bmd{(\tau,u)}<\ze r$. Here $\md{(\tau,u)}$ is computed
using the cone metric $\iota_i^*(g')$ on~$\Si_i\t(0,R')$.
\label{cu3thm5}
\end{thm}

\section{Asymptotically Conical SL $m$-folds}
\label{cu4}

Let $C$ be an SL cone in $\C^m$ with an isolated singularity
at 0. Section \ref{cu3} considered SL $m$-folds with conical
singularities, which are asymptotic to $C$ at 0. We now
discuss {\it Asymptotically Conical\/} SL $m$-folds $L$ in
$\C^m$, which are asymptotic to $C$ at infinity. Here is
the definition.

\begin{dfn} Let $C$ be an SL cone in $\C^m$ with isolated
singularity at 0 for $m>2$, as in Definition \ref{cu3def1}, and
let $\Si=C\cap{\cal S}^{2m-1}$, so that $\Si$ is a compact,
nonsingular $(m-1)$-manifold, not necessarily connected.
Let $g_\sSi$ be the metric on $\Si$ induced by the
metric $g'$ on $\C^m$ in \eq{cu2eq1}, and $r$ the radius
function on $\C^m$. Define $\iota:\Si\t(0,\iy)\ra\C^m$ by
$\iota(\si,r)=r\si$. Then the image of $\iota$ is $C\sm\{0\}$,
and $\iota^*(g')=r^2g_\sSi+\d r^2$ is the cone metric
on~$C\sm\{0\}$.

Let $L$ be a closed, nonsingular SL $m$-fold in $\C^m$ and
$\la<2$. We call $L$ {\it Asymptotically Conical (AC)}
with {\it rate} $\la$ and {\it cone} $C$ if there exists a
compact subset $K\subset L$ and a diffeomorphism $\vp:\Si\t
(T,\iy)\ra L\sm K$ for some $T>0$, such that
\e
\bmd{\na^k(\vp-\iota)}=O(r^{\la-1-k})
\quad\text{as $r\ra\iy$ for $k=0,1$.}
\label{cu4eq1}
\e
Here $\na,\md{\,.\,}$ are computed using the cone metric
$\iota^*(g')$ on~$\Si\t(T,\iy)$.
\label{cu4def1}
\end{dfn}

This is very similar to Definition \ref{cu3def3}, and in fact
there are strong parallels between the theories of SL $m$-folds
with conical singularities and of AC SL $m$-folds. We recall
some results from \cite[\S 7]{Joyc4}, including versions of
the material in \S\ref{cu33}. We continue to assume $m>2$
throughout.

\subsection{Cohomological invariants of AC SL $m$-folds}
\label{cu41}

When $Y$ is a manifold, write $H^k(Y,\R)$ for the $k^{\rm th}$
{\it de Rham cohomology group} of $Y$, and $H_k(Y,\R)$ for the
$k^{\rm th}$ {\it real singular homology group} of $Y$, defined
using smooth simplices. Then the pairing between homology and
cohomology is defined at the chain level by integrating $k$-forms
over $k$-simplices. We can also define {\it relative} homology
and cohomology groups in the usual way. The {\it Betti numbers}
of $Y$ are~$b^k(Y)=\dim H^k(Y,\R)$.

Let $L$ be an AC SL $m$-fold in $\C^m$ with cone $C$, and set
$\Si=C\cap{\cal S}^{2m-1}$. As $\Si$ is in effect the boundary
of $L$, there is a natural map $H^k(L,\R)\ra H^k(\Si,\R)$.
Following \cite[Def.~7.2]{Joyc4} we define {\it cohomological
invariants\/} $Y(L),Z(L)$ of~$L$.

\begin{dfn} Let $L$ be an AC SL $m$-fold in $\C^m$ with cone $C$,
and let $\Si=C\cap{\cal S}^{2m-1}$. As $\om',\Im\Om'$ in \eq{cu2eq1}
are closed forms with $\om'\vert_L\equiv\Im\Om'\vert_L=0$, they
define classes in the relative de Rham cohomology groups $H^k(\C^m;
L,\R)$ for $k=2,m$. For $k>1$ we have the exact sequence
\begin{equation*}
0=H^{k-1}(\C^m,\R)\ra H^{k-1}(L,\R){\buildrel\cong\over\longra}
H^k(\C^m;L,\R)\ra H^k(\C^m,\R)=0.
\end{equation*}
Define $Y(L)\in H^1(\Si,\R)$ to be the image of $[\om']$ in
$H^2(\C^m;L,\R)\cong H^1(L,\R)$ under $H^1(L,\R)\ra H^1(\Si,R)$,
and $Z(L)\in H^{m-1}(\Si,\R)$ to be the image of $[\Im\Om']$ in
$H^m(\C^m;L,\R)\cong H^{m-1}(L,\R)$ under~$H^{m-1}(L,\R)\ra
H^{m-1}(\Si,R)$.
\label{cu4def2}
\end{dfn}

Here are some conditions for $Y(L)$ or $Z(L)$ to be
zero,~\cite[Prop.~7.3]{Joyc4}.

\begin{prop} Let\/ $L$ be an AC SL\/ $m$-fold in $\C^m$ with
cone $C$ and rate $\la$, and let\/ $\Si=C\cap{\cal S}^{2m-1}$.
If\/ $\la<0$ or $b^1(L)=0$ then $Y(L)=0$. If\/ $\la<2-m$ or
$b^0(\Si)=1$ then~$Z(L)=0$.
\label{cu4prop}
\end{prop}

In this paper we will consider only AC SL $m$-folds $L_i$ with
rates $\la_i<0$. These all have $Y(L_i)=0$ by the proposition.
Because of this we shall avoid some tricky issues of global
symplectic topology in defining Lagrangian $m$-folds $N^t$ by
gluing $tL_i$ in at a singular point $x_i$ of an SL $m$-fold
$X$ with conical singularities, so \S\ref{cu6} and \S\ref{cu7}
are simplified. The case $Y(L_i)\ne 0$ will be considered
in the sequel \cite{Joyc6}. Here is a (trivial) lemma on
{\it dilations} of AC SL $m$-folds.

\begin{lem} Let\/ $L$ be an AC SL\/ $m$-fold in $\C^m$ with rate
$\la$ and cone $C$, and let\/ $t>0$. Then $tL=\{t\,{\bf x}:{\bf x}
\in L\}$ is also an AC SL\/ $m$-fold in $\C^m$ with rate $\la$
and cone $C$, satisfying $Y(tL)=t^2Y(L)$ and\/~$Z(tL)=t^mZ(L)$.
\label{cu4lem}
\end{lem}

\subsection{Lagrangian Neighbourhood Theorems and regularity}
\label{cu42}

Next we generalize \S\ref{cu33} to AC SL $m$-folds.
Here is the analogue of Theorem \ref{cu3thm3}, proved
in~\cite[Th.~7.4]{Joyc4}.

\begin{thm} Let\/ $C$ be an SL cone in $\C^m$ with isolated
singularity at\/ $0$, and set\/ $\Si=C\cap{\cal S}^{2m-1}$.
Define $\iota:\Si\t(0,\iy)\ra\C^m$ by $\iota(\si,r)=r\si$. Let\/
$\ze$, $U_\sC\subset T^*\bigl(\Si\t(0,\iy)\bigr)$ and\/
$\Phi_\sC:U_\sC\ra\C^m$ be as in Theorem~\ref{cu3thm2}.

Suppose $L$ is an AC SL\/ $m$-fold in $\C^m$ with cone $C$
and rate $\la<2$. Then there exists a compact\/ $K\subset L$
and a diffeomorphism $\vp:\Si\t(T,\iy)\ra L\sm K$ for some
$T>0$ satisfying \eq{cu4eq1}, and a closed\/ $1$-form $\chi$ on
$\Si\t(T,\iy)$ written $\chi(\si,r)=\chi^1(\si,r)+\chi^2(\si,r)\d r$
for $\chi^1(\si,r)\in T_\si^*\Si$ and\/ $\chi^2(\si,r)\in\R$,
satisfying
\e
\begin{gathered}
\bmd{\chi(\si,r)}<\ze r,\quad \vp(\si,r)\equiv
\Phi_\sC\bigl(\si,r,\chi^1(\si,r),\chi^2(\si,r)\bigr)\\
\text{and}\quad\bmd{\na^k\chi}=O(r^{\la-1-k})
\quad\text{as $r\ra\iy$ for $k=0,1,$}
\end{gathered}
\label{cu4eq2}
\e
computing $\na,\md{\,.\,}$ using the cone metric~$\iota^*(g')$.
\label{cu4thm1}
\end{thm}

Now suppose that the rate $\la$ of $L$ satisfies $\la<0$. Then
$Y(L)=0$ by Proposition \ref{cu4prop}, and the results of
\cite[\S 7.3]{Joyc4} simplify. Combining \cite[Prop.~7.6]{Joyc4},
\cite[Th.~7.7]{Joyc4} and \cite[Th.~7.11]{Joyc4} gives an analogue
of Theorem~\ref{cu3thm4}.

\begin{thm} In the situation of Theorem \ref{cu4thm1}, suppose
$\la<0$. Then $\chi=\d E$, where $E\in C^\iy\bigl(\Si\t(T,\iy)
\bigr)$ is given by $E(\si,r)=-\int_r^\iy\chi^2(\si,s)\d s$.
If either $\la=\la'$, or $\la'\in(2-m,0)$, or $\la,\la'$ lie
in the same connected component of\/ $\R\sm\D_\sSi$, then $L$
is an AC SL\/ $m$-fold with rate $\la'$ and
\e
\begin{gathered}
\bmd{\na^k(\vp-\iota)}=O(r^{\la'-1-k}),\quad
\bmd{\na^k\chi}=O(r^{\la'-1-k})\quad\text{and}\\
\bmd{\na^kE}=O(r^{\la'-k})
\quad\text{as $r\ra\iy$ for all\/ $k\ge 0$.}
\end{gathered}
\label{cu4eq3}
\e
Here $\na,\md{\,.\,}$ are computed using the cone
metric $\iota^*(g')$ on~$\Si\t(T,\iy)$.
\label{cu4thm2}
\end{thm}

In particular, this shows that any AC SL $m$-fold $L$ with
rate $\la<0$ is {\it also} an AC SL $m$-fold with rate
$\la'$ for $\la'\in\bigl(2-m,\ha(2-m)\bigr)$. This will be
important in \S\ref{cu62}, where we need to assume that
$\la<\ha(2-m)$ to make an error term sufficiently small.
Here \cite[Th.~7.5]{Joyc4} is the analogue of
Theorem~\ref{cu3thm5}.

\begin{thm} Suppose $L$ is an AC SL\/ $m$-fold in $\C^m$
with cone $C$. Let\/ $\Si,\iota,
\allowbreak
\ze,
\allowbreak
U_\sC,
\allowbreak
\Phi_\sC,
\allowbreak
K,T,\vp,\chi,
\chi^1,\chi^2$ be as in Theorem \ref{cu4thm1}. Then making
$T,K$ larger if necessary, there exists an open tubular
neighbourhood\/ $U_\sL\subset T^*L$ of the zero section $L$ in $T^*L$,
such that under $\d\vp:T^*\bigl(\Si\t(T,\iy)\bigr)\ra T^*L$ we have
\e
(\d\vp)^*(U_\sL)=\bigl\{(\si,r,\tau,u)\in
T^*\bigl(\Si\t(T,\iy)\bigr):\bmd{(\tau,u)}<\ze r\bigr\},
\label{cu4eq4}
\e
and there exists an embedding $\Phi_\sL:U_\sL\ra\C^m$ with\/
$\Phi_\sL\vert_L=\id:L\ra L$ and\/ $\Phi_\sL^*(\om')=\hat\om$,
where $\hat\om$ is the canonical symplectic structure on $T^*L$,
such that
\e
\Phi_\sL\circ\d\vp(\si,r,\tau,u)\equiv
\Phi_\sC\bigl(\si,r,\tau+\chi^1(\si,r),u+\chi^2(\si,r)\bigr)
\label{cu4eq5}
\e
for all\/ $(\si,r,\tau,u)\!\in\!T^*\bigl(\Si\!\t\!(T,\iy)\bigr)$
with\/ $\md{(\tau,u)}<\ze r$, computing $\md{\,.\,}$
using~$\iota^*(g')$.
\label{cu4thm3}
\end{thm}

In \cite[Th.~7.10]{Joyc4} we study the {\it bounded harmonic
functions\/} on~$L$.

\begin{thm} Suppose $L$ is an AC SL\/ $m$-fold in $\C^m$, with cone
$C$. Let\/ $\Si,T$ and\/ $\vp$ be as in Theorem \ref{cu4thm1}. Let\/
$l=b^0(\Si)$, and\/ $\Si^1,\ldots,\Si^l$ be the connected components
of\/ $\Si$. Let\/ $V$ be the vector space of bounded harmonic functions
on $L$. Then $\dim V=l$, and for each\/ ${\bf c}=(c^1,\ldots,c^l)\in
\R^l$ there exists a unique $v^{\bf c}\in V$ such that for all\/
$j=1,\ldots,l$, $k\ge 0$ and\/ $\be\in(2-m,0)$ we have
\begin{equation*}
\na^k\bigl(\vp^*(v^{\bf c})-c^j\,\bigr)=O\bigl(\md{{\bf c}}
r^{\be-k}\bigr) \quad\text{on $\Si^j\t(T,\iy)$ as $r\ra\iy$.}
\end{equation*}
Note also that\/ $V=\{v^{\bf c}:{\bf c}\in\R^l\}$
and\/~$v^{(1,\ldots,1)}\equiv 1$.
\label{cu4thm4}
\end{thm}

\section{An analytic existence result for SL $m$-folds}
\label{cu5}

We shall now use analysis to prove that under certain conditions
a compact, nonsingular Lagrangian $m$-fold $N$ in an almost
Calabi--Yau $m$-fold $M$ which is approximately special Lagrangian
can be deformed to a nearby special Lagrangian $m$-fold $\ti N$ in
$M$. We begin in \S\ref{cu51} with some background material from
analysis. The main result, Theorem \ref{cu5thm2}, is stated in
\S\ref{cu52}, and proved in~\S\ref{cu53}--\S\ref{cu55}.

Theorem \ref{cu5thm2} and its proof are based on results by the
author \cite[Th.~11.6.1 \& Th.~13.6.1]{Joyc2}, which are used
to construct compact 7- and 8-manifolds $M$ with holonomy $G_2$
and Spin(7) by deforming a $G_2$- or Spin(7)-structure with
small torsion on $M$. The geometry is rather different, but the
underlying conception and structure of the proof is the same.

In each of \S\ref{cu6} and \S\ref{cu7} we will construct a
family of compact, nonsingular Lagrangian $m$-folds $N^t$ in
$M$ for $t\in(0,\de)$ by gluing AC SL $m$-folds $L_1,\ldots,L_n$
in at the singular points $x_1,\ldots,x_n$ of a compact SL
$m$-fold $X$ in $M$ with conical singularities. We then apply
Theorem \ref{cu5thm2} to show that $N^t$ can be deformed to
a nearby compact, nonsingular SL $m$-fold $\smash{\ti N^t}$
in $M$ for small~$t$.

The proof of Theorem \ref{cu5thm2} is long and technical, and
some readers may prefer to skip over it. The rest of the paper
will use only the statement of Theorem \ref{cu5thm2}, and not
refer to its proof in~\S\ref{cu53}--\S\ref{cu55}.

\subsection{Banach spaces of functions}
\label{cu51}

Let $(N,g)$ be a Riemannian manifold. To establish notation,
we shall define various Banach spaces of functions on $N$. 
Some references for these spaces are Aubin \cite{Aubi} and
Gilbarg and Trudinger \cite{GiTr}. For each integer $k\ge 0$,
define $C^k(N)$ to be the vector space of continuous, bounded
functions $f$ on $N$ that have $k$ continuous, bounded
derivatives, and define the norm $\cnm{.}k$ on $C^k(N)$ by
$\cnm{f}k=\sum_{j=0}^k\sup_N\bmd{\na^jf}$, where $\na$ is the
Levi-Civita connection. Then $C^k(N)$ is a Banach space.
Let~$C^\iy(N)=\bigcap_{k\ge 0}C^k(N)$.

For $k\ge 0$ and $\al\in(0,1)$, define the {\it H\"older space}
$C^{k,\al}(N)$ to be the subset of $f\in C^k(N)$ for which
\begin{equation*}
[\na^kf]_\al=\sup_{\substack{x\ne y\in N\\ d(x,y)<\de(g)}}
\frac{\bmd{\na^kf(x)-\na^kf(y)}}{d(x,y)^\al}
\end{equation*}
is finite. Here $d(x,y)$ is the geodesic distance between $x$ and $y$
and $\de(g)>0$ the {\it injectivity radius}. Note that $\na^kf(x)$
and $\na^kf(y)$ lie in different vector spaces $\ot^kT_x^*N$,
$\ot^kT_y^*N$ when $k>0$, but we identify them by parallel
translation using $\na$ along the unique geodesic $\ga$ of
length $d(x,y)$ joining $x$ and $y$. The {\it H\"older norm}
is~$\cnm{f}{k,\al}=\cnm{f}{k}+[\na^kf]_\al$.

For $q\ge 1$, define the {\it Lebesgue space} $L^q(N)$ to be the set
of locally integrable functions $f$ on $N$ for which the norm
\begin{equation*}
\lnm{f}q=\left(\int_N\md{f}^q\,\d V_g\right)^{1/q}
\end{equation*}
is finite. Here $\d V_g$ is the volume form of $g$. Suppose
$r,s,t\ge 1$ with $1/r=1/s+1/t$. If $\phi\in L^s(N)$ and
$\psi\in L^t(N)$ then $\phi\psi\in L^r(N)$, and $\lnm{\phi\psi}r
\le\lnm{\phi}s\lnm{\psi}t$; this is {\it H\"older's inequality}.

Let $q\ge 1$ and $k\ge 0$ be an integer. Define the 
{\it Sobolev space} $L^q_k(N)$ to be the set of $f\in L^q(N)$ 
such that $f$ is $k$ times weakly differentiable and 
$\md{\na^jf}\in L^q(N)$ for $j\le k$. Define the {\it Sobolev 
norm} on $L^q_k(N)$ to be
\begin{equation*}
\snm{f}qk=\biggl(\sum_{j=0}^k\int_N\md{\na^jf}^q\,\d V_g\biggr)^{1/q}.
\end{equation*}
Then $L^q_k(N)$ is a Banach space, and $L^2_k(N)$ a Hilbert space.

The {\it Sobolev Embedding Theorem} \cite[Th.~2.30]{Aubi} gives
inclusions between the spaces $L^q_k(N)$ and~$C^{l,\al}(N)$.

\begin{thm} Suppose $(N,g)$ is a compact Riemannian $n$-manifold,
$k\ge l\ge 0$ are integers, $\al\in(0,1)$ and\/ $q,r\ge 1$. If\/
$\frac{1}{q}\le\frac{1}{r}+\frac{k-l}{n}$, then $L^q_k(N)$ is
continuously embedded in $L^r_l(N)$ by inclusion. If\/
$\frac{1}{q}\le\frac{k-l-\al}{n}$, then $L^q_k(N)$ is
continuously embedded in $C^{l,\al}(N)$ by inclusion.
\label{cu5thm1}
\end{thm}

\subsection{Statement of the result}
\label{cu52}

The following definition sets up the notation we shall use.

\begin{dfn} Let $(M,J,\om,\Om)$ be an almost Calabi--Yau $m$-fold,
with metric $g$. Let $N$ be a compact, oriented, immersed,
Lagrangian $m$-submanifold in $M$, with immersion $\iota:N\ra M$,
so that $\iota^*(\om)\equiv 0$. Define $h=\iota^*(g)$, so that
$(N,h)$ is a Riemannian manifold. Let $\d V$ be the volume form
on $N$ induced by the metric $h$ and orientation.

Let $\psi:M\ra(0,\iy)$ be the smooth function given in \eq{cu2eq3}.
Then $\Om\vert_N$ is a complex $m$-form on $N$, and using \eq{cu2eq3}
and the Lagrangian condition we find that $\bmd{\Om\vert_N}=\psi^m$,
calculating $\md{\,.\,}$ using $h$ on $N$. Therefore we may write
\e
\Om\vert_N=\psi^m{\rm e}^{i\th}\,\d V \quad\text{on $N$,}
\label{cu5eq1}
\e
for some phase function ${\rm e}^{i\th}$ on $N$. Suppose that
$\cos\th\ge\ha$ on $N$. Then we can choose $\th$ to be a smooth
function $\th:N\ra(-\frac{\pi}{3},\frac{\pi}{3})$. Suppose that
$[\iota^*(\Im\Om)]=0$ in $H^m(N,\R)$. Then
$\int_N\psi^m\sin\th\,\d V=0$, by~\eq{cu5eq1}.

Suppose we are given a finite-dimensional vector subspace
$W\subset C^\iy(N)$ with $1\in W$. Define $\pi_\sW:L^2(N)\ra W$
to be the projection onto $W$ using the $L^2$-inner product.

For $r>0$, define $\B_r\subset T^*N$ to be the bundle of 1-forms
$\al$ on $N$ with $\md{\al}<r$. Regard $\B_r$ as a noncompact
$2m$-manifold with natural projection $\pi:\B_r\ra N$, whose fibre
at $x\in N$ is the ball of radius $r$ about 0 in $T_x^*N$. We will
sometimes identify $N$ with the zero section of $\B_r$, and
write~$N\subset\B_r$.

At each $y\in\B_r$ with $\pi(y)=x\in N$, the Levi-Civita connection
$\na$ of $h$ on $T^*N$ defines a splitting $T_y\B_r=H\op V$ into
horizontal and vertical subspaces $H,V$, with $H\cong T_xN$ and
$V\cong T_x^*N$. Write $\hat\om$ for the natural symplectic
structure on $\B_r\subset T^*N$, defined using $T\B_r\cong H\op V$
and $H\cong V^*$. Define a natural Riemannian metric $\hat h$ on
$\B_r$ such that the subbundles $H,V$ are orthogonal, and
$\hat h\vert_H=\pi^*(h)$, $\hat h\vert_V=\pi^*(h^{-1})$.

Let $\hat\na$ be the connection on $T\B_r\cong H\op V$ given by
the lift of the Levi-Civita connection $\na$ of $h$ on $N$ in
the horizontal directions $H$, and by partial differentiation
in the vertical directions $V$, which is well-defined as
$T\B_r$ is naturally trivial along each fibre. Then $\hat\na$
preserves $\hat h,\hat\om$ and the splitting $T\B_r\cong H\op V$.
It is {\it not\/} torsion-free in general, but has torsion
$T(\hat\na)$ depending linearly on the Riemann curvature~$R(h)$.

As $N$ is a Lagrangian submanifold of $M$, by Theorem \ref{cu3thm1}
the symplectic manifold $(M,\om)$ is locally isomorphic near $N$ to
$T^*N$ with its canonical symplectic structure. That is, for some
small $r>0$ there exists an immersion $\Phi:\B_r\ra M$ such that
$\Phi^*(\om)=\hat\om$ and $\Phi\vert_N=\iota$. Define an $m$-form
$\be$ on $\B_r$ by~$\be=\Phi^*(\Im\Om)$.

If $\al\in C^\iy(T^*N)$ with $\md{\al}<r$, write $\Ga(\al)$ for the
{\it graph} of $\al$ in $\B_r$. Then $\Phi_*(\Ga(\al))$ is a compact,
immersed submanifold in $M$ diffeomorphic to~$N$.
\label{cu5def1}
\end{dfn}

With this notation, we can state our main result.

\begin{thm} Let\/ $\ka>1$ and\/ $A_1,\ldots,A_8>0$ be real,
and\/ $m\ge 3$ an integer. Then there exist\/ $\ep,K>0$ depending
only on $\ka,A_1,\ldots,A_8$ and\/ $m$ such that the following
holds.

Suppose $0<t\le\ep$ and Definition \ref{cu5def1} holds with\/
$r=A_1t$, and
\begin{itemize}
\item[{\rm(i)}] $\lnm{\psi^m\sin\th}{2m/(m+2)}\le A_2t^{\ka+m/2}$,
$\cnm{\psi^m\sin\th}{0}\le A_2t^{\ka-1}$,
\newline
$\lnm{\d(\psi^m\sin\th)}{2m}\le A_2t^{\ka-3/2}$
and\/~$\lnm{\pi_\sW(\psi^m\sin\th)}{1}\le A_2t^{\ka+m-1}$.
\item[{\rm(ii)}] $\psi\ge A_3$ on $N$.
\item[{\rm(iii)}] $\cnm{\hat\na^k\be}{0}\le A_4t^{-k}$ for
$k=0,1,2$ and\/~$3$.
\item[{\rm(iv)}] The injectivity radius $\de(h)$
satisfies~$\de(h)\ge A_5t$.
\item[{\rm(v)}] The Riemann curvature $R(h)$
satisfies~$\cnm{R(h)}{0}\le A_6t^{-2}$.
\item[{\rm(vi)}] If\/ $v\in L^2_1(N)$ with\/ $\pi_\sW(v)=0$,
then $v\in L^{2m/(m-2)}(N)$ by Theorem \ref{cu5thm1},
and\/~$\lnm{v}{2m/(m-2)}\le A_7\lnm{\d v}{2}$.
\item[{\rm(vii)}] For all\/ $w\in W$ we have
$\lnm{\d^*\d w}{2m/(m+2)}\le\ha A_7^{-1}\lnm{\d w}{2}$.
\newline
For all\/ $w\in W$ with\/ $\int_Nw\,\d V=0$
we have~$\cnm{w}{0}\le A_8t^{1-m/2}\lnm{\d w}{2}$.
\end{itemize}
Here norms are computed using the metric $h$ on $N$ in {\rm(i)},
{\rm(v)}, {\rm(vi)} and\/ {\rm(vii)}, and the metric $\hat h$ on
$\B_{A_1t}$ in {\rm(iii)}. Then there exists $f\in C^\iy(N)$ with\/
$\int_Nf\,\d V=0$, such that\/ $\cnm{\d f}{0}\le Kt^\ka<A_1t$
and\/ $\ti N=\Phi_*\bigl(\Ga(\d f)\bigr)$ is an immersed special
Lagrangian $m$-fold in\/~$(M,J,\om,\Om)$.
\label{cu5thm2}
\end{thm}

The theorem will be proved in \S\ref{cu53}--\S\ref{cu55}. In
the rest of the section we work in the situation of Theorem
\ref{cu5thm2}, so we suppose $M,J,\om,\Om$ and $N$ are given,
we use the notation of Definition \ref{cu5def1}, and we suppose
that $\ka>1,A_1,\ldots,A_8>0$ and $t>0$ are given such that
parts (i)--(vii) of Theorem \ref{cu5thm2} hold.

\subsection{Special Lagrangian submanifolds close to $N$}
\label{cu53}

We begin the proof by studying the conditions for a submanifold
$\ti N$ of $M$ close to $N$ to be special Lagrangian. We write
$\ti N$ as $\Phi\bigl(\Ga(\al)\bigr)$, where $\al$ is a small
1-form on $N$ and $\Ga(\al)$ its graph in~$\B_{A_1t}\subset T^*N$.

\begin{lem} In the situation above, let\/ $\al\in C^\iy(T^*N)$ be
a smooth\/ $1$-form with\/ $\cnm{\al}{0}<A_1t$, and\/ $\Ga(\al)$ the
graph of\/ $\al$ in $\B_{A_1t}$. Then $\ti N=\Phi\bigl(\Ga(\al)\bigr)$
is a special Lagrangian $m$-fold in $M$ if and only if\/ $\d\al=0$
and\/~$\pi_*\bigl(\be\vert_{\Ga(\al)}\bigr)=0$.
\label{cu5lem1}
\end{lem}

\begin{proof} Note that $\pi:\Ga(\al)\ra N$ is a diffeomorphism and
$\Phi:\Ga(\al)\ra M$ an immersion. By Definition \ref{cu2def4},
$\ti N$ is an SL $m$-fold in $M$ if and only if $\om\vert_{\ti N}
\equiv\Im\Om\vert_{\ti N}\equiv 0$. Pulling back by $\Phi$, this
holds if and only if $\hat\om\vert_{\Ga(\al)}\equiv\be
\vert_{\Ga(\al)}\equiv 0$, since $\Phi^*(\om)=\hat\om$
and~$\Phi^*(\Im\Om)=\be$.

Pushing forward by $\pi:\Ga(\al)\ra N$, we see that $\ti N$ is
special Lagrangian if and only if $\pi_*\bigl(\hat\om\vert_{\Ga(\al)}
\bigr)\equiv\pi_*\bigl(\be\vert_{\Ga(\al)}\bigr)\equiv 0$. But as
$\B_{A_1t}\subset T^*N$ and $\hat\om$ is the natural symplectic
structure on $T^*N$ we have $\pi_*\bigl(\hat\om\vert_{\Ga(\al)}
\bigr)=-\d\al$ by a well-known piece of symplectic geometry, and
the lemma follows.
\end{proof}

We rewrite the condition $\pi_*\bigl(\be\vert_{\Ga(\al)}\bigr)=0$
in terms of a function~$F$.

\begin{dfn} Define $\A=\bigl\{\al\in C^\iy(T^*N):
\cnm{\al}{0}<A_1t\bigr\}$, and define $F:\A\ra C^\iy(N)$
by $\pi_*\bigl(\be\vert_{\Ga(\al)}\bigr)=F(\al)\,\d V$.
Then Lemma \ref{cu5lem1} shows that if $\al\in\A$
then $\Phi\bigl(\Ga(\al)\bigr)$ is special Lagrangian
if and only if~$\d\al=F(\al)=0$.

The value of $F(\al)$ at $x\in N$ depends on the tangent
space $T_y\Ga(\al)$, where $y\in\Ga(\al)$ with $\pi(y)=x$.
But $T_y\Ga(\al)$ depends on both $\al\vert_x$ and
$\na\al\vert_x$. Hence $F(\al)$ depends pointwise on both $\al$
and $\na\al$, rather than just $\al$. Therefore we may write
\e
\begin{split}
&F(\al)[x]=F'\bigl(x,\al(x),\na\al(x)\bigr)
\quad\text{for all $x\in N$, where}\\
&F':\bigl\{(x,\ga,\de):x\in N,\;\> \ga\in T_x^*N,\;\>
\md{\ga}<A_1t,\;\> \de\in \ot^2T_x^*N\bigr\}\ra\R
\end{split}
\label{cu5eq2}
\e
is a smooth, nonlinear function. Note that $F$ maps between
infinite-dimensional spaces $\A\ra C^\iy(N)$, but $F'$
maps between finite-dimensional spaces.

For fixed $x\in N$ the variables $\ga,\de$ in the domain
of $F'$ lie in vector spaces $T_x^*N$, $\ot^2T_x^*N$. Thus we
may take partial derivatives in these directions (without using
a connection), with values in the dual spaces $T_xN,\ot^2T_xN$.
Write $\pd_1,\pd_2$ for the partial derivatives in the $\ga,\de$
directions respectively. Then we have
\begin{gather*}
\pd_1F'(x,\ga,\de)\in T_xN,\quad
\pd_2F'(x,\ga,\de)\in \ot^2T_xN,\quad
\pd_1^2F'(x,\ga,\de)\in S^2(T_xN),\\
\pd_1\pd_2F'(x,\ga,\de)\in \ot^3T_xN\quad\text{and}\quad
\pd_2^2F'(x,\ga,\de)\in S^2(\ot^2T_xN).
\end{gather*}
\label{cu5def2}
\end{dfn}

We compute the expansion of $F$ up to first order in~$\al$.

\begin{prop} This function $F$ may be written
\e
F(\al)=\psi^m\sin\th-\d^*\bigl(\psi^m\cos\th\,\al\bigr)+Q(\al),
\label{cu5eq3}
\e
where $Q:\A\ra C^\iy(N)$ is smooth with\/ $Q(\al)=O\bigl(
\ms{\al}\!+\!\ms{\na\al}\bigr)$ for small\/~$\al$.
\label{cu5prop1}
\end{prop}

\begin{proof} It is easy to see that $F$ depends smoothly
on $\al$. Therefore by Taylor's theorem we may expand $F$
about $\al=0$ up to second order, and get an equation with
the general form of \eq{cu5eq3}, with $Q$ smooth. Since
$F(\al)$ depends pointwise on $\al,\na\al$, the second-order
remainder term $Q(\al)$ is of the form $O\bigl(\ms{\al}\!+
\!\ms{\na\al}\bigr)$, and the estimate valid when $\md{\al},
\md{\na\al}$ are small, that is, when $\al$ is small in~$C^1$.

So to prove \eq{cu5eq3} we need to compute $F(0)$ and
$\d F(0)$ and show they coincide with the first two terms
on the right hand side of \eq{cu5eq3}. When $\al=0$ we
have $\Phi_*\bigl(\Ga(0)\bigr)=N$ in $M$, and therefore
$\pi_*\bigl(\be\vert_{\Ga(0)}\bigr)=\Im\Om\vert_N=
\psi^m\sin\th\,\d V$ by \eq{cu5eq1}. Thus
$F(0)=\psi^m\sin\th$ by Definition \ref{cu5def2}, giving
the zeroth order term in~\eq{cu5eq3}.

Next we compute the first order term in $\al$. Let $v$ be
the vector field on $T^*N$ with $v\cdot\hat\om=-\pi^*(\al)$.
Then $v$ is tangent to the fibres of $\pi:T^*N\ra N$, and
$\exp(v)$ maps $T^*N\ra T^*N$ taking $\be\mapsto\al+\be$ for
1-forms $\be$ on $N$. Identifying $N$ with the zero section
of $T^*N$, the image $\exp(v)[N]$ of $N$ under $\exp(v)$ is
$\Ga(\al)\subset\B_{A_1t}\subset T^*N$. More generally,
$\exp(sv)[N]=\Ga(s\al)$ for~$s\in[0,1]$.

Therefore $F(s\al)\,\d V=\exp(sv)^*(\be)$ for $s\in[0,1]$.
Differentiating gives
\e
\begin{split}
\d F(0)[\al]\,\d V&=\frac{\d}{\d s}\bigl(F(s\al)\bigr)\Big\vert_{s=0}
\d V=\frac{\d}{\d s}\bigl(\exp(sv)^*(\be)\bigr)\Big\vert_{s=0}\\
&=\bigl({\cal L}_v\be\bigr)\Big\vert_N
=\bigl(\d(v\cdot\be)+v\cdot(\d\be)\bigr)\Big\vert_N
=\d\bigl((v\cdot\be)\vert_N\bigr),
\end{split}
\label{cu5eq4}
\e
where ${\cal L}_v$ is the Lie derivative, `$\,\cdot\,$'
contracts together vector fields and forms in the usual
way, and we have used the fact that $\d\be=0$ since
$\Om$ is closed and~$\be=\Phi^*(\Im\Om)$.

Fix $x\in N$. We may choose local real coordinates
$(x_1,\ldots,x_m,y_1,\ldots,y_m)$ on $M$ near $x$ such
that at $x$ we have
\begin{align*}
g&=\sum_{j=1}^m(\d x_j^2+\d y_j^2),\quad
T_xN=\Big\langle\frac{\pd}{\pd x_1},\ldots,\frac{\pd}{\pd x_m}
\Big\rangle,\quad
\d V=\d x_1\w\cdots\w\d x_m\big\vert_{T_xN},\\
\om&=\sum_{j=1}^m\d x_j\w\d y_j
\quad\text{and}\quad
\Om=\psi^m{\rm e}^{i\th}(\d x_1+i\d y_1)\w\cdots\w(\d x_m+i\d y_m).
\end{align*}
Identify $T_xM$ with $T_x\B_{A_1t}$ using $\d\Phi$, so that
$\hat\om=\om$ and $\be=\Im\Om$.

Write $\al=\sum_{j=1}^ma_j\,\d x_j$ at $x$ for $a_j\in\R$. Then
$v=\sum_{j=1}^ma_j\frac{\pd}{\pd y_j}$ at $x$ as
$v\cdot\om=-\al$. Calculation with the above expression for
$\Om$ then shows that
\begin{align*}
(v\cdot\be)\vert_N&=
\psi^m\cos\th\sum_{j=1}^m(-1)^{j-1}a_j\,\d x_1\w\cdots\w\d x_{j-1}
\w\d x_{j+1}\w\cdots\w\d x_m\\
&=\psi^m\cos\th\,*_{\sst T_xN}\al
\qquad\qquad\text{at $x$,}
\end{align*}
where $*_{\sst T_xN}$ is the Hodge star on $T_xN$, computed using the
explicit expressions for $g$ and $\d V$ at $x$. Since $*\d V=1$ and
$*\d *=-\d^*$ on 1-forms, equation \eq{cu5eq4} gives
\begin{equation*}
\d F(0)[\al]\,\d V\!=\!\d(\psi^m\cos\th\,*\al)
\!=\!\bigl(*\d*(\psi^m\cos\th\,\al)\bigr)\,\d V
\!=\!\bigl(-\d^*(\psi^m\cos\th\,\al)\bigr)\,\d V.
\end{equation*}
This shows that $\d F(0):\al\mapsto-\d^*(\psi^m\cos\th\,\al)$,
which yields the first order term in \eq{cu5eq3}, and completes
the proof.
\end{proof}

Here are some properties of~$Q$.

\begin{lem} This function $Q$ satisfies $Q(0)=\d Q(0)=0$ and\/
$\int_NQ(\al)\,\d V=0$ for all\/ $\al\in\A$, and\/
$\Phi\bigl(\Ga(\al)\bigr)$ is special Lagrangian if and only if
\e
\d\al=0 \quad\text{and}\quad
\d^*\bigl(\psi^m\cos\th\,\al\bigr)=\psi^m\sin\th+Q(\al).
\label{cu5eq5}
\e
\label{cu5lem2}
\end{lem}

\begin{proof} Proposition \ref{cu5prop1} gives $Q(\al)=
O\bigl(\ms{\al}\!+\!\ms{\na\al}\bigr)$, which implies that
$Q(0)=\d Q(0)=0$. By definition $\pi_*\bigl(\be\vert_{\Ga(\al)}
\bigr)=F(\al)\,\d V$ for $\al\in\A$, so
\begin{equation*}
\int_NF(\al)\,\d V=\int_{\Ga(\al)}\be=\int_{\Ga(0)}\be
=\int_N\iota^*(\Im\Om)=0,
\end{equation*}
as $\be$ is closed, $\Ga(\al)$ and $\Ga(0)$ are homologous,
and $[\iota^*(\Im\Om)]=0$ in $H^m(N,\R)$ by Definition
\ref{cu5def1}. Now $F(0)=\psi^m\sin\th$ by \eq{cu5eq3}, so
$\int_N\psi^m\sin\th\,\d V=0$, and $\int_N\d^*\bigl(\psi^m
\cos\th\,\al\bigr)\,\d V=0$ by integration by parts.
Therefore multiplying \eq{cu5eq3} by $\d V$ and integrating
over $N$ gives $\int_NQ(\al)\,\d V=0$. Finally
$\Phi\bigl(\Ga(\al)\bigr)$ is special Lagrangian if and
only if $\d\al=F(\al)=0$ by Definition \ref{cu5def2},
and by \eq{cu5eq3} this is equivalent to~\eq{cu5eq5}.
\end{proof}

The notation $Q(\al)$ was chosen because $Q$ is {\it approximately
quadratic} for small $\al$. The following estimates of $Q$ are modelled
on the fact that if $q$ is a homogeneous quadratic polynomial on $\R^n$
then $\bmd{q({\bf x})-q({\bf y})}\le C\md{{\bf x}-{\bf y}}\bigl(\md{\bf
x}+\md{\bf y}\bigr)$ for some $C\ge 0$ and all~${\bf x},{\bf y}\in\R^n$.

\begin{prop} There exist\/ $C_1,\ldots,C_4>0$ depending only
on $A_1,A_4,A_6,m$ such that\/ $C_1<A_1$ and if\/ $\al,\be\in\A$
with\/ $\cnm{\al}{0},\cnm{\be}{0}\le C_1t$ and\/ $\cnm{\na\al}{0},
\allowbreak
\cnm{\na\be}{0}\le C_2$ then
\ea
\begin{split}
\bmd{Q(\al)-Q(\be)}\le C_3&\bigl(t^{-1}\md{\al-\be}+
\md{\na\al-\na\be}\bigr)\cdot\\
&\bigl(t^{-1}\md{\al}+t^{-1}\md{\be}+
\md{\na\al}+\md{\na\be}\bigr) \qquad\text{and}
\end{split}
\label{cu5eq6}\\
\begin{split}
\bmd{\d\bigl(Q(\al)-Q(\be)\bigr)}\le C_4&
\Bigl(
t^{-3}\md{\al-\be}\bigl(\md{\al}\!+\!\md{\be}\bigr)\!+\!
t^{-2}\md{\al-\be}\bigl(\md{\na\al}\!+\!\md{\na\be}\bigr)\\
+t^{-1}\md{\al&-\be}\bigl(\md{\na^2\al}\!+\!\md{\na^2\be}\bigr)\!+\!
t^{-2}\md{\na\al-\na\be}\bigl(\md{\al}\!+\!\md{\be}\bigr)\\
+t^{-1}\md{\na\al-&\na\be}\bigl(\md{\na\al}\!+\!\md{\na\be}\bigr)\!+\!
\md{\na\al-\na\be}\bigl(\md{\na^2\al}\!+\!\md{\na^2\be}\bigr)\\
+t^{-1}\md{\na^2\al&-\na^2\be}\bigl(\md{\al}\!+\!\md{\be}\bigr)\!+\!
\md{\na^2\al-\na^2\be}\bigl(\md{\na\al}\!+\!\md{\na\be}\bigr)\Bigr).
\end{split}
\label{cu5eq7}
\ea
\label{cu5prop2}
\end{prop}

\begin{proof} Let $\al,\be\in\A$, fix $x\in N$, and
define a real function $P$ on the triangle $\bigl\{(r,s):0\le
r\le s\le 1\bigr\}$ by $P(r,s)=Q\bigl(r(\al-\be)+s\be\bigr)[x]$.
This is well-defined as $\A$ is convex and contains 0, so
$r(\al-\be)+s\be\in\A$ when $0\le r\le s\le 1$. Then
\begin{equation*}
\bigl(Q(\al)-Q(\be)\bigr)[x]=P(1,1)-P(0,1)
=\int_0^1\frac{\pd P}{\pd r}(u,1)\d u.
\end{equation*}

Lemma \ref{cu5lem2} gives $\d Q(0)=0$, so
$\frac{\pd P}{\pd r}(0,0)=0$. Therefore
\begin{equation*}
\frac{\pd P}{\pd r}(u,1)=\int_0^1\frac{\d}{\d s}\Bigl(
\frac{\pd P}{\pd r}(us,s)\Bigr)\d s=
\int_0^1\Bigl[u\frac{\pd^2P}{\pd r^2}(us,s)
+\frac{\pd^2P}{\pd r\pd s}(us,s)
\Bigr]\d s.
\end{equation*}
Substituting this into the previous equation and changing
variables to $r=us$ and $s$, we obtain
\e
\bigl(Q(\al)-Q(\be)\bigr)[x]=
\int_0^1\int_0^s\Bigl[\frac{r}{s^2}\,\frac{\pd^2P}{\pd r^2}(r,s)
+\frac{1}{s}\,\frac{\pd^2P}{\pd r\pd s}(r,s)\Bigr]\d r\,\d s.
\label{cu5eq8}
\e

By the definitions of $P,Q$ and $F'$ we have
\begin{align*}
P(r,s)=\,&F'\bigl(x,r(\al(x)-\be(x))+s\be(x),
r(\na\al(x)-\na\be(x))+s\na\be(x)\bigr)\\
&-(\psi^m\sin\th)[x]
+r\,\d^*\bigl(\psi^m\cos\th\,(\al-\be)\bigr)[x]
+s\,\d^*\bigl(\psi^m\cos\th\,\be\bigr)[x].
\end{align*}
Taking second derivatives, the last line drops out to give
\begin{align*}
\frac{\pd^2P}{\pd r^2}(r,s)&=
\ot^2(\al-\be)\cdot\pd_1^2F'
+\ot^2(\na\al-\na\be)\cdot\pd_2^2F'\\
&\qquad\qquad
+2(\al-\be)\ot(\na\al-\na\be)\cdot\pd_1\pd_2F'\qquad\text{and}\\
\frac{\pd^2P}{\pd r\pd s}(r,s)&=
(\al-\be)\ot\be\cdot\pd_1^2F'
+(\na\al-\na\be)\ot\na\be\cdot\pd_2^2F'\\
&\qquad\qquad
+\bigl((\al-\be)\ot\na\be+\be\ot(\na\al-\na\be)\bigr)\cdot\pd_1\pd_2F',
\end{align*}
evaluating $\pd_j\pd_kF'$ at $\bigl(x,r(\al\!-\!\be)\!+\!s\be,
r(\na\al\!-\!\na\be)\!+\!s\na\be\bigr)$ and $\al,\be,\na\al,\na\be$
at~$x$.

Substituting these two equations into \eq{cu5eq8} and taking mods gives
\e
\begin{split}
&\bmd{Q(\al)\!-\!Q(\be)}[x]\le
\int_0^1\int_0^s\Bigl[
\bigl(rs^{-2}\ms{\al\!-\!\be}+s^{-1}\md{\al\!-\!\be}\,
\md{\be}\bigr)\bmd{\pd_1^2F'}\\
&+\bigl(2rs^{-2}\md{\al\!-\!\be}\md{\na\al\!-\!\na\be}
\!+\!s^{-1}\md{\al\!-\!\be}\md{\na\be}
\!+\!s^{-1}\md{\be}\md{\na\al\!-\!\na\be}\bigr)\bmd{\pd_1\pd_2F'}\\
&+\bigl(rs^{-2}\ms{\na\al\!-\!\na\be}+s^{-1}\md{\na\al\!-\!\na\be}\,
\md{\na\be}\bigr)\bmd{\pd_2^2F'}\,\,\Bigr]\d r\,\d s.
\end{split}
\label{cu5eq9}
\e
Here $\al,\be,\na\al,\na\be$ are independent of $r,s$ and so
$\md{\al-\be},\ldots,\md{\na\be}$ are constants, but $\pd_j\pd_kF'$
is evaluated at $\bigl(x,r(\al-\be)+s\be,r(\na\al-\na\be)+
s\na\be\bigr)$, so $\bmd{\pd_j\pd_kF'}$ is a function of~$r,s$.

Let us interpret $F'(x,\ga,\de)$ in terms of $\hat\be$. Regard
$(x,\ga)$ as a point in $\B_{A_1t}\subset T^*N$, with $\ga\in T_x^*N$.
Then $T_{(x,\ga)}\B_{A_1t}\cong T_xN\op T_x^*N$ as in Definition
\ref{cu5def1}. Using $\de\in\ot^2T_x^*N$ we define a map
$I_\de:T_xN\ra T_xN\op T_x^*N=T_{(x,\ga)}\B_{A_1t}$ by
$v\mapsto(v,\de\cdot v)$. Then $F'(x,\ga,\de)\,\d V\vert_x$ is
the pullback to $T_xN$ under $I_\de$ of the restriction of
$\hat\be$ to~$T_{(x,\ga)}\B_{A_1t}$.

Because of this, estimates on the derivatives of $\hat\be$
imply estimates on the derivatives of $F'$. In particular,
as $\cnm{\hat\na^k\be}{0}\le A_4t^{-k}$ for $k=0,1,2$ by
part (iii) of Theorem \ref{cu5thm2} we can show that there
exist $C_1,C_2,C>0$ depending only on $A_4,m$ such that
\e
\bmd{\pd_1^2F'}\le Ct^{-2},\quad
\bmd{\pd_1\pd_2F'}\le Ct^{-1}\quad\text{and}\quad
\bmd{\pd_2^2F'}\le C
\quad\text{at $(x,\ga,\de)$},
\label{cu5eq10}
\e
provided $\md{\ga}\le C_1t$ and $\md{\de}\le C_2$. Here the power
of $t$ is determined by the number of derivatives $\pd_1$. This
is because changing $\de$ does not affect the point $(x,\ga)$ in
$\B_{A_1t}$, so $\pd_2$ does not involve differentiating $\hat\be$
on $\B_{A_1t}$. Note that $\hat\na,\pd_1$ are the same in the fibre
directions, both given by partial differentiation.

Substituting \eq{cu5eq10} into \eq{cu5eq9} and integrating we
prove \eq{cu5eq6}, for some $C_3>0$ depending only on $A_4,m$.
Equation \eq{cu5eq7} can be proved by a similar but rather
more complicated argument, which we leave to the reader. The
extra derivative on $Q$ means that we also use the inequalities
$\cnm{R(h)}{0}\le A_6t^{-2}$ and $\cnm{\hat\na^3\be}{0}\le
A_4t^{-3}$ in Theorem~\ref{cu5thm2}.
\end{proof}

\subsection{Some analytic estimates on $N$}
\label{cu54}

Section \ref{cu53} studied the geometry of $M$ near $N$. We now
give some estimates on $N$ itself, Propositions \ref{cu5prop4}
and \ref{cu5prop6} below, depending only on the Riemannian
manifold $(N,h)$. The proofs are based on that of Theorem G1
in the author's book~\cite[\S 11.7]{Joyc2}.

These estimates are all proved by considering small balls in
$N$ of radius $O(t)$, and comparing them with balls of the same
radius in $\R^m$. We begin by showing that the metric $h$ on
balls of radius $O(t)$ in $N$ is close to the Euclidean metric
$g_0$ on $\R^m$ in the $L^{2m}_2$ norm.
                                       
\begin{prop} Let\/ $D_1>0$ be smaller than a positive bound depending
on $m$. Then there exist\/ $D_2,D_3,D_4>0$ depending only on
$A_5,A_6,m$ and\/ $D_1$ such that the following holds. Let\/ $B_2,B_3$
be the balls of radii $2,3$ about\/ $0$ in $\R^m$, and\/ $g_0$ the
Euclidean metric on $B_3$. Set\/ $r=D_2t$. Then for each\/
$x\in N$ we have $D_3t^m\le\vol\bigl(B_r(x)\bigr)\le\vol\bigl(B_{4r}(x)
\bigr)\le D_4t^m$, where $B_r(x)$ is the geodesic ball of radius $r$ 
about\/ $x$, and there is a smooth, injective map $\Psi_x:B_3\ra N$
satisfying $\bsnm{r^{-2}\Psi_x^*(h)-g_0}{2m}{2}\le D_1$
and\/~$B_r(x)\subset\Psi_x(B_2)\subset\Psi_x(B_3)\subset B_{4r}(x)$.
\label{cu5prop3}
\end{prop}

\begin{proof} For simplicity, first suppose that $t=1$. We
require systems of coordinates on open balls in $N$, in which 
the metric $h$ appears close to the Euclidean metric $g_0$ in
the $L^{2m}_2$ norm. These are provided by Jost and Karcher's
theory of {\it harmonic coordinates} \cite{JoKa}. Jost and
Karcher show that if the injectivity radius is bounded below
and the sectional curvature is bounded above, then there exist
coordinate systems on all balls of a given radius, in which
the $C^{1,\al}$ norm of the metric is bounded in terms of
$\al$ for each~$\al\in(0,1)$.

The $C^{1,\al}$ norm is not quite strong enough for our
purposes, but fortunately Jost and Karcher's results can be
improved to the $L^p_2$ norm, for $p>m/2$. This was done mainly by
Anderson, and is described in Petersen \cite[\S 4--\S 5]{Pete}.
From \cite[Th.~5.1, p.~185]{Pete} we deduce that since
$\de(h)\ge A_5$ and $\cnm{R(h)}0\le A_6$ (as $t=1$),
for $D_2>0$ depending only on $A_5,A_6,m$ and $D_1$, there
exists a coordinate system $\Psi_x$ about $x$ for each
$x\in N$, which we may write as a map $\Psi_x:B_3\ra N$ with
$\Psi_x(0)=x$, such that $\snm{D_2^{-2}\Psi_x^*(h)-g_0}{2m}{2}
\le D_1$, as we have to prove.

Now the radius and volume of balls are controlled by the $C^0$
norm of the metric on the balls, which is controlled by the
$L^{2m}_2$ norm by Theorem \ref{cu5thm1}. Thus if $D_1$ is
small enough in terms of $m$, the balls $\Psi_x(B_2)$,
$\Psi_x(B_3)$ in $N$ must have volume and radius close to
those of the balls of radius $2D_2$ and $3D_2$ in~$\R^m$.

By making $D_1$ and $D_2$ smaller if necessary, we can ensure 
that $D_3\le\vol(B_{D_2}(x))$ and $\vol(B_{4D_2}(x))\le D_4$
for some $D_3,D_4>0$ depending only on $A_5,A_6,m$ and $D_1$,
and that $B_{D_1}(x)\subset\Psi_x(B_2)$ and $\Psi_x(B_3)\subset 
B_{4D_1}(x)$, for all $x\in N$. This completes the proof when
$t=1$. To prove the proposition for general $t>0$, apply the
case $t=1$ to the rescaled metric~$t^{-2}h$.
\end{proof}

By the Sobolev Embedding Theorem, Theorem \ref{cu5thm1},
$L^{2m}_1$ embeds in $C^0$. Using this we may prove the
following result on balls in $\R^m$, following
\cite[Lem.~2.22]{Aubi}. It is easy to modify the proof
to get a bound involving $\lnm{u}2$ rather than~$\lnm{u}{2m}$.

\begin{lem} Let\/ $B_2,B_3$ be the balls of radii $2,3$
about\/ $0$ in $\R^m$. Then there exist $D_5,D_6>0$
depending only on $m$ such that if\/ $u\in C^1(B_3)$
and\/ $v\in L^{2m}_1(B_3)$ then $\cnm{u\vert_{B_2}}0\le D_5
\bigl(\cnm{\d u}{0}+\lnm{u}2\bigr)$ and\/~$\cnm{v\vert_{B_2}}0
\le D_6\bigl(\lnm{\d v}{2m}+\lnm{v}2\bigr)$.
\label{cu5lem3}
\end{lem}

We can now prove a {\it Sobolev embedding result\/} for 1-forms on~$N$.

\begin{prop} There exist\/ $C_5,C_6>0$ depending only on $A_5,A_6$ and\/
$m$ such that if\/ $\al\in L^{2m}_2(T^*N)$ then $\al\in C^1(T^*N)$ and
\ea
\cnm{\al}{0}&\le C_5\bigl(t\cnm{\na\al}{0}
+t^{-m/2}\lnm{\al}{2}\bigr),\quad\text{and}
\label{cu5eq11}\\
\cnm{\na\al}{0}&\le C_6\bigl(t^{1/2}\lnm{\na^2\al}{2m}
+t^{-m/2}\lnm{\na\al}{2}\bigr).
\label{cu5eq12}
\ea
\label{cu5prop4}
\end{prop}

\begin{proof} Let $D_1>0$ be sufficiently small in terms of
$m$, and let $x\in N$. Then Proposition \ref{cu5prop3} gives
$D_2,D_3,D_4>0$ and $\Psi_x:B_3\ra N$. Define $u\in L^{2m}_2(B_3)$
and $v\in L^{2m}_1(B_3)$ by $u=\Psi_x^*\bigl(\md{\al}_h\bigr)$
and $v=\Psi_x^*\bigl(\md{\na\al}_h\bigr)$, where $\md{\,.\,}_h$
is taken using the metric~$h$.

Lemma \ref{cu5lem3} then applies to $u$ and $v$, as $u\in
C^1(B_3)$ by Theorem \ref{cu5thm1}. The norms in Lemma
\ref{cu5lem3} are calculated w.r.t.\ $g_0$ on $B_3$.
But $\snm{r^{-2}\Psi_x^*(h)-g_0}{2m}{2}\le D_1$ by
Proposition \ref{cu5prop3}, so if $D_1$ is small then
the metrics $r^{-2}\Psi_x^*(h)$ and $g_0$ are close in
$C^0$. Hence we can increase $D_5,D_6$ to $D_5',D_6'$
depending only on $D_1,m$ such that $\cnm{u\vert_{B_2}}0
\le D_5'\bigl(\cnm{\d u}{0}+\lnm{u}2\bigr)$ and
$\cnm{v\vert_{B_2}}0\le D_6'\bigl(\lnm{\d v}{2m}+\lnm{v}2\bigr)$,
where now all norms are taken w.r.t.\ the metric
$r^{-2}\Psi_x^*(h)$ on $B_3$.

Pushing these forward via $\Psi_x$ we deduce that
\begin{align*}
\bcnm{\md{\al}\big\vert_{\Psi_x(B_2)}}{0}&\le
D_5'\bigl(r\bcnm{\d\md{\al}\big\vert_{\Psi_x(B_3)}}{0}
+r^{-m/2}\blnm{\md{\al}\big\vert_{\Psi_x(B_3)}}{2}\bigr)
\;\>\text{and}\\
\bcnm{\md{\na\al}\big\vert_{\Psi_x(B_2)}}{0}&\le
D_6'\bigl(r^{1/2}\blnm{\,\d\md{\na\al}\big\vert_{\Psi_x(B_3)}}{2}
+r^{-m/2}\blnm{\md{\al}\big\vert_{\Psi_x(B_3)}}{2}\bigr),
\end{align*}
where all mods and norms are taken w.r.t.\ $h$, and the powers of
$r$ compensate for the change from $r^{-2}h$ to $h$. Substituting
$r=D_2t$, noting that $\bmd{\d\md{\al}}\le\md{\na\al}$ and
$\bmd{\d\md{\na\al}}\le\md{\na^2\al}$, and taking the
supremum of these two inequalities over all $x\in N$, we
quickly prove \eq{cu5eq11} and \eq{cu5eq12} with $C_5=D_5'
\max(D_2,D_2^{-m/2})$ and~$C_6=D_6'\max(D_2^{1/2},D_2^{-m/2})$.
\end{proof}

Next we prove some interior elliptic regularity estimates on~$B_2,B_3$.

\begin{prop} Suppose $E_1,E_2>0$. Then there exist\/ $E_3,E_4>0$
depending only on $E_1,E_2$ and\/ $m$ such that the following holds.

Let\/ $B_2,B_3$ be the balls of radii $2,3$ about\/ $0$ in $\R^m$,
and suppose $a^{ij}\in L^{2m}_2(B_3)$ for $i,j=1,\ldots,m$
and\/ $b^i\in L^{2m}_1(B_3)$ for $i=1,\ldots,m$ such that\/
\e
\begin{gathered}
E_1\sum_{i=1}^m\xi_i^2\le -\sum_{i,j=1}^ma^{ij}\xi_i\xi_j
\quad\text{on $B_3$ for all\/ $(\xi_1,\ldots,\xi_m)\in\R^m$,}\\
\text{and\/}\quad
\snm{a^{ij}}{2m}{2},\snm{b^i}{2m}{1}\le E_2
\quad\text{for all\/ $i,j=1,\ldots,m$.}
\end{gathered}
\label{cu5eq13}
\e
Then whenever $\si\in L^{2m}_3(B_3)$ and $\tau\in L^{2m}_1(B_3)$ with
\e
\sum_{i,j=1}^ma^{ij}\frac{\pd^2\si}{\pd x_i\pd x_j}
+\sum_{i=1}^mb^i\frac{\pd\si}{\pd x_i}=\tau,
\label{cu5eq14}
\e
we have
\ea
\blnm{\na^2\si\vert_{B_2}}{2}&\le E_3\bigl(\lnm{\d\si}{2}
+\lnm{\tau}{2}\bigr),\quad\text{and}
\label{cu5eq15}\\
\blnm{\na^3\si\vert_{B_2}}{2m}&\le E_4\bigl(\lnm{\d\si}{2}
+\snm{\tau}{2m}{1}\bigr).
\label{cu5eq16}
\ea
\label{cu5prop5}
\end{prop}

\begin{proof} Aubin \cite[Cor.~4.3]{Aubi} shows that
if $(M,g)$ is a compact Riemannian manifold and
$\vp\in L^2_1(M)$ with $\int_M\vp\,\d V=0$ then
$\lnm{\vp}{2}\le C\lnm{\d\vp}{2}$ for $C>0$ depending
on $(M,g)$. Using the technique of `doubling' we see
this also holds for compact Riemannian manifolds
with boundary. Therefore if $\vp\in L^2_1(B_3)$ with
$\int_{B_3}\vp\,\d V_{g_0}=0$ then $\lnm{\vp}{2}\le
E_5\lnm{\d\vp}{2}$ for $E_5>0$ depending only on $m$.
Since \eq{cu5eq14}--\eq{cu5eq16} are unchanged by adding
a constant to $\si$, we may assume that $\int_{B_3}\si\,\d
V_{g_0}=0$, and thus we have~$\lnm{\si}{2}\le E_5\lnm{\d\si}{2}$.

Equation \eq{cu5eq14} is a {\it second-order linear elliptic
equation}, with {\it coefficients} $a^{ij},b^i$. As
$L^{2m}_1\hookra C^{0,1/2}$ by Theorem \ref{cu5thm1}, the bounds
\eq{cu5eq13} imply a $C^{0,1/2}$ bound on the coefficients of
\eq{cu5eq14}, and also show that \eq{cu5eq14} is {\it
uniformly elliptic}. By the interior elliptic regularity
estimates of Gilbarg and Trudinger \cite[Th.~9.11, p.~235]{GiTr}
there exists $E_6>0$ depending only on $E_1,E_2$ and $m$ such
that $\snm{\si\vert_{B_2}}{2}{2}\le E_6\bigl(\lnm{\si}{2}+
\lnm{\tau}{2}\bigr)$. Combining this with
$\lnm{\na^2\si\vert_{B_2}}{2}\le\snm{\si\vert_{B_2}}{2}{2}$
and $\lnm{\si}{2}\le E_5\lnm{\d\si}{2}$ we deduce \eq{cu5eq15},
with~$E_3=E_6\max(E_5,1)$.

We would like to apply the same approach to prove \eq{cu5eq16}.
However, there is a problem: to get an $L^{2m}_3$ estimate of
$\si$ we need a $C^{1,1/2}$ bound on the coefficients of \eq{cu5eq14},
which we do not have. So we instead rewrite \eq{cu5eq14} as
\e
\sum_{i,j=1}^ma^{ij}\frac{\pd^2\si}{\pd x_i\pd x_j}=\tau-
\sum_{i=1}^mb^i\frac{\pd\si}{\pd x_i}=\ti\tau.
\label{cu5eq17}
\e
Now the l.h.s.\ is a uniformly elliptic operator with coefficients
$a^{ij}$, which are bounded in $C^{1,1/2}$ by $\snm{a^{ij}}{2m}{2}
\le E_2$ and the Sobolev embedding~$L^{2m}_2\hookra C^{1,1/2}$.

Let $B_{5/2}$ be the ball of radius $5/2$ about 0 in $\R^m$.
Applying \cite[Th.~9.19, p.~243]{GiTr} to \eq{cu5eq17} on $B_{5/2}$
we get a bound of the form
\e
\bsnm{\si\vert_{B_2}}{2m}{3}\le E_7\bigl(
\blnm{\si\vert_{B_{5/2}}}{2}+\snm{\ti\tau\vert_{B_{5/2}}}{2m}{1}\bigr),
\label{cu5eq18}
\e
for $E_7>0$ depending on $E_1,E_2$ and $m$. Now as
$L^{2m}_1\hookra C^0$ by Theorem \ref{cu5thm1} we can show
that multiplication is a {\it continuous} map
$L^{2m}_1(B_{5/2})\t L^{2m}_1(B_{5/2})\ra L^{2m}_1(B_{5/2})$,
and so there exists $C''>0$ depending only on $m$ with
\e
\snm{\vp\,\eta}{2m}{1}\le C''\snm{\vp}{2m}{1}\snm{\eta}{2m}{1}
\quad\text{for all $\vp,\eta\in L^{2m}_1(B_{5/2})$.}
\label{cu5eq19}
\e

But $\snm{b^i}{2m}{1}\le E_2$, and by the method used to prove
\eq{cu5eq15} we can bound $\snm{\si\vert_{B_{5/2}}}{2m}{2}$,
and hence $\snm{\frac{\pd\si}{\pd x_j}\vert_{B_{5/2}}}{2m}{1}$,
in terms of $\lnm{\si}{2}$ and $\lnm{\tau}{2m}$. Therefore using
\eq{cu5eq17} and \eq{cu5eq19} we can bound $\snm{\ti\tau
\vert_{B_{5/2}}}{2m}{1}$ in terms of $\lnm{\si}{2}$ and
$\snm{\tau}{2m}{1}$. Combining this with \eq{cu5eq18} and
$\lnm{\si}{2}\le E_5\lnm{\d\si}{2}$, we prove \eq{cu5eq16}
for some $E_4>0$ depending only on $E_1,E_2$ and~$m$.
\end{proof}

Finally we prove an elliptic regularity result for
$u\mapsto\d^*\bigl(\psi^m\cos\th\,\d u\bigr)$ on~$N$.

\begin{prop} There exist\/ $C_7,C_8>0$ depending only on
$A_3,A_4,A_5,A_6$ and\/ $m$ such that if\/ $u\in L^{2m}_3(N)$,
$v\in L^{2m}_1(N)$ with\/ $\d^*\bigl(\psi^m\cos\th\,\d u\bigr)=v$ then
\ea
\lnm{\na^2u}{2}&\le C_7\bigl(t^{-1}\lnm{\d u}{2}
+\lnm{v}{2}\bigr),\quad\text{and}
\label{cu5eq20}\\
\lnm{\na^3u}{2m}&\le C_8\bigl(t^{-(m+3)/2}\lnm{\d u}{2}
+t^{-1}\lnm{v}{2m}+\lnm{\d v}{2m}\bigr).
\label{cu5eq21}
\ea
\label{cu5prop6}
\end{prop}

\begin{proof} Let $x\in N$, and define $\si=\Psi_x^*(u)\in
L^{2m}_3(B_3)$ and $\tau=r^2\Psi_x^*(v)\in L^{2m}_1(B_3)$. Then
pulling the equation $\d^*\bigl(\psi^m\cos\th\,\d u\bigr)=v$
back using $\Psi_x$ and rewriting it in the standard coordinates
$(x_1,\ldots,x_m)$ on $B_3$, calculation shows that \eq{cu5eq14}
holds, where
\e
\begin{split}
a^{ij}&=-\Psi_x^*(\psi^m\cos\th)\bigl[\bigl(r^{-2}
\Psi_x^*(h)\bigr)^{-1}\bigr]^{ij}
\quad\text{and}\\
b^i&=\sum_{j=1}^m\Bigl(\frac{\pd a^{ij}}{\pd x_j}+\frac{a^{ij}}{2}
\frac{\pd}{\pd x_j}\bigl(
\log\det\bigl[r^{-2}\Psi_x^*(h)\bigr]_{kl}\bigr)\Bigr).
\end{split}
\label{cu5eq22}
\e
Here $(r^{-2}\Psi_x^*(h))^{-1}$ denotes the inverse in $S^2TB_3$
of the metric $r^{-2}\Psi_x^*(h)$, and $[\ldots]^{ij},$
$[\ldots]_{kl}$ is the index notation for tensors, and
$[r^{-2}\Psi_x^*(h)]_{kl}$ is regarded as an $m\t m$ matrix,
so that we can take its determinant.

Now by Proposition \ref{cu5prop3}, if $D_1$ is small then
$r^{-2}\Psi_x^*(h)$ is $L^{2m}_2$ close to $g_0$, and
hence $C^1$ close as $L^{2m}_2\hookra C^1$ by Theorem
\ref{cu5thm1}. Therefore $(r^{-2}\Psi_x^*(h))^{-1}$ is
$L^{2m}_2$ and $C^1$ close to $g_0^{-1}$, and thus
$[(r^{-2}\Psi_x^*(h))^{-1}]^{ij}$ is $L^{2m}_2$ and $C^1$
close to $\de^i_j$.

Since the component of $\be\vert_N$ in $\La^{m-1}H^*\ot V^*$ is
$\psi^m\cos\th$, the inequality $\cnm{\be}{0}\le A_4$ in
Theorem \ref{cu5thm2} implies that $\psi^m\cos\th\le A_4$.
Using $\cos\th\ge\ha$ from Definition \ref{cu5def1} and
$\psi\ge A_3>0$ by Theorem \ref{cu5thm2} then gives
$0<\ha A_3^m\le\Psi_x^*(\psi^m\cos\th)\le A_4$ on $B_3$.
Combining this with the fact that $[(r^{-2}\Psi_x^*(h))^{-1}]^{ij}$
is $C^0$ close to $\de^i_j$, we can find $E_1>0$ depending only on
$D_1,A_3,A_4$ and $m$ such that the first equation of \eq{cu5eq13} holds.

Again, the inequality $\cnm{\hat\na^k\be}{0}\le A_4t^{-k}$
for $k=0,1,2$ in Theorem \ref{cu5thm2} implies that
$\bmd{\na^k(\psi^m\cos\th)\vert_N}_h\le A_4t^{-k}$ for $k=0,1,2$ on
$N$. As $r=D_2t$, this implies that $\bmd{\na^k(\psi^m\cos\th)
\vert_N}_{r^{-2}h}\le A_4D_2^k$ for $k=0,1,2$ on $N$, taking
$\md{\,.\,}$ using the metric $r^{-2}h$ rather than $h$. Pulling
back to $B_3$ using $\Psi_x$ and using the fact that $r^{-2}
\Psi_x^*(h)$ is $C^1$ close to $g_0$, we can find $E_8>0$ depending
on $A_4,D_1,D_2$ and $m$ such that $\bmd{\pd^k\Psi_x^*(\psi^m\cos
\th)}_{g_0}\le E_8$ on $B_3$ for~$k=0,1,2$.

Combining this with \eq{cu5eq22} and $\snm{r^{-2}\Psi_x^*(h)-
g_0}{2m}{2}\le D_1$ we can find $E_2>0$ depending only on
$A_3,A_4,D_1,D_2$ and $m$ such that the second equation of
\eq{cu5eq13} holds. Therefore Proposition \ref{cu5prop5} gives
$E_3,E_4>0$ such that \eq{cu5eq15} and \eq{cu5eq16} hold. Here
$\na$ and all norms are taken w.r.t.\ $g_0$ on $B_3$. But as
$\snm{r^{-2}\Psi_x^*(h)-g_0}{2m}{2}\le D_1$ we can increase
$E_3,E_4$ to $E_3',E_4'$ depending only on $E_3,E_4,D_1$ and
$m$ such that \eq{cu5eq15} and \eq{cu5eq16} hold with $\na$
and all norms taken w.r.t.\ $r^{-2}\Psi_x^*(h)$ on~$B_3$.

Pushing these inequalities forward with $\Psi_x$ and remembering
that $\si=\Psi_x^*(u)$ and $\tau=r^2\Psi_x^*(v)$ gives
\ea
\blnm{\na^2u\vert_{\Psi_x(B_2)}}{2}&\le E_3'\bigl(
r^{-1}\blnm{\,\d u\vert_{\Psi_x(B_3)}}{2}
+\blnm{v\vert_{\Psi_x(B_3)}}{2}\bigr)\quad\text{and}
\label{cu5eq23}\\
\begin{split}
\blnm{\na^3u\vert_{\Psi_x(B_2)}}{2m}&\le E_4'\bigl(
r^{-(m+3)/2}\blnm{\,\d u\vert_{\Psi_x(B_3)}}{2}\\
&\qquad\quad +r^{-1}\blnm{v\vert_{\Psi_x(B_3)}}{2m}
+\blnm{\,\d v\vert_{\Psi_x(B_3)}}{2m}\bigr).
\end{split}
\label{cu5eq24}
\ea
Here $\na$ and all norms are taken w.r.t.\ $h$, and the powers
of $r$ compensate for the change of metrics from $r^{-2}h$ to
$h$, and the $r^2$ factor in~$\tau=r^2\Psi_x^*(v)$.

Raising \eq{cu5eq23} and \eq{cu5eq24} to the powers $2$ and $2m$
respectively, we deduce
\begin{align*}
\int_{\Psi_x(B_2)}\ms{\na^2u}\d V
&\le 2(E_3')^2\int_{\Psi_x(B_3)}
\bigl(r^{-2}\ms{\d u}+\ms{v}\bigr)\d V\quad\text{and}\\
\int_{\Psi_x(B_2)}\md{\na^3u}^{2m}
&\le 3^{2m-1}(E_4')^{2m}
r^{-m(m+3)}\blnm{\,\d u\vert_{\Psi_x(B_3)}}{2}^{2m-2}
\int_{\Psi_x(B_3)}\ms{\d u}\d V\\
&+3^{2m-1}(E_4')^{2m}
\int_{\Psi_x(B_3)}\bigl(r^{-2m}\md{v}^{2m}+\md{\d v}^{2m}\bigr)\d V,
\end{align*}
since $(a+b)^2\le 2(a^2+b^2)$ and $(a+b+c)^{2m}\le 3^{2m-1}
(a^{2m}+b^{2m}+c^{2m})$. As $B_r(x)\subset\Psi_x(B_2)\subset
\Psi_x(B_3)\subset B_{4r}(x)$ by Proposition \ref{cu5prop3} and
$\lnm{\d u\vert_{\Psi_x(B_3)}}{2}\le\lnm{\d u}{2}$, this gives
\ea
\int_{B_r(x)}\ms{\na^2u}\d V
&\le 2(E_3')^2\int_{B_{4r}(x)}
\bigl(r^{-2}\ms{\d u}+\ms{v}\bigr)\d V\quad\text{and}
\label{cu5eq25}\\
\begin{split}
\int_{B_r(x)}\md{\na^3u}^{2m}
&\le 3^{2m-1}(E_4')^{2m}r^{-m(m+3)}\lnm{\d u}{2}^{2m-2}
\int_{B_{4r}(x)}\ms{\d u}\d V\\
&+3^{2m-1}(E_4')^{2m}
\int_{B_{4r}(x)}\bigl(r^{-2m}\md{v}^{2m}+\md{\d v}^{2m}\bigr)\d V.
\end{split}
\label{cu5eq26}
\ea

Integrate \eq{cu5eq25} and \eq{cu5eq26} over $x\in N$. The left
hand side of \eq{cu5eq25} gives
\begin{align*}
\int_{x\in N}\int_{y\in B_r(x)}\ms{\na^2u}(y)\d V_y\d V_x&=
\int_{y\in N}\int_{x\in B_r(y)}\ms{\na^2u}(y)\d V_x\d V_y\\
&=\int_{y\in N}\vol\bigl(B_r(y)\bigr)\ms{\na^2u}(y)\d V_y,
\end{align*}
exchanging the order of integration of $x,y$ and noting that
$y\in B_r(x)$ if and only if $x\in B_r(y)$. Using
$D_3t^m\le\vol\bigl(B_r(y)\bigr)\le\vol\bigl(B_{4r}(y)
\bigr)\le D_4t^m$ from Proposition \ref{cu5prop3}, we get
\begin{align*}
\int_ND_3t^m\ms{\na^2u}\d V&\le 2(E_3')^2\int_ND_4t^m
\bigl(r^{-2}\ms{\d u}+\ms{v}\bigr)\d V\quad\text{and}\\
\int_ND_3t^m\md{\na^3u}^{2m}
&\le 3^{2m-1}(E_4')^{2m}r^{-m(m+3)}\lnm{\d u}{2}^{2m-2}
\int_ND_4t^m\ms{\d u}\d V\\
&+3^{2m-1}(E_4')^{2m}
\int_ND_4t^m\bigl(r^{-2m}\md{v}^{2m}+\md{\d v}^{2m}\bigr)\d V,
\end{align*}
or equivalently, dividing by $D_3t^m$ and substituting $r=D_2t$,
\ea
\lnm{\na^2u}{2}^2&\le 2D_3^{-1}D_4(E_3')^2
\bigl(D_2^{-2}t^{-2}\lnm{\d u}{2}^2+\lnm{v}{2}^2\bigr)\quad\text{and}
\label{cu5eq27}\\
\begin{split}
\lnm{\na^3u}{2m}^{2m}&\le
3^{2m-1}D_3^{-1}D_4(E_4')^{2m}
D_2^{-m(m+3)}t^{-m(m+3)}\lnm{\d u}{2}^{2m}\\
&+3^{2m-1}D_3^{-1}D_4(E_4')^{2m}
\bigl(D_2^{-2m}t^{-2m}\lnm{v}{2m}^{2m}+\lnm{\d v}{2m}^{2m}\bigr).
\end{split}
\label{cu5eq28}
\ea

Raising \eq{cu5eq27}, \eq{cu5eq28} to the powers $\ha$, $\frac{1}{2m}$
and using $(a\!+\!b)^{1/2}\le a^{1/2}\!+\!b^{1/2}$ and
$(a\!+\!b\!+\!c)^{1/2m}\le a^{1/2m}\!+\!b^{1/2m}\!+\!c^{1/2m}$ for
$a,b,c\ge 0$ yields \eq{cu5eq20} and \eq{cu5eq21} with
\begin{align*}
C_7&=2^{1/2}D_3^{-1/2}D_4^{1/2}E_3'\max(D_2^{-1},1)
\quad\text{and}\\
C_8&=3^{1-1/2m}D_3^{-1/2m}D_4^{1/2m}E_4'\max(D_2^{-(m+3)/2},D_2^{-1},1).
\end{align*}
These depend only on $A_3,A_4,A_5,A_6$ and $m$, and the proof is complete.
\end{proof}

\subsection{The proof of Theorem \ref{cu5thm2}}
\label{cu55}

We shall now prove Theorem \ref{cu5thm2}. Let $f\in C^\iy(N)$,
and apply Lemma \ref{cu5lem2} to the 1-form $\al=\d f$. As $\d\al=0$
automatically, the lemma shows that if $\cnm{\d f}{0}<A_1t$ then
$\Phi\bigl(\Ga(\d f)\bigr)$ is special Lagrangian if and only if
\e
\d^*\bigl(\psi^m\cos\th\,\d f\bigr)=\psi^m\sin\th+Q(\d f).
\label{cu5eq29}
\e
The idea of the proof is to construct by induction a sequence
$(f_n)_{n=0}^\iy$ in $C^\iy(N)$ satisfying $\cnm{\d f_n}{0}<A_1t$
and
\e
\d^*\bigl(\psi^m\cos\th\,\d f_n\bigr)=\psi^m\sin\th+Q(\d f_{n-1})
\label{cu5eq30}
\e
for $n\ge 1$. Then we prove using a priori estimates that the $f_n$
converge in $L^{2m}_3(N)$ to $f$ which satisfies \eq{cu5eq29}, and
finally we show that $f\in C^\iy(N)$ by elliptic regularity.

We start with two lemmas. The first follows from Aubin
\cite[Th.~4.7]{Aubi}, and will give existence for $f_n$ in
\eq{cu5eq30} by induction.

\begin{lem} For each\/ $v\in C^\iy(N)$ with\/ $\int_Nv\,\d V=0$
there exists a unique $u\in C^\iy(N)$ with\/ $\int_Nu\,\d V=0$
and\/~$\d^*\bigl(\psi^m\cos\th\,\d u\bigr)=v$.
\label{cu5lem4}
\end{lem}

\begin{lem} Let\/ $v\in C^\iy(N)$ with\/ $\pi_\sW(v)=0$ and\/
$w\in W$. Then $\lnm{\d v}{2}+\lnm{\d w}{2}\le 2\lnm{\d v+\d w}{2}$.
\label{cu5lem5}
\end{lem}

\begin{proof} Using integration by parts and H\"older's inequality
gives
\begin{align*}
\lnm{\d v+\d w}{2}^2
&=\lnm{\d v}{2}^2\!+\!\lnm{\d w}{2}^2\!+\!2\an{\d v,\d w}
\!=\!\lnm{\d v}{2}^2\!+\!\lnm{\d w}{2}^2\!+\!2\an{v,\d^*\d w}\\
&\ge\lnm{\d v}{2}^2+\lnm{\d w}{2}^2-2\lnm{v}{2m/(m-2)}
\lnm{\d^*\d w}{2m/(m+2)}\\
&\ge\lnm{\d v}{2}^2+\lnm{\d w}{2}^2
-2A_7\lnm{\d v}{2}\cdot \ha A_7^{-1}\lnm{\d w}{2}\\
&=\ts\frac{1}{4}\bigl(\lnm{\d v}{2}+\lnm{\d w}{2}\bigr)^2+
\frac{3}{4}\bigl(\lnm{\d v}{2}-\lnm{\d w}{2}\bigr)^2,
\end{align*}
by parts (vi), (vii) of Theorem \ref{cu5thm2}. The lemma follows.
\end{proof}

The following proposition constructs the sequence $(f_n)_{n=0}^\iy$
and proves the a priori estimates we need. At various points in its
proof we shall need $t$ to be smaller than some positive constant
defined in terms of $\ka$ and $A_1,\ldots,A_8$. As a shorthand we
will simply say that this holds as $t\le\ep$, and suppose without
remark that $\ep>0$ has been chosen so that the relevant restriction
holds.

The constants $C_1,\ldots,C_8>0$ appearing in the proposition and
proof are those of Propositions \ref{cu5prop2}, \ref{cu5prop4} and
\ref{cu5prop6}. Note that as $C_1,\ldots,C_8$ depend only on
$A_1,\ldots,A_7$ and $m$, it is all right for $\ep,K,F_2,\ldots,F_5$
to depend on them.

\begin{prop} There exist\/ $\ep,K,F_1,\ldots,F_4>0$ depending only
on $m,\ka$ and\/ $A_1,\ldots,A_8$, such that if\/ $0<t\le\ep$
then there is a unique sequence $(f_n)_{n=0}^\iy$ in $C^\iy(N)$
with\/ $f_0=0$ satisfying \eq{cu5eq30} and\/ $\int_Nf_n\,\d V=0$
for all\/ $n\ge 1$ and
\begin{equation*}
\begin{array}{llll}
{\rm(A)\!\!\!}& \lnm{\d f_n}{2}\le F_1t^{\ka+m/2},
&{\rm(a)\!\!\!}& \lnm{\d f_n-\d f_{n-1}}{2}\le F_12^{-n}t^{\ka+m/2},\\
{\rm(B)\!\!\!}& \cnm{\d f_n}{0}\le Kt^\ka\le C_1t,
&{\rm(b)\!\!\!}& \cnm{\d f_n-\d f_{n-1}}{0}\le K2^{-n}t^\ka,\\
{\rm(C)\!\!\!}& \lnm{\na^2f_n}{2}\le F_2t^{\ka+m/2-1},
&{\rm(c)\!\!\!}& \lnm{\na^2f_n-\na^2f_{n-1}}{2}\le
F_22^{-n}t^{\ka+m/2-1},\\
{\rm(D)\!\!\!}& \cnm{\na^2f_n}{0}\!\le\!F_3t^{\ka-1}\!\le\!C_2,
&{\rm(d)\!\!\!}& \cnm{\na^2f_n-\na^2f_{n-1}}{0}\le F_32^{-n}t^{\ka-1},\\
{\rm(E)\!\!\!}& \lnm{\na^3f_n}{2m}\!\le\!F_4t^{\ka-3/2},
&{\rm(e)\!\!\!}& \lnm{\na^3f_n-\na^3f_{n-1}}{2m}\le F_42^{-n}t^{\ka-3/2}.
\end{array}
\end{equation*}
\label{cu5prop7}
\end{prop}

\begin{proof} First note that as $f_0=0$ we have
$f_k=\sum_{n=1}^k(f_n-f_{n-1})$. Suppose that (a) holds
for $n=1,\ldots,k$. Then we have
\begin{equation*}
\lnm{\d f_k}{2}\le\sum_{n=1}^k
\lnm{\d f_n-\d f_{n-1}}{2}
\le F_1t^{\ka+m/2}\sum_{n=1}^k2^{-n}
\le F_1t^{\ka+m/2}.
\end{equation*}
Therefore (a) for $n=1,\ldots,k$ implies (A) for $n=k$.
In the same way, (b)--(e) for $n=1,\ldots,k$ imply (B)--(E)
for $n=k$. The extra inequalities $Kt^\ka\le C_1t$ and
$F_3t^{\ka-1}\le C_2$ in (B), (D) hold as $\ka>1$ and~$t\le\ep$.

Next suppose that $f_1,\ldots,f_n$ exist and (a), (c), (e)
hold for $n$, for some $F_1,F_2,F_4>0$ depending only on
$m,\ka,A_1,\ldots,A_8$. We shall prove (b), (d) for $n$.
Apply Proposition \ref{cu5prop4} to $\al=\d f_n-\d f_{n-1}$.
Equation \eq{cu5eq12} yields
\begin{align*}
\cnm{\na^2f_n\!-\!\na^2f_{n\!-\!1}}{0}&\!\le\!
C_6\bigl(t^{1/2}\lnm{\na^3f_n\!-\!\na^3f_{n\!-\!1}}{2m}\!+\!
t^{-\!m/2}\lnm{\na^2f_n\!-\!\na^2f_{n\!-\!1}}{2}\bigr)\\
&\!\le\!C_6\bigl(t^{1/2}F_42^{-n}t^{\ka\!-\!3/2}\!+\!
t^{-\!m/2}F_22^{-\!n}t^{\ka\!+\!m/2\!-\!1}\bigr)\!
=\!F_32^{-\!n}t^{\ka\!-\!1}
\end{align*}
for $n$ by parts (c), (e), where $F_3=C_6(F_4+F_2)$. This
proves (d). Similarly, \eq{cu5eq11} and parts (a), (d)
prove part (b), with~$K=C_5(F_3+F_1)$.

Therefore, to complete the proof we only need to show that
a unique sequence $(f_n)_{n=0}^\iy$ exists and satisfies
\eq{cu5eq30}, $\int_Nf_n\,\d V=0$ and parts (a), (c) and
(e) for all $n\ge 1$. We will do this by induction on $n$.
The first step is the case $n=1$. As $f_0=0$ and $Q(0)=0$
by Lemma \ref{cu5lem2}, equation \eq{cu5eq30} gives
\e
\d^*\bigl(\psi^m\cos\th\,\d f_1\bigr)=\psi^m\sin\th.
\label{cu5eq31}
\e
As in the proof of Lemma \ref{cu5lem2} $\int_N\psi^m\sin\th\,\d V=0$,
so Lemma \ref{cu5lem4} shows there exists a unique $f_1\in C^\iy(N)$
satisfying \eq{cu5eq31} and~$\int_Nf_1\,\d V=0$.

Let $w=\pi_\sW(f_1)$. Multiplying \eq{cu5eq31} by $f_1$ and
integrating over $N$ yields
\begin{align*}
\ha A_3^m&\lnm{\d f_1}{2}^2\le
\int_N\psi^m\cos\th\ms{\d f_1}\d V
=\int_Nf_1\psi^m\sin\th\,\d V\\
&=\int_N(f_1-w)\psi^m\sin\th\,\d V+
\int_Nw\,\pi_\sW(\psi^m\sin\th)\,\d V\\
&\le\lnm{f_1\!-\!w}{2m/(m-2)}\lnm{\psi^m\sin\th}{2m/(m+2)}
\!+\!\cnm{w}{0}\blnm{\pi_\sW(\psi^m\sin\th)}{1}\\
&\le A_7\lnm{\d f_1-\d w}{2}\cdot A_2t^{\ka+m/2}+
A_8t^{1-m/2}\lnm{\d w}{2}\cdot A_2t^{\ka+m-1}\\
&\le A_7\cdot 2\lnm{\d f_1}{2}\cdot A_2t^{\ka+m/2}+
A_8t^{1-m/2}\cdot 2\lnm{\d f_1}{2}\cdot A_2t^{\ka+m-1}.
\end{align*}
Here the first line uses part (ii) of Theorem \ref{cu5thm2} and
$\cos\th\ge\ha$ from Definition \ref{cu5def1}, the second the
fact that $\an{w,\psi^m\sin\th}=\an{w,\pi_\sW(\psi^m\sin\th)}$
as $w\in W$, the third H\"older's inequality, the fourth parts
(i), (vi) and (vii) of Theorem \ref{cu5thm2}, and the fifth
Lemma \ref{cu5lem5} with $v=f_1-w$. Therefore $\lnm{\d f_1}{2}
\le 4A_2A_3^{-m}(A_7\!+\!A_8)t^{\ka+m/2}$, which proves part
(a) for $n=1$ with~$F_1=8A_2A_3^{-m}(A_7\!+\!A_8)$.

Now apply Proposition \ref{cu5prop6} with $u=f_1$ and
$v=\psi^m\sin\th$ by \eq{cu5eq31}, to get
\begin{align*}
\lnm{\na^2f_1}{2}&\le C_7\bigl(t^{-1}\lnm{\d f_1}{2}
+\lnm{\psi^m\sin\th}{2}\bigr)\\
&\le C_7\bigl(t^{-1}\ha F_1t^{\ka+m/2}
+A_2t^{\ka+m/2-1}\bigr)
=\ha F_2t^{\ka+m/2-1}\quad\text{and}\\
\lnm{\na^3f_1}{2m}
&\le C_8\bigl(t^{-(m+3)/2}\lnm{\d f_1}{2}
+t^{-1}\lnm{\psi^m\sin\th}{2m}+\lnm{\d\psi^m\sin\th}{2m}\bigr)\\
&\le\!C_8\bigl(t^{-(m+3)/2}\ha F_1t^{\ka+m/2}
\!+\!t^{-1}A_2t^{\ka-1/2}\!+\!A_2t^{\ka-3/2}\bigr)\!=\!\ha F_4t^{\ka-3/2},
\end{align*}
where $F_2=C_7(F_1+2A_2)$ and $F_4=C_8(F_1+4A_2)$. Here we have used
(a) when $n=1$ and part (i) of Theorem \ref{cu5thm2}, noting that
$\lnm{\psi^m\sin\th}{2m/(m+2)}\le A_2t^{\ka+m/2}$ and $\cnm{\psi^m
\sin\th}{0}\le A_2t^{\ka-1}$ imply that $\lnm{\psi^m\sin\th}{2}\le
A_2t^{\ka+m/2-1}$ and $\lnm{\psi^m\sin\th}{2m}\le A_2t^{\ka-1/2}$
by interpolation. This proves (c) and (e) for $n=1$, completing the
first step.

For the inductive step, suppose by induction that $k\ge 1$ and
that $f_1,\ldots,f_k$ exist and satisfy \eq{cu5eq30},
$\int_Nf_n\,\d V=0$, (A)--(E) and (a)--(e) for $n=1,\ldots,k$.
We shall show that there exists a unique $f_{k+1}$ satisfying
\eq{cu5eq30}, $\int_Nf_n\,\d V=0$ and parts (a), (c) and (e)
for~$n=k+1$.

By (B) we have $\cnm{\d f_k}{0}\le C_1t$, and $C_1<A_1$ by
Proposition \ref{cu5prop2}, so $\cnm{\d f_k}{0}<A_1t$, and
$Q(\d f_k)$ is well-defined. Also $\int_NQ(\d f_k)\,\d V=0$
by Lemma \ref{cu5lem2}, and $\int_N\psi^m\sin\th\,\d V=0$
as above. Therefore by Lemma \ref{cu5lem4} there exists a
unique $f_{k+1}\in C^\iy(N)$ satisfying \eq{cu5eq30} for
$n=k+1$ and~$\int_Nf_{k+1}\,\d V=0$.

Let $u=f_{k+1}-f_k$. Subtracting \eq{cu5eq30} for $n=k+1$,
$n=k$ gives
\e
\d^*\bigl(\psi^m\cos\th\,\d u\bigr)\!=\!
\d^*\bigl(\psi^m\cos\th\,(\d f_{k+1}\!-\!\d f_k)\bigr)\!=\!
Q(\d f_k)\!-\!Q(\d f_{k-1})\!=\!v,
\label{cu5eq32}
\e
say. We shall estimate some norms of~$v$.

\begin{lem} There exist\/ $G_1,G_2>0$ depending only
on $C_3,C_4,K,F_1,\ldots,F_4$ and\/ $m$ such that\/
$\lnm{v}{p}\le G_12^{-k}t^{2\ka+m/p-2}$ for $p\ge 1$,
$\cnm{v}{0}\le G_12^{-k}t^{2\ka-2}$
and\/~$\lnm{\d v}{2m}\le G_22^{-k}t^{2\ka-5/2}$.
\label{cu5lem6}
\end{lem}

\begin{proof} Observe that (B) and (D) for $n=k,k-1$ give
\begin{equation*}
\cnm{\d f_k}{0},\cnm{\d f_{k-1}}{0}\le C_1t
\quad\text{and}\quad
\cnm{\na^2f_k}{0},\cnm{\na^2f_{k-1}}{0}\le C_2.
\end{equation*}
So we may apply Proposition \ref{cu5prop2} with
$\al=\d f_k$ and $\be=\d f_{k-1}$, and \eq{cu5eq6} gives
\e
\begin{split}
\md{v}\le C_3&\bigl(t^{-1}\md{\d f_k-\d f_{k-1}}+
\md{\na^2f_k-\na^2f_{k-1}}\bigr)\cdot\\
&\bigl(t^{-1}\md{\d f_k}+t^{-1}\md{\d f_{k-1}}+
\md{\na^2f_k}+\md{\na^2f_{k-1}}\bigr).
\end{split}
\label{cu5eq33}
\e

Integrating \eq{cu5eq33} over $N$ and using parts (A), (C), (a),
(c) for $n=k,k-1$ shows that $\lnm{v}{1}\le 4C_3(F_1+F_2)^22^{-k}
t^{2\ka+m-2}$. Taking the supremum of \eq{cu5eq33} over $N$ and
using parts (B), (D), (b), (d) for $n=k,k-1$ gives $\cnm{v}{0}\le
4C_3(K+F_3)^22^{-k}t^{2\ka-2}$. Define $G_1=4C_3\max\bigl(
(F_1+F_2)^2,(K+F_3)^2\bigr)$. Then $\lnm{v}{1}\le G_12^{-k}
t^{2\ka+m-2}$ and $\cnm{v}{0}\le G_12^{-k}t^{2\ka-2}$, so
$\lnm{v}{p}\le G_12^{-k}t^{2\ka+m/p-2}$ for $p\ge 1$ by
interpolation. A similar but longer proof using
\eq{cu5eq7} and (A)--(E), (a)--(e) for $n=k,k-1$ gives
$G_2$ with~$\lnm{\d v}{2m}\le G_22^{-k}t^{2\ka-5/2}$.
\end{proof}

Let $w=\pi_\sW(u)$. Multiplying \eq{cu5eq32} by $u$ and
integrating over $N$ yields
\begin{align*}
\ha A_3^m&\lnm{\d u}{2}^2\!\le\!\int_N\psi^m\cos\th\ms{\d u}\d V
\!=\!\int_Nuv\,\d V\!=\!\int_N(u-w)v\,\d V\!+\!\int_Nwv\,\d V\\
&\le\lnm{u-w}{2m/(m-2)}\lnm{v}{2m/(m+2)}
+\cnm{w}{0}\lnm{v}{1}\\
&\le A_7\lnm{\d u\!-\!\d w}{2}\cdot G_12^{-k}t^{2\ka+m/2-1}
\!+\!A_8t^{1-m/2}\lnm{\d w}{2}\cdot G_12^{-k}t^{2\ka+m-2}\\
&\le A_7\cdot 2\lnm{\d u}{2}\cdot G_12^{-k}t^{2\ka+m/2-1}
\!+\!A_8t^{1-m/2}\cdot 2\lnm{\d u}{2}\cdot G_12^{-k}t^{2\ka+m-2}.
\end{align*}
Here the first line uses part (ii) of Theorem \ref{cu5thm2} and
$\cos\th\ge\ha$ from Definition \ref{cu5def1}, the second
H\"older's inequality, the third parts (vi) and (vii) of
Theorem \ref{cu5thm2} and Lemma \ref{cu5lem6}, and the fourth
Lemma \ref{cu5lem5} with $u-w$ in place of~$v$.

Cancelling shows that $\lnm{\d u}{2}\le 4G_1A_3^{-m}(A_7\!+\!A_8)
2^{-k}t^{2\ka+m/2-1}$. As $u=f_{k+1}-f_k$ and $F_1=8A_2A_3^{-m}
(A_7\!+\!A_8)$, this implies part (a) with $n=k+1$ if
$t^{\ka-1}G_1\le A_2$, which is true as $\ka>1$ and $t\le\ep$.
Now apply Proposition \ref{cu5prop6} to \eq{cu5eq32}, so that
\eq{cu5eq20}, part (a) with $n=k+1$ and Lemma \ref{cu5lem6} give
\begin{align*}
\lnm{\na^2f_{k+1}-\na^2f_k}{2}&\le
C_7\bigl(t^{-1}\lnm{\d f_{k+1}-\d f_k}{2}+\lnm{v}{2}\bigr)\\
&\le C_7\bigl(t^{-1}F_12^{-k-1}t^{\ka+m/2}
+G_12^{-k}t^{2\ka+m/2-2}\bigr).
\end{align*}
Since $F_2=C_7(F_1+2A_2)$, this implies part (c) with $n=k+1$ if
$G_1t^{\ka-1}\le A_2$, which is true as $\ka>1$ and~$t\le\ep$.

In the same way, from \eq{cu5eq21}, part (a) with $n=k+1$ and
Lemma \ref{cu5lem6} we get
\begin{gather*}
\lnm{\na^3f_{k+1}\!-\!\na^3f_k}{2m}
\le C_8\bigl(t^{-(m+3)/2}
\lnm{\d f_{k+1}-\d f_k}{2}
\!+\!t^{-1}\lnm{v}{2m}\!+\!\lnm{\d v}{2m}\bigr)\\
\le C_8\bigl(t^{-(m+3)/2}F_12^{-k-1}t^{\ka+m/2}
\!+\!t^{-1}G_12^{-k}t^{2\ka-3/2}
\!+\!G_22^{-k}t^{2\ka-5/2}\bigr).
\end{gather*}
As $F_4=C_8(F_1+4A_2)$, this implies (e) with $n=k+1$ if
$(G_1+G_2)t^{\ka-1}\le 2A_2$, which holds as $\ka>1$ and
$t\le\ep$. This completes the inductive step, and the proof.
\end{proof}

We can now finish the proof of Theorem \ref{cu5thm2}.
If $u\in C^1(N)$ and $\int_Nu\,\d V=0$ then $u(x)=0$
for some $x\in N$. If $y\in N$, then $x$ and $y$ are
joined by a smooth path $\ga$ in $N$ of length no more
than $\diam(N)$, the diameter of $N$. Integrating along
$\ga$ shows that $\bmd{u(y)}\le\diam(N)\cnm{\d u}{0}$.
Hence $\cnm{u}{0}\le\diam(N)\cnm{\d u}{0}$.

Applying this to part (b) of Proposition \ref{cu5prop7}
gives
\begin{equation*}
\cnm{f_n-f_{n-1}}{0}\le K\diam(N)2^{-n}t^\ka
\quad\text{for $n\ge 1$.}
\end{equation*}
Combining this with parts (b) and (d) shows that
$(f_n)_{n=0}^\iy$ is a {\it Cauchy sequence} in
the Banach space $C^2(N)$. Let $f$ be the limit of
$(f_n)_{n=0}^\iy$ in~$C^2(N)$.

Part (B) gives $\cnm{\d f_n}{0}\le Kt^\ka\le C_1t<A_1t$, as
$C_1<A_1$ by Proposition \ref{cu5prop2}. Taking the limit as
$n\ra\iy$ gives $\cnm{\d f}{0}\le Kt^\ka<A_1t$. Therefore
$Q(\d f)$ is well-defined. Since $Q(\al)$ depends pointwise
on $\al,\na\al$ and $f_n\ra f$ in $C^2(N)$ it follows that
$Q(\d f_n)\ra Q(\d f)$ in $C^0(N)$ as $n\ra\iy$. So
taking the limit $n\ra\iy$ in \eq{cu5eq30} we see that
\eq{cu5eq29} holds for~$f$.

By Definition \ref{cu5def2} and Proposition \ref{cu5prop1},
equation \eq{cu5eq29} is equivalent to
\e
F'\bigl(x,\d f(x),\na^2f(x)\bigr)\equiv 0.
\label{cu5eq34}
\e
This is a {\it second-order nonlinear elliptic equation}
on $f$. Note that $F'$ is not linear in $\na^2f(x)$, so
that \eq{cu5eq34} is not quasilinear, and that $F'$ is a
smooth function of its arguments. Now Aubin
\cite[Th.~3.56]{Aubi} gives a regularity result for $C^2$
solutions of such equations, which implies that $f\in
C^\iy(N)$, as we wish.

We have constructed $\ep,K>0$ in Proposition \ref{cu5prop7}
depending only on $m,\ka$ and $A_1,\ldots,A_8$, and
$f\in C^\iy(N)$ with $\cnm{\d f}{0}\le Kt^\ka<A_1t$
satisfying \eq{cu5eq29}. Taking the limit $n\ra\iy$ in
$\int_Nf_n\,\d V=0$ yields $\int_Nf\,\d V=0$. Finally,
\eq{cu5eq29} and Lemma \ref{cu5lem2} show that
$\Phi_*\bigl(\Ga(\d f)\bigr)$ is an immersed special
Lagrangian $m$-fold in $(M,J,\om,\Om)$. This concludes
the proof of Theorem~\ref{cu5thm2}.

\section{Desingularization: the simplest case}
\label{cu6}

Let $(M,J,\om,\Om)$ be an almost Calabi--Yau $m$-fold, $X$
be a compact SL $m$-fold in $M$ with conical singularities
at $x_1,\ldots,x_n$ with cones $C_i$, and $L_1,\ldots,L_n$
be AC SL $m$-folds in $\C^m$, where $L_i$ has cone $C_i$
and rate $\la_i$ for $i=1,\ldots,n$. The goal of the rest
of the paper is to {\it desingularize} $X$ by `gluing'
$L_1,\ldots,L_n$ in at $x_1,\ldots,x_n$, to produce a
family of compact, nonsingular SL $m$-folds $\ti N^t$
for $t\in(0,\ep]$, which converge to $X$ as $t\ra 0$
in an appropriate sense.

Very briefly, we do this by first shrinking $L_i$ by a small
factor $t>0$ and gluing $tL_i$ into $X$ at $x_i$ to make a family
of compact Lagrangian $m$-folds $N^t$ for $t\in(0,\de)$, and
then applying Theorem \ref{cu5thm2} to show that $N^t$ can be
deformed to a nearby SL $m$-fold $\ti N^t$ when $0<t\le\ep<\de$.
Now to do this in full generality is rather complex. Therefore
we begin in this section with the easiest case, in which $\la_i<0$
for all $i$ and $X'=X\sm\{x_1,\ldots,x_n\}$ is connected.

As explained in \S\ref{cu1}, this simplifies the problem,
avoiding issues of small eigenvalues and obstructions to the
existence of $N^t$ from global symplectic topology. Section
\ref{cu7} will extend the results to the case when $X'$ is
not connected. The sequel \cite{Joyc6} will study the case
when $\la_i=0$ and $Y(L_i)\ne 0$, and extend the results to
{\it families} of almost Calabi--Yau $m$-folds $(M,J^s,\om^s,\Om^s)$
for~$s\in{\cal F}$.

\subsection{Setting up the problem}
\label{cu61}

We shall consider the following situation.

\begin{dfn} Let $(M,J,\om,\Om)$ be an almost Calabi--Yau
$m$-fold with metric $g$, and define $\psi:M\ra(0,\iy)$
as in \eq{cu2eq3}. Let $X$ be a compact SL $m$-fold in
$M$ with conical singularities at $x_1,\ldots,x_n$ with
identifications $\up_1,\ldots,\up_n$, cones $C_1,\ldots,C_n$
and rates $\mu_1,\ldots,\mu_n\in(2,3)$, as in Definition
\ref{cu3def3}. Define $\Si_i=C_i\cap{\cal S}^{2m-1}$ for
$i=1,\ldots,n$. Let $L_1,\ldots,L_n$ be AC SL $m$-folds in
$\C^m$, where $L_i$ has cone $C_i$ and rate $\la_i$ for
$i=1,\ldots,n$, as in Definition \ref{cu4def1}. Suppose
that $\la_i<\ha(2-m)$ for~$i=1,\ldots,n$.

We use the following notation:
\begin{itemize}
\setlength{\parsep}{0pt}
\setlength{\itemsep}{0pt}
\item Let $R,B_R,X'$ and $\iota_i,\Up_i$ for $i=1,\ldots,n$ be
as in Definition~\ref{cu3def3}.
\item Let $\ze$ and $U_\sCi,\Phi_\sCi$ for
$i=1,\ldots,n$ be as in Theorem~\ref{cu3thm2}.
\item Let $R',K$ and $\phi_i,\eta_i,\eta_i^1,\eta_i^2,S_i$
for $i=1,\ldots,n$ be as in Theorem~\ref{cu3thm3}.
\item Let $U_\sXp,\Phi_\sXp$ be as in Theorem~\ref{cu3thm5}.
\item Let $A_i$ be as in Theorem \ref{cu3thm4} for $i=1,\ldots,n$,
so that~$\eta_i=\d A_i$.
\item Apply Theorem \ref{cu4thm1} to $L_i$ with $\ze,U_\sCi,
\Phi_\sCi$ as above, for $i=1,\ldots,n$. Let $T>0$ be as in
the theorem, the same for all $i$. Let the subset $K_i\subset L_i$,
the diffeomorphism $\vp_i:\Si_i\t(T,\iy)\ra L_i\sm K_i$ and the
1-form $\chi_i$ on $\Si_i\t(T,\iy)$ with components $\chi_i^1,
\chi_i^2$ be as in Theorem~\ref{cu4thm1}.
\item Let $U_\sLi,\Phi_\sLi$ be as in Theorem \ref{cu4thm3}
for~$i=1,\ldots,n$.
\item Let $E_i\in C^\iy\bigl(\Si_i\t(T,\iy)\bigr)$ be as in
Theorem \ref{cu4thm2} for~$i=1,\ldots,n$.
\end{itemize}
\label{cu6def1}
\end{dfn}

In \S\ref{cu64}--\S\ref{cu65} we will also suppose that $X'$ is
{\it connected}, and in \S\ref{cu65} we will relax the condition
$\la_i<\ha(2-m)$ to $\la_i<0$ using Theorem \ref{cu4thm2}. With
this notation, we define a family of Lagrangian $m$-folds $N^t$
in $(M,\om)$ for~$t\in(0,\de)$.

\begin{dfn} In the situation of Definition \ref{cu6def1}, choose a
smooth, increasing function $F:(0,\iy)\ra[0,1]$ with $F(r)\equiv 0$
for $r\in(0,1)$ and $F(r)\equiv 1$ for $r>2$. Write $F'$ for
$\d F/\d r$. Let $t>0$ act as a {\it dilation} ${\bf x}\mapsto
t\,{\bf x}$ on $\C^m$. Write $tK_i,tL_i$ for the images of $K_i,L_i$
under $t$, so that $tK_i=\{t\,{\bf x}:{\bf x}\in K_i\}$, and so on.
Let $\tau\in(0,1)$ satisfy
\e
\max_{i=1,\ldots,n}\Bigl(\frac{m+2}{2\mu_i+m-2}\Bigr)<\tau<1,
\label{cu6eq1}
\e
which is possible as $\mu_i>2$ implies~$(m+2)/(2\mu_i+m-2)<1$.

For $i=1,\ldots,n$ and small enough $t>0$, define $P_i^t=\Up_i(tK_i)$.
This is well-defined if $tK_i\subset B_R\subset\C^m$, and is a
compact submanifold of $M$ with boundary, diffeomorphic to $K_i$. As
$K_i$ is Lagrangian in $(\C^m,\om')$ and $\Up_i^*(\om)=\om'$,
we see that $P_i^t$ is {\it Lagrangian} in~$(M,\om)$.

For $i=1,\ldots,n$ and $t>0$ with $tT<t^\tau<2t^\tau<R'$,
define a 1-form $\xi_i^t$ on $\Si_i\t(tT,R')$ by 
\e
\begin{split}
\xi_i^t(\si,r)&=\d\bigl[F(t^{-\tau}r)A_i(\si,r)+
t^2(1-F(t^{-\tau}r))E_i(\si,t^{-1}r)\bigr]\\
&=F(t^{-\tau}r)\eta_i(\si,r)+
t^{-\tau}F'(t^{-\tau}r)A_i(\si,r)\d r\\
&\quad+t^2(1-F(t^{-\tau}r))\chi_i(\si,t^{-1}r)-
t^{2-\tau}F'(t^{-\tau}r)E_i(\si,t^{-1}r)\d r.
\end{split}
\label{cu6eq2}
\e
Let $\xi_i^{1,t},\xi_i^{2,t}$ be the components of $\xi_i^t$
in $T^*\Si$ and $\R$, as for $\eta_i,\chi_i$ in Theorems
\ref{cu3thm3} and \ref{cu4thm1}. Note that when $r\ge 2t^\tau$
we have $F(t^{-\tau}r)\equiv 1$ so that $\xi_i^t(\si,r)=\eta_i
(\si,r)$, and when $r\le t^\tau$ we have $F(t^{-\tau}r)\equiv 0$,
so that $\xi_i^t(\si,r)=t^2\chi_i(\si,t^{-1}r)$. Thus $\xi_i^t$
is an exact 1-form which interpolates between $\eta_i(\si,r)$
near $r=R'$ and $t^2\chi_i(\si,t^{-1}r)$ near~$r=tT$.

Choose $\de\in(0,1]$ with $\de T\le\de^\tau<2\de^\tau\le R'$
and $\de K_i\!\subset\!B_R\!\subset\!\C^m$ and
\e
\bmd{\xi_i^t(\si,r)}<\ze r\quad\text{on $\Si_i\t(tT,R')$ for
all $i=1,\ldots,n$ and $t\in(0,\de)$.}
\label{cu6eq3}
\e
Here $t\in(0,\de)$ and $\de T\!\le\!\de^\tau\!<\!2\de^\tau
\!\le\!R'$ imply that $tT\!<\!t^\tau\!<\!2t^\tau\!<\!R'$,
so $\xi_i^t$ exists. As $\md{\eta_i(\si,r)}<\ze r$ in Theorem
\ref{cu3thm3} and $\xi_i^t\equiv\eta_i$ when $r\ge 2t^\tau$,
equation \eq{cu6eq3} holds automatically on $\Si_i\t[2t^\tau,R')$.
Similarly, as $\bmd{\chi_i(\si,r)}<\ze r$ in Theorem \ref{cu4thm1}
and $\xi_i^t(\si,r)=t^2\chi_i(\si,t^{-1}r)$ when $r\le t^\tau$,
equation \eq{cu6eq3} holds on $\Si_i\t(tT,t^\tau]$. We can show
using \eq{cu3eq4}, \eq{cu4eq3} and \eq{cu6eq2} that \eq{cu6eq3}
holds on $\Si_i\t(t^\tau,2t^\tau)$ for small enough $t>0$, so
$\de$ exists.

For $i=1,\ldots,n$ and $t\in(0,\de)$, define
$\Xi_i^t:\Si_i\t(tT,R')\ra M$ by
\e
\Xi_i^t(\si,r)=\Up_i\circ\Phi_\sCi\bigl(\si,r,
\xi_i^{1,t}(\si,r),\xi_i^{2,t}(\si,r)\bigr).
\label{cu6eq4}
\e
Then $\Phi_\sCi\bigl(\ldots)$ is well-defined as
$\bmd{\xi_i^t(\si,r)}<\ze r$ by \eq{cu6eq3}, and making $R'$
smaller if necessary we can ensure that $\Phi_\sCi\bigl(
\ldots)$ lies in $B_R$, so $\Xi_i^t$ is well-defined, and is
an embedding as $\Up_i,\Phi_\sCi$ are.

Define $Q_i^t=\Xi_i^t\bigl(\Si_i\t(tT,R')\bigr)$ for $i=1,\ldots,n$
and $t\in(0,\de)$. As $\Up_i^*(\om)=\om'$, $\Phi_\sCi^*(\om')
=\hat\om$ and $\xi_i^t$ is a closed 1-form we see that $(\Xi_i^t)^*
(\om)\equiv 0$. Thus $Q_i^t$ is {\it Lagrangian} in $(M,\om)$, and is
a noncompact embedded submanifold diffeomorphic to $\Si_i\t(tT,R')$.
For $t\in(0,\de)$, define $N^t$ to be the disjoint union of $K$,
$P_1^t,\ldots,P_n^t$ and $Q_1^t,\ldots,Q_n^t$, where $K\subset
X'$ is as above.

Then $N^t$ is {\it Lagrangian} in $(M,\om)$, as $K,P_i^t$ and
$Q_i^t$ are. We claim that $N^t$ is a compact, smooth submanifold
of $M$ {\it without boundary}. That is, the boundary $\pd P_i^t
\cong\Si_i$ joins smoothly onto $Q_i^t\cong\Si_i\t(tT,R')$ at the
$\Si_i\t\{tT\}$ end, and the boundary $\pd K$ is the disjoint
union of pieces $\Si_i$ for $i=1,\ldots,n$ which join smoothly
onto $Q_i^t\cong\Si_i\!\t\!(tT,R')$ at the $\Si_i\!\t\!\{R'\}$ end.

To see this, note that $\xi_i^t(\si,r)=t^2\chi_i(\si,t^{-1}r)$
on $\Si_i\!\t\!(tT,t^\tau]$, and so
\begin{align*}
\Xi_i^t(\si,r)&=\Up_i\circ\Phi_\sCi\bigl(\si,r,
t^2\chi_i^1(\si,t^{-1}r),t\chi_i^2(\si,t^{-1}r)\bigr)\\
&=\Up_i\bigl(t\,\Phi_\sCi\bigl(\si,t^{-1}r,
\chi_i^1(\si,t^{-1}r),\chi_i^2(\si,t^{-1}r)\bigr)\bigr)
=\Up_i\bigl(t\,\vp_i(\si,t^{-1}r)\bigr)
\end{align*}
on $\Si_i\!\t\!(tT,t^\tau]$, using \eq{cu4eq2} and the dilation
equivariance of $\Phi_{\sst C_i}$ in Theorem \ref{cu3thm2}. Thus
the end $\Xi_i^t\bigl(\Si_i\!\t\!(tT,t^\tau]\bigr)$ of $Q_i^t$
is $\Up_i\bigl(t\,\vp_i(\Si_i\!\t\!(T,t^{\tau-1}])\bigr)\subset
\Up_i(tL_i)$, and as $\vp_i\bigl(\Si_i\!\t\!(T,t^{\tau-1}]\bigr)
\subset L_i$ joins smoothly onto $K_i\subset L_i$ we see that
$\Xi_i^t\bigl(\Si_i\!\t\!(tT,t^\tau]\bigr)\subset Q_i^t$ joins
smoothly onto~$P_i^t=\Up_i(tK_i)$.

In the same way, as $\xi_i^t\equiv\eta_i$ on $\Si_i\t[2t^\tau,R')$
Theorem \ref{cu3thm3} gives
\begin{equation*}
\Xi_i^t(\si,r)=\Up_i\circ\Phi_\sCi\bigl(\si,r,
\eta_i^1(\si,r),\eta_i^2(\si,r)\bigr)=\Up_i\circ\phi_i(\si,r)
\;\>\text{on}\;\>
\Si_i\t[2t^\tau,R'),
\end{equation*}
so $\Xi_i^t\bigl(\Si_i\!\t\![2t^\tau,R')\bigr)\!\subset\!Q_i^t$
is $\Up_i\circ\phi_i(\Si_i\!\t\![2t^\tau,R'))\bigr)\subset
S_i\subset X'$, which joins smoothly onto $K$. Therefore $N^t$
is compact and smooth without boundary.
\label{cu6def2}
\end{dfn}

Here $N^t$ depends smoothly on $t$, and converges to the singular
SL $m$-fold $X$ in $M$ as $t\ra 0$, in the sense of currents in
Geometric Measure Theory. Also $N^t$ is equal to $X'$ in $K$ and
the annuli $\Xi_i^t\bigl(\Si_i\!\t\![2t^\tau,R')\bigr)\subset
Q_i^t$, and equal to $\Up_i(tL_i)$ on $P_i^t=\Up_i(tK_i)$ and the
annuli $\Xi_i^t\bigl(\Si_i\!\t\!(tT,t^\tau]\bigr)\subset Q_i^t$.
In between, on the annuli $\Xi_i^t\bigl(\Si_i\!\t\!(t^\tau,2t^\tau)
\bigr)\subset Q_i^t$, $N^t$ interpolates smoothly between $X'$ and
$\Up_i(tL_i)$ as a Lagrangian submanifold.

At several points later on we shall make $\de>0$ smaller if
necessary to ensure that something works for all $t\in(0,\de)$.
This is for simplicity, to avoid introducing a series of further
constants $\de',\de'',\ldots$ in~$(0,\de)$.

\subsection{Estimating $\Im\Om\vert_{N^t}$}
\label{cu62}

We now prove estimates for $\Im\Om\vert_{N^t}$, to use in
part (i) of Theorem \ref{cu5thm2}. Here is some more notation.

\begin{dfn} In the situation of Definitions \ref{cu6def1} and
\ref{cu6def2}, let $h^t$ be the restriction of $g$ to $N^t$
for $t\in(0,\de)$, so that $(N^t,h^t)$ is a compact Riemannian
manifold. As $X',L_i$ are SL $m$-folds they are oriented, and
$N^t$ is made by gluing $X',L_1,\ldots,L_n$ together in an
orientation-preserving way, so $N^t$ is also oriented. Let
$\d V^t$ be the volume form on $N^t$ induced by $h^t$ and this
orientation. As in \eq{cu5eq1} we may write $\Om\vert_{N^t}=
\psi^m{\rm e}^{\smash{i\th^t}}\,\d V^t$ for some phase
function ${\rm e}^{\smash{i\th^t}}$ on $N^t$. Write
$\ve^t=\psi^m\sin\th^t$, so that $\Im\Om\vert_{N^t}=\ve^t\,
\d V^t$ for~$t\in(0,\de)$.
\label{cu6def3}
\end{dfn}

We compute bounds for $\ve^t$ at each point in~$N^t$.

\begin{prop} In the situation above, making $\de>0$ smaller if
necessary, there exists $C>0$ such that for all\/ $t\in(0,\de)$
we have $\ve^t=0$ on $K$, the pull-back\/ $(\Xi_i^t)^*(\ve^t)$
of\/ $\ve^t$ on $Q_i^t$ satisfies
\ea
\bmd{(\Xi_i^t)^*(\ve^t)}(\si,r)&\le
\begin{cases}
Cr, & r\in(tT,t^\tau],\\
Ct^{\tau(\mu_i-2)}+Ct^{(1-\tau)(2-\la_i)},
{}^{\phantom{jn}}   %%%%% to get spacing right in cd6eq4, cd6eq5.
& r\in(t^\tau,2t^\tau),\\
0, & r\in[2t^\tau,R'),
\end{cases}
\label{cu6eq5}\\
\bmd{(\Xi_i^t)^*(\d\ve^t)}(\si,r)&\le
\begin{cases}
C, & r\in(tT,t^\tau],\\
Ct^{\tau(\mu_i-3)}+Ct^{(1-\tau)(2-\la_i)-\tau}, & r\in(t^\tau,2t^\tau),\\
0, & r\in[2t^\tau,R'),
\end{cases}
\label{cu6eq6}\\
\text{and}\quad
\md{\ve^t}&\le Ct,\quad
\md{\d\ve^t}\le C
\quad\text{on $P_i^t$ for all\/ $i=1,\ldots,n$.}
\label{cu6eq7}
\ea
Here $\md{\,.\,}$ is computed using the metrics $(\Xi_i^t)^*(h^t)$
in \eq{cu6eq6} and\/ $h^t$ in~\eq{cu6eq7}.
\label{cu6prop1}
\end{prop}

\begin{proof} Since $N^t$ coincides with $X'$ in $K$ and
$\Xi_i^t\bigl(\Si_i\t[2t^\tau,R')\bigr)$, and $\Im\Om\vert_{X'}$
as $X'$ is special Lagrangian, we see that $\ve^t\equiv 0$ on
$K$ and $\Xi_i^t\bigl(\Si_i\t[2t^\tau,R')\bigr)$, giving the
bottom lines of \eq{cu6eq5} and~\eq{cu6eq6}.

As $\Up_i^*(\Im\Om)$ is a smooth $m$-form on $B_R$ and
$\Up_i^*(\Im\Om)\vert_0=\up_i^*(\Im\Om)=\psi(x_i)^m\Im\Om'$
by Definition \ref{cu3def3}, we see that $\Up_i^*(\Im\Om)=
\psi(x_i)^m\Im\Om'+O(r)$ on $B_R$, by Taylor's Theorem.
Since $tL_i$ is special Lagrangian in $\C^m$ we have
$\Im\Om'\vert_{tL_i}\equiv 0$. Thus
\e
\bmd{\Up_i^*(\Im\Om)\vert_{tL_i}}=O(r)\quad\text{on $tL_i\cap B_R$,}
\label{cu6eq8}
\e
computing $\md{\,.\,}$ using the metric $\Up_i^*(g)$ on $B_R$,
restricted to~$tL_i$.

Now $N^t$ coincides with $\Up_i(tL_i)$ on $P_i^t$ and
$\Xi_i^t\bigl(\Si_i\t(tT,t^\tau]\bigr)$, so
$\ve^t\d V^t=\Im\Om\vert_{\Up_i(tL_i)}$ on these regions.
As $h^t$ is the restriction of $g$ to $N^t$ we have
$\md{\d V^t}=1$, computing $\md{\,.\,}$ using $g$, so
\e
\bmd{\Up_i^*(\ve^t)}=\bmd{\Up_i^*(\Im\Om)\vert_{tL_i}}
\quad\text{on $t\bigl(K\cup\vp_i(\Si_i\t(T,t^{\tau-1}]\bigr)
\subset tL_i\cap B_R$.}
\label{cu6eq9}
\e

Combining \eq{cu6eq8} and \eq{cu6eq9} gives $\bmd{\ve^t}=
O\bigl((\Up_i)_*(r)\bigr)$ on $P_i^t$ and $\Xi_i^t\bigl(
\Si_i\t(tT,t^\tau]\bigr)$. As $(\Up_i)_*(r)=O(t)$ on
$P_i^t$, we see that
\e
\bmd{(\Xi_i^t)^*(\ve^t)}(\si,r)=O(r)
\quad\text{for $r\in(tT,t^\tau]$, and}\quad
\md{\ve^t}=O(t)\quad\text{on $P_i^t$.}
\label{cu6eq10}
\e
A similar argument for the derivative $\d\ve^t$ gives
\e
\bmd{(\Xi_i^t)^*(\d\ve^t)}(\si,r)=O(1)
\quad\text{for $r\in(tT,t^\tau]$, and}\quad
\md{\d\ve^t}=O(1)\quad\text{on $P_i^t$.}
\label{cu6eq11}
\e

Next we estimate $\ve^t$ and $\d\ve^t$ on the annuli
$\Xi_i^t\bigl(\Si_i\t(t^\tau,2t^\tau)\bigr)$. First
we bound $\xi_i^t$ and its derivatives on $\Si_i\t
(t^\tau,2t^\tau)$. From \eq{cu3eq4} and \eq{cu4eq3}
we find that
\begin{equation*}
\bmd{\na^kA_i(\si,r)}=O(t^{\tau(\mu_i-k)})
\quad\text{and}\quad
\bmd{\na^kE_i(\si,t^{-1}r)}=O(t^{-\la_i+\tau(\la_i-k)})
\end{equation*}
for $r\in(t^\tau,2t^\tau)$, computing $\na,\md{\,.\,}$ using
the cone metric $\iota_i^*(g')$ on $\Si_i\t(t^\tau,2t^\tau)$.
Substituting these into \eq{cu6eq2} gives
\e
\bmd{\na^k\xi_i^t(\si,r)}=
O(t^{\tau(\mu_i-1-k)})+O(t^{2-\la_i+\tau(\la_i-1-k)})
\quad\text{for $r\in(t^\tau,2t^\tau)$,}
\label{cu6eq12}
\e
computing $\na,\md{\,.\,}$ using the cone metric~$\iota_i^*(g')$.

Now $\Up_i^*(\ve^t)$ depends pointwise on $\xi_i^t$ and $\na\xi_i^t$,
and when $\xi_i^t,\na\xi_i^t$ are small the dominant error terms in
$\ve^t$ are linear in $\na\xi_i^t$. Similarly, $\d\ve^t$ depends
pointwise on $\xi_i^t,\na\xi_i^t$ and $\na^2\xi_i^t$, and when
$\xi_i^t,\na\xi_i^t$ are small the dominant error terms in $\d\ve^t$
are linear in $\na^2\xi_i^t$. The metrics $\iota_i^*(g')$ and
$(\Xi_i^t)^*(h^t)$ on $\Si_i\t(t^\tau,2t^\tau)$ are equivalent
uniformly in $t$, so \eq{cu6eq12} also holds with $\na,\md{\,.\,}$
computed using~$(\Xi_i^t)^*(h^t)$.

Putting all this together we see that
\e
\begin{aligned}
\bmd{(\Xi_i^t)^*(\ve^t)}&=
O(t^{\tau(\mu_i-2)})+O(t^{(1-\tau)(2-\la_i)})
&\quad&\text{and}\\
\bmd{(\Xi_i^t)^*(\d\ve^t)}&=
O(t^{\tau(\mu_i-3)})+O(t^{(1-\tau)(2-\la_i)-\tau})
&\quad&\text{for $r\in(t^\tau,2t^\tau)$,}
\end{aligned}
\label{cu6eq13}
\e
computing $\md{\,.\,}$ using $(\Xi_i^t)^*(h^t)$. Here we have
used $\mu_i<3$ and $\tau>0$ to absorb error terms $O(t^\tau)$
and $O(1)$ respectively into the first term on the right hand side
of each line of \eq{cu6eq13}, from the same source as \eq{cu6eq10}
and \eq{cu6eq11}. Making $\de>0$ smaller if necessary, the rest
of \eq{cu6eq5}--\eq{cu6eq7} now follow from \eq{cu6eq10},
\eq{cu6eq11} and \eq{cu6eq13}, for some $C>0$ independent of~$t$.
\end{proof}

Now we can estimate norms of $\ve^t$ and $\d\ve^t$, as in
part (i) of Theorem~\ref{cu5thm2}.

\begin{prop} There exists $C'>0$ such that for all\/ $t\in(0,\de)$
we have
\ea
\lnm{\ve^t}{2m/(m+2)}&\le C't^{\tau(1+m/2)}
\sum_{i=1}^n\bigl(t^{\tau(\mu_i-2)}+t^{(1-\tau)(2-\la_i)}\bigr),
\label{cu6eq14}\\
\cnm{\ve^t}{0}&\le C'
\sum_{i=1}^n\bigl(t^{\tau(\mu_i-2)}+t^{(1-\tau)(2-\la_i)}\bigr),
\label{cu6eq15}\\
\text{and}\quad\lnm{\d\ve^t}{2m}&\le C't^{-\tau/2}
\sum_{i=1}^n\bigl(t^{\tau(\mu_i-2)}+t^{(1-\tau)(2-\la_i)}\bigr),
\label{cu6eq16}
\ea
computing norms with respect to the metric $h^t$ on~$N^t$.
\label{cu6prop2}
\end{prop}

\begin{proof} As $g'$ and $\Up_i^*(g)$ are equivalent metrics
on $B_R$, it is easy to see that there exist $D_1,D_2,D_3>0$
such that for all $t\in(0,\de)$ we have
\e
\begin{gathered}
\vol\bigl(P_i^t\bigr)\le D_1t^m,\quad
\vol\bigl(\Xi_i^t\bigl(\Si_i\t(t^\tau,2t^\tau)\bigr)\bigr)\le
D_2t^{m\tau}\\
\text{and}\quad
(\Xi_i^t)^*(\d V^t)\le D_3\,\d V_{g'}
\quad\text{on $\Si_i\t(tT,t^\tau]$,}
\end{gathered}
\label{cu6eq17}
\e
where $\d V_{g'}$ is the volume form of the cone metric
$\iota_i^*(g')$ on~$\Si_i\t(tT,t^\tau]$.

Combining \eq{cu6eq5}, \eq{cu6eq7} and \eq{cu6eq17} we find that
\begin{align*}
\int_{N^t}\md{\ve^t}^{2m/(m+2)}&\,\d V^t\le D_1t^m(Ct)^{2m/(m+2)}\\
&+D_2t^{m\tau}\sum_{i=1}^n\bigl(Ct^{\tau(\mu_i-2)}
+Ct^{(1-\tau)(2-\la_i)}\bigr)^{2m/(m+2)}\\
&+D_3\sum_{i=1}^n\vol(\Si_i)
\int_{tT}^{t^{\tau}}(Cr)^{2m/(m+2)}r^{m-1}\d r.
\end{align*}
Raising this to the power $(m+2)/2m$ and manipulating, we can
prove \eq{cu6eq14} for some $C'>0$ depending only on $C,D_1,
D_2,D_3,m,n$ and $\vol(\Si_i)$. Equations \eq{cu6eq15} and
\eq{cu6eq16} follow by similar arguments.
\end{proof}

Now for part (i) of Theorem \ref{cu5thm2} to hold, we want
$\lnm{\ve^t}{2m/(m+2)}\le A_2t^{\ka+m/2}$, $\cnm{\ve^t}{0}\le
A_2t^{\ka-1}$ and $\lnm{\d\ve^t}{2m}\le A_2t^{\ka-3/2}$ for some
$\ka>1$. Clearly, as $t<1$ from \eq{cu6eq14}--\eq{cu6eq16} these
hold with $A_2=2nC'$ provided for all $i=1,\ldots,n$ we have
\ea
\tau(1\!+\!m/2)\!+\!\tau(\mu_i\!-\!2)&\!\ge\!\ka\!+\!m/2,&\;\>
\tau(1\!+\!m/2)\!+\!(1\!-\!\tau)(2\!-\!\la_i)&\!\ge\!\ka\!+\!m/2,
\label{cu6eq18}\\
\tau(\mu_i\!-\!2)&\!\ge\!\ka\!-\!1,&\;\>
(1\!-\!\tau)(2\!-\!\la_i)&\!\ge\!\ka\!-\!1,
\label{cu6eq19}\\
-\tau/2\!+\!\tau(\mu_i\!-\!2)&\ge\ka\!-\!3/2,
&\text{and}\;\>
-\tau/2\!+\!(1\!-\!\tau)(2\!-\!\la_i)&\ge\ka\!-\!3/2.
\label{cu6eq20}
\ea

Elementary calculations using $0<\tau<1$, $\mu_i>2$ and $\la_i<2$
show that the first equation of \eq{cu6eq18} admits a solution
$\ka>1$ provided $\tau>(2+m)/(2\mu_i-2+m)$, and the second equation
of \eq{cu6eq18} admits a solution $\ka>1$ provided $\la_i-1+m/2<0$, 
that is, provided $\la_i<\ha(2-m)$ as in Definition \ref{cu6def1}.
Also \eq{cu6eq18} implies \eq{cu6eq19} and \eq{cu6eq20} as~$\tau\le 1$.

Therefore, in Definition \ref{cu6def2} we choose $\tau\in(0,1)$
to satisfy \eq{cu6eq1}. Then there exists $\ka>1$ satisfying
\eq{cu6eq18}--\eq{cu6eq20} for all $i=1,\ldots,n$, and we have proved:

\begin{thm} Making $\de>0$ smaller if necessary, there exist\/
$A_2>0$ and\/ $\ka>1$ such that the functions $\ve^t=\psi^m\sin\th^t$
on $N^t$ satisfy $\lnm{\ve^t}{2m/(m+2)}\le A_2t^{\ka+m/2}$,
$\cnm{\ve^t}{0}\le A_2t^{\ka-1}$ and\/ $\lnm{\d\ve^t}{2m}\le
A_2t^{\ka-3/2}$ for all\/ $t\in(0,\de)$, as in part\/ {\rm(i)}
of Theorem~\ref{cu5thm2}.
\label{cu6thm1}
\end{thm}

Here is why the condition $\la_i<\ha(2-m)$ is needed in Definition
\ref{cu6def1}. In \eq{cu4eq4} the term $Ct^{(1-\tau)(2-\la_i)}$ when
$r\in(t^\tau,2t^\tau)$ contributes $C't^{2-\la_i+\tau(\la_i-1+m/2)}$
to the bound for $\lnm{\ve^t}{2m/(m+2)}$ in \eq{cu6eq14}. Now if
$\la_i\ge\ha(2-m)$ then $2-\la_i+\tau(\la_i-1+m/2)\le 1+m/2$ for
all $\tau\in[0,1]$, so whatever value of $\tau$ we choose our
estimate for $\lnm{\ve^t}{2m/(m+2)}$ will be at least as big
as~$O(t^{1+m/2})$.

But for Theorem \ref{cu5thm2} to work we need $\lnm{\ve^t}{2m/(m+2)}
=O(t^{\ka+m/2})$ for $\ka>1$, that is, we need $\lnm{\ve^t}{2m/(m+2)}$
to be smaller than $O(t^{1+m/2})$. If $\la_i\ge\ha(2-m)$ then $L_i$
does not decay to the cone $C_i$ fast enough at infinity, so the
errors we make in tapering $L_i$ off to $C_i$ are too great for
the method of \S\ref{cu5} to cope with.

The condition $\la_i<\ha(2-m)$ in Definition \ref{cu6def1}
will be relaxed to $\la_i<0$ in \S\ref{cu65} using Theorem
\ref{cu4thm2}, so it is not as strong a restriction as it appears.

\subsection{Lagrangian neighbourhoods and bounds on $R(h^t),\de(h^t)$}
\label{cu63}

Next we show that parts (ii)--(v) of Theorem \ref{cu5thm2} hold for
$N^t$ when $t\in(0,\de)$, with appropriate $A_1,A_3,\ldots,A_6>0$
independent of $t$. We begin by gluing together the Lagrangian
neighbourhoods $U_\sXp,\Phi_\sXp$ for $X'$ in Theorem \ref{cu3thm5}
and $U_\sLi,\Phi_\sLi$ for $L_i$ in Theorem \ref{cu4thm3} to get a
Lagrangian neighbourhood $U_\sNt,\Phi_\sNt$ for $N^t$, which we
will use to define the $m$-form $\be^t$ on $\B_{A_1t}\subset T^*N^t$
in Definition \ref{cu5def1}, and so prove part (iii) of
Theorem~\ref{cu5thm2}.

\begin{dfn} Define an open neighbourhood $U_\sNt\subset
T^*N^t$ of the zero section $N^t$ in $T^*N^t$ and a smooth map
$\Phi_\sNt:U_\sNt\ra M$ as follows. Let $\pi:T^*N^t\ra N^t$ be
the natural projection. As $N^t$ is the disjoint union of $K$
and $P_i^t,Q_i^t$ for $i=1,\ldots,n$ we shall define $U_\sNt$
and $\Phi_\sNt$ separately over $K,P_i^t$ and~$Q_i^t$.

Define $U_\sNt\cap\pi^*(K)$ and $\Phi_\sNt\vert_{U_\sNt\cap\pi^*(K)}$ by
\e
U_\sNt\cap\pi^*(K)=U_\sXp\cap\pi^*(K)\;\>\text{and}\;\>
\Phi_\sNt\vert_{U_\sNt\cap\pi^*(K)}=\Phi_\sXp\vert_{U_\sXp\cap\pi^*(K)},
\label{cu6eq21}
\e
recalling that $K$ is part of $N^t$ and $X'$, so $U_\sNt\cap\pi^*(K)$
and $U_\sXp\cap\pi^*(K)$ are both subsets of $T^*K$. For $i=1,\ldots,n$,
define $U_\sNt\cap\pi^*(P_i^t)$ and $\Phi_\sNt\vert_{U_\sNt\cap
\pi^*(P_i^t)}$ by
\e
\begin{gathered}
U_\sNt\cap\pi^*(P_i^t)=\d(\Up_i\circ t)
\bigl(\{\al\in T^*K_i:t^{-2}\al\in U_\sLi\}\bigr)\\
\text{and}\quad \Phi_\sNt\circ\d(\Up_i\circ t)(\al)=
\Up_i\circ t\circ\Phi_\sLi(t^{-2}\al).
\end{gathered}
\label{cu6eq22}
\e
Here the diffeomorphism $\Up_i\circ t:K_i\ra P_i^t$ induces
$\d(\Up_i\circ t):T^*K_i\ra T^*P_i^t$, and $\al\mapsto t^{-2}\al$
is multiplication by $t^{-2}$ in the vector space fibres
of~$T^*K_i\ra K_i$.

As in \eq{cu3eq5}--\eq{cu3eq6} and \eq{cu4eq4}--\eq{cu4eq5}, define
$U_\sNt\cap\pi^*(Q_i^t)$ and $\Phi_\sNt\vert_{U_\sNt\cap\pi^*(Q_i^t)}$~by
\begin{gather}
(\d\Xi_i^t)^*(U_\sNt)=\bigl\{(\si,r,\varsigma,u)
\in T^*\bigl(\Si_i\t(tT,R')\bigr):\bmd{(\varsigma,u)}<\ze r\bigr\}
\quad\text{and}
\label{cu6eq23}\\
\Phi_\sNt\circ\d\Xi_i^t(\si,r,\varsigma,u)\equiv\Up_i\circ\Phi_\sCi
\bigl(\si,r,\varsigma+\xi_i^{1,t}(\si,r),u+\xi_i^{2,t}(\si,r)\bigr)
\label{cu6eq24}
\end{gather}
for all $(\si,r,\varsigma,u)\in T^*\bigl(\Si_i\t(tT,R')\bigr)$
with $\bmd{(\varsigma,u)}<\ze r$, computing $\na,\md{\,.\,}$
using $\iota_i^*(g')$. Here $\Xi_i^t:\Si_i\t(tT,R')\ra Q_i^t$
is a diffeomorphism, and induces an isomorphism~$\d\Xi_i^t:
T^*\bigl(\Si_i\t(tT,R')\bigr)\ra T^*Q_i^t$.

Careful consideration shows that $U_\sNt$ is well-defined, and
$\Phi_\sNt$ is well-defined in \eq{cu6eq22} and \eq{cu6eq24}
for small $t$, so making $\de>0$ smaller if necessary $\Phi_\sNt$
is well-defined for $t\in(0,\de)$. Clearly $\Phi_\sNt$ is smooth
on each of $U_\sNt\cap\pi^*(K)$, $U_\sNt\cap\pi^*(P_i^t)$ and
$U_\sNt\cap\pi^*(Q_i^t)$, but we must still show that $\Phi_\sNt$
is smooth over the joins between them.

As $\xi_i^t\equiv\eta_i$ on $\Si_i\t[2t^\tau,R')$, comparing
\eq{cu3eq5}--\eq{cu3eq6} and \eq{cu6eq23}--\eq{cu6eq24} shows
that over $\Xi_i^t\bigl(\Si_i\t[2t^\tau,R')\bigr)\subset N^t
\cap X'$ we have $U_\sNt=U_\sXp$ and $\Phi_\sNt=\Phi_\sXp$.
Therefore $U_\sNt$ and $\Phi_\sNt$ join smoothly over the
$r=R'$ end of $Q_i^t$ and the corresponding end of $K$.
Similarly, as $\xi_i^t(\si,r)=t^2\chi_i(\si,t^{-1}r)$ when
$r\le t^\tau$, comparing \eq{cu4eq4}--\eq{cu4eq5} and
\eq{cu6eq23}--\eq{cu6eq24} shows that $U_\sNt=\d(\Up_i\circ
t)(t^2U_\sLi)$ and $\Phi_\sNt\circ\d(\Up_i\circ t)=\Up_i\circ
t\circ\Phi_\sLi\circ t^{-2}$ over $\Xi_i^t\bigl(\Si_i\t(tT,
t^\tau]\bigr)$. So $U_\sNt,\Phi_\sNt$ join smoothly at the
$r=tT$ end of $Q_i^t$ and $\pd P_i^t$, by~\eq{cu6eq22}.

Therefore $U_\sNt$ is an open tubular neighbourhood of $N^t$
in $T^*N^t$, and $\Phi_\sNt:U_\sNt\ra M$ is well-defined and
smooth. We shall show that $\Phi_\sNt^*(\om)=\hat\om$. On
$U_\sNt\cap\pi^*(K)$ this follows from $\Phi_\sXp^*(\om)=\hat\om$.
On $U_\sNt\cap\pi^*(P_i^t)$ it follows from $\Up_i^*(\om)=\om'$,
$\Phi_\sLi^*(\om')=\hat\om$, and the fact that the powers of $t$
in \eq{cu6eq22} cancel out in their effect on $\Phi_\sNt^*(\om)$.
On $U_\sNt\cap\pi^*(Q_i^t)$ it follows from $\Up_i^*(\om)=\om'$,
$\Phi_\sCi^*(\om')=\hat\om$ and the fact that $\xi_i^t$ is closed.

Define an $m$-form $\be^t$ on $U_\sNt$ by $\be^t=\Phi_\sNt^*(\Im\Om)$,
as in Definition~\ref{cu5def1}.
\label{cu6def4}
\end{dfn}

We now prove that parts (ii)--(v) of Theorem \ref{cu5thm2} hold for
$N^t$ when~$t\in(0,\de)$.

\begin{thm} Making $\de>0$ smaller if necessary, there exist\/
$A_1,A_3,\ldots,A_6>0$ such that for all\/ $t\in(0,\de)$,
as in parts {\rm(ii)--(v)} of Theorem \ref{cu5thm2} we have
\begin{itemize}
\item[{\rm(ii)}] $\psi\ge A_3$ on $N^t$.
\item[{\rm(iii)}] The subset\/ $\B_{A_1t}\subset T^*N^t$
of Definition \ref{cu5def1} lies in $U_\sNt$, and
$\cnm{\hat\na^k\be^t}{0}\le A_4t^{-k}$ on $\B_{A_1t}$ for
$k=0,1,2$ and\/~$3$.
\item[{\rm(iv)}] The injectivity radius $\de(h^t)$
satisfies~$\de(h^t)\ge A_5t$.
\item[{\rm(v)}] The Riemann curvature $R(h^t)$
satisfies~$\cnm{R(h^t)}{0}\le A_6t^{-2}$.
\end{itemize}
Here part\/ {\rm(iii)} uses the notation of Definition \ref{cu5def1},
and parts {\rm(iv)} and\/ {\rm(v)} refer to the compact Riemannian
manifold\/~$(N^t,h^t)$.
\label{cu6thm2}
\end{thm}

\begin{proof} All of (ii)--(v) are in fact elementary, following
from obvious facts about the behaviour of $N^t,h^t,U_\sNt,\Phi_\sNt$
for small $t$. Let $A_3=\inf_M\psi$. Then $A_3>0$ as $M$ is compact
and $\psi:M\ra(0,\iy)$ is continuous, and $\psi\ge A_3$ on $N^t$ for
all $t\in(0,\de)$ as $N^t\subset M$. This proves~(ii).

For part (iii), under $\Up_i\circ t:K_i\ra P_i^t$ we have
$(\Up_i\circ t)^*(U_\sNt)=t^2U_\sLi$. Hence the radii of the fibres
of $\pi:(\Up_i\circ t)^*(U_\sNt)\ra K_i$ scale like $t^2$ in the
metric $g'$ on $K_i$, and thus approximately like $t$ in the
metric $(\Up_i\circ t)^*(h^t)\approx t^2g'$ on $K_i$. Therefore
$U_\sNt$ contains all the balls of radius $A_1t$ in $T^*P_i^t$
for small~$A_1>0$.

The fibre of $\pi:U_\sNt\cap\pi^*(Q_i^t)\ra Q_i^t$ over
$\Xi_i^t(\si,r)$ for $(\si,r)\in\Si_i\t(tT,R')$ is a ball of radius
$\ze r$ about 0 in $T^*_{\smash{\sst\Xi_i^t(\si,r)}}Q_i^t$ w.r.t.\
the metric $(\Xi_i^t)_*\bigl(\iota_i^*(g')\bigr)$ on $Q_i^t$. As
$h^t$ and $(\Xi_i^t)_*\bigl(\iota_i^*(g')\bigr)$ are equivalent 
metrics uniformly in $t$, we see that $U_\sNt$ contains all the
balls of radius $A_1t$ in $T^*Q_i^t$ for small $A_1>0$. And $U_\sNt$
and $h^t$ are independent of $t$ over $K$, so this is obviously true
over $K$. Thus for small $A_1>0$ we have $\B_{A_1t}\subset U_\sNt$
for all~$t\in(0,\de)$.

To see that $\cnm{\hat\na^k\be^t}{0}\le A_4t^{-k}$ for $k=0,\ldots,3$
in (iii), that $\de(h^t)\ge A_5t$ in (iv) and $\cnm{R(h^t)}{0}\le A_6
t^{-2}$ in (v), consider the r\^ole of $t$ in defining $N^t,h^t,U_\sNt,
\Phi_\sNt$ and $\be^t$ in Definitions \ref{cu6def2}, \ref{cu6def3} and
\ref{cu6def4}. We make $N^t$ by shrinking $L_i$ by a factor $t>0$
and gluing it into $X'$ using $\Up_i$. Therefore for small $t>0$ the
metric $h^t$ on $P_i^t$, and on $Q_i^t$ near $P_i^t$, approximates
the metric $t^2g'$ on $L_i$, and the metric $h^t$ on $K$, and on
$Q_i^t$ near $K$, is $g\vert_{X'}$ and independent of~$t$.

Now under the homothety $g'\mapsto t^2g'$ for $g'$ on $L_i$ we have
$\de(t^2g')=t\,\de(g')$ and $\cnm{R(t^2g')}{0}=t^{-2}\cnm{R(g')}{0}$.
Thus on and near $P_i^t$ we have $\de(h^t)=O(t)$ and $\cnm{R(h^t)}{0}
=O(t^{-2})$, and it is easy to see that these are the dominant
contributions to $\de(h^t),\cnm{R(h^t)}{0}$. (In particular, the
derivatives of $F(t^{-\tau}r)$ on the annuli $\Si_i\t(t^\tau,2t^\tau)$
in \eq{cu6eq2} only contribute terms $O(t^\tau),O(t^{-2\tau})$
respectively.) Thus making $\de$ smaller if necessary, there exist
$A_5,A_6>0$ such that (iv), (v) hold for all~$t\in(0,\de)$.

In a similar way, on and near $P_i^t$ we can identify
$U_\sNt$ with $U_\sLi$ by construction, and then for small
$t\in(0,\de)$ we have $\smash{\hat h^t\approx t^2\hat h}$ and
$\be^t\approx t^m\Phi_\sLi^*(\Im\Om')$, where $\hat h^t$ and
$\hat h$ are the metrics constructed on $T^*N^t$ using $h^t$ and
on $T^*L_i$ using $g'\vert_{L_i}$ in Definition \ref{cu5def1}.
It then follows that $\cnm{\hat\na^k\be^t}{0}=O(t^{-k})$ on and
near $P_i^t$ for small $t$ and all $k\ge 0$. This is the
dominant contribution to $\cnm{\hat\na^k\be^t}{0}$ on $N^t$,
so making $\de$ smaller if necessary, part (iii) holds for
some~$A_4>0$.
\end{proof}

\subsection{Sobolev embedding inequalities on $N^t$}
\label{cu64}

We now prove that parts (vi) and (vii) of Theorem \ref{cu5thm2}
hold in the simplest case that $X'$ is connected. The author
learned the idea behind the next three results, and the
reference \cite{MiSi}, from Lee~\cite[\S 3]{Lee}.

Write $C^k_{\rm cs}(S)$ for the vector subspace of
{\it compactly-supported\/} functions in $C^k(S)$. In
\cite{MiSi}, Michael and Simon prove a Sobolev inequality
for submanifolds $S$ of $\R^l$, depending on their {\it mean
curvature vector} $H$ in $\R^l$. Applying \cite[Th.~2.1]{MiSi}
with $h=\md{u}^{2(m-1)/(m-2)}$ and using H\"older's inequality,
we easily prove:

\begin{thm} Let\/ $S$ be an $m$-submanifold of\/ $\R^l$ for
$m>2$, and\/ $u\in C^1_{\rm cs}(S)$. Then $\lnm{u}{2m/(m-2)}\le
D_1\bigl(\lnm{\d u}{2}+\lnm{u\,H}{2}\bigr)$, where $D_1>0$ depends
only on $m$, and\/ $H$ is the mean curvature of\/ $S$ in~$\R^l$.
\label{cu6thm3}
\end{thm}

If $S$ is {\it minimal\/} in $\R^n$ then $H\equiv 0$, and
$\lnm{u}{2m/(m-2)}\le D_1\lnm{\d u}{2}$. The SL $m$-folds
$tL_i$ in $\C^m$ are automatically minimal, so we deduce:

\begin{cor} There exists $D_1>0$ such that\/ $\lnm{u}{2m/(m-2)}
\le D_1\lnm{\d u}{2}$ for all\/ $t>0$, $i=1,\ldots,n$ and\/~$u
\in C^1_{\rm cs}(tL_i)$.
\label{cu6cor1}
\end{cor}

Next we prove a similar inequality for $X'$, when it is connected.

\begin{prop} Suppose $X'$ is connected. Then there exists
$D_2>0$ such that for all\/ $v\in C^1_{\rm cs}(X')$ we have
\e
\lnm{v}{2m/(m-2)}\le D_2\bigl(\lnm{\d v}{2}+
\bmd{\ts\int_{X'}v\,\d V_g\,}\,\bigr).
\label{cu6eq25}
\e
\label{cu6prop3}
\end{prop}

\begin{proof} By the Nash Embedding Theorem we can choose an
isometric embedding of $(M,g)$ in some $\R^l$. Then $X'$ is
also isometrically embedded in $\R^l$. The mean curvature
of $X'$ in $M$ is zero, as $X'$ is minimal, and the mean
curvature of $M$ in $\R^l$ is bounded, as $M$ is compact.
Thus the mean curvature $H$ of $X'$ in $\R^l$ is bounded,
say $\md{H}\le D_3$. Applying Theorem \ref{cu6thm3} to
$X'$ in $\R^l$ then gives
\e
\lnm{u}{2m/(m-2)}\le D_1\bigl(\lnm{\d u}{2}+D_3\lnm{u}{2}\bigr)
\quad\text{for all $u\in C^1_{\rm cs}(X')$.}
\label{cu6eq26}
\e

By studying the eigenvalues of $\De$ on $X'$, we show in
\cite[Th.~2.17]{Joyc4} that if $X'$ is connected then
for some $D_4>0$ and all $u\in C^2_{\rm cs}(X')$ with
$\int_{X'}u\,\d V_g=0$, we have $\lnm{u}{2}\le D_4\lnm{\d
u}{2}\le D_4^2\lnm{\De u}{2}$. As $C^2_{\rm cs}(X')$ is
dense in $C^1_{\rm cs}(X')$, the first inequality
$\lnm{u}{2}\le D_4\lnm{\d u}{2}$ holds if $u\in C^1_{\rm
cs}(X')$ with $\int_{X'}u\,\d V_g=0$. Combining this with
\eq{cu6eq26} and setting $D_5=D_1(1+D_3D_4)$ proves that
\e
\lnm{u}{2m/(m-2)}\le D_5\lnm{\d u}{2}
\quad\text{for all $u\in C^1_{\rm cs}(X')$ with
$\ts\int_{X'}u\,\d V_g=0$.}
\label{cu6eq27}
\e

Fix some $w\in C^1_{\rm cs}(X')$ with $\int_{X'}w\,\d V_g=1$.
For any $v\in C^1_{\rm cs}(X')$, define $u=v-w\int_{X'}v\,\d V_g$.
Then $u\in C^1_{\rm cs}(X')$ with $\int_{X'}u\,\d V_g=0$, so
\eq{cu6eq27} gives
\begin{align*}
\lnm{v}{2m/(m-2)}&-\bmd{\ts\int_{X'}v\,\d V_g\,}\cdot\lnm{w}{2m/(m-2)}
\le\lnm{u}{2m/(m-2)}\\
&\le D_5\lnm{\d u}{2}\le D_5\bigl(\lnm{\d v}{2}+
\bmd{\ts\int_{X'}v\,\d V_g\,}
\cdot\lnm{\d w}{2}\bigr).
\end{align*}
Equation \eq{cu6eq25} then follows with~$D_2=\max\bigl(D_5,
\lnm{w}{2m/(m-2)}\!+\!D_5\lnm{\d w}{2}\bigr)$.
\end{proof}

We now combine the inequalities on $L_i$ and $X'$ in Corollary
\ref{cu6cor1} and Proposition \ref{cu6prop3} to prove an
equality on $N^t$ for small $t$. The theorem implies parts
(vi) and (vii) of Theorem \ref{cu5thm2} with $W=\an{1}$, as
(vii) is trivial when~$W=\an{1}$.

\begin{thm} Suppose $X'$ is connected. Making $\de>0$ smaller if
necessary, there exists $A_7>0$ such that for all\/ $t\in(0,\de),$
if\/ $v\in L^2_1(N^t)$ with\/ $\int_{N^t}v\,\d V^t=0$ then $v\in
L^{2m/(m-2)}(N^t)$ and\/~$\lnm{v}{2m/(m-2)}\le A_7\lnm{\d v}{2}$.
\label{cu6thm4}
\end{thm}

\begin{proof} Choose $a,b\in\R$ with $0<a<b<\tau$. Then for
small $t>0$ we have
\e
2t^\tau<t^b<t^a<\min(1,R').
\label{cu6eq28}
\e
Making $\de>0$ smaller if necessary, we suppose that \eq{cu6eq28}
holds for all $t\in(0,\de)$. Let $G:(0,\iy)\ra[0,1]$ be a smooth,
decreasing function with $G(s)=1$ for $s\in(0,a]$ and $G(s)=0$
for $s\in[b,\iy)$. Write $G'$ for $\d G/\d s$. For $t\in(0,\de)$,
define a function $F^t:N^t\ra[0,1]$ by
\e
F^t(x)=\begin{cases} 1, & x\in K, \\
G\bigl((\log r)/(\log t)\bigr), & x=\Xi_i^t(\si,r)\in Q_i^t,
\quad i=1,\ldots,n, \\
0, & x\in P_i^t,\quad i=1,\ldots,n.
\end{cases}
\label{cu6eq29}
\e

Then $F^t$ is smooth, and $F^t\equiv 0$ on $P_i^t$ and $\Xi_i^t
\bigl(\Si_i\t(tT,t^b]\bigr)$ for $i=1,\ldots,n$, and $F^t\equiv 1$
on $K$ and $\Xi_i^t\bigl(\Si_i\t[t^a,R')\bigr)$ for $i=1,\ldots,n$.
It changes only on the annuli $\Xi_i^t\bigl(\Si_i\t(t^b,t^a)\bigr)$
for $i=1,\ldots,n$, and there we have
\e
(\Xi_i^t)^*(\d F^t)=(\log t)^{-1}\cdot G'\bigl((\log r)/(\log t)
\bigr)r^{-1}\d r \quad\text{for $r\in(t^b,t^a)$.}
\label{cu6eq30}
\e

Suppose now that $t\in(0,\de)$ and $v\in C^1(N^t)$ with $\int_{N^t}v\,
\d V^t=0$. The main idea of the proof is to write $v=F^tv+(1-F^t)v$,
where we treat $F^tv$ as a compactly-supported function on $X'$ and
apply Proposition \ref{cu6prop3} to it, and we treat $(1-F^t)v$ as
a sum over $i=1,\ldots,n$ of compactly-supported functions on $tL_i$,
and apply Corollary \ref{cu6cor1} to them.

The first is straightforward. As $2t^\tau<t^b$ by \eq{cu6eq28},
the support of $F^t$ lies in the union of $K$ and $\Xi_i^t\bigl(
\Si_i\t[2t^\tau,R')\bigr)$ for $i=1,\ldots,n$, which is part of
both $N^t$ and $X'$. Thus, extending $F^tv$ by zero to the rest
of $X'$ we can regard it as an element of $C^1_{\rm cs}(X')$, 
and since $\int_{X'}F^tv\,\d V_g=\int_{N^t}F^tv\,\d V^t$ equation 
\eq{cu6eq25} gives
\e
\begin{split}
\blnm{&F^tv}{2m/(m-2)}\le D_2\bigl(\blnm{\,\d(F^tv)}{2}+
\bmd{\ts\int_{N^t}F^tv\,\d V^t\,}\,\bigr)\\
&\le D_2\bigl(\lnm{F^t\d v}{2}+\lnm{v}{2m/(m-2)}\cdot
(\lnm{\d F^t}{m}+\lnm{1-F^t}{2m/(m+2)})\bigr).
\end{split}
\label{cu6eq31}
\e
Here all functions are on $N^t$ and all norms computed using $h^t$.
This is valid because the metrics $g$ on $X'$ and $h^t$ on $N^t$
coincide on the support of $F^t$. We have also used H\"older's
inequality upon $\lnm{v\,\d F^t}{2}$, and
\begin{equation*}
\ts\bmd{\int_{N^t}F^tv\,\d V^t}=\bmd{\int_{N^t}(1-F^t)v\,\d V^t}\le
\lnm{v}{2m/(m-2)}\cdot\lnm{1-F^t}{2m/(m+2)},
\end{equation*}
since $\d V_g=\d V^t$ on the support of $F^t$ as $g=h^t$,
and~$\int_{N^t}v\,\d V^t=0$.

For the second, identify the subset $tK_i\cup t\circ\vp_i\bigl(
\Si_i\t(T,t^{-1}R')\bigr)$ in $tL_i\subset\C^m$ with the subset
$P_i^t\cup Q_i^t$ in $N^t\subset M$ for $i=1,\ldots,n$ by
\begin{gather*}
tK_i\ni tx\mapsto\Up_i(tx)\in P_i^t\quad\text{and}\\
t\circ\vp_i\bigl(\Si_i\t(T,t^{-1}R')\bigr)\ni
t\circ\vp_i(\si,r)\mapsto\Xi_i^t(\si,tr)\in Q_i^t.
\end{gather*}
The restriction of $(1-F^t)v$ to $P_i^t\cup Q_i^t$ is
compactly-supported, so using this identification we can
regard $(1-F^t)v$ as an element of $C^1_{\rm cs}(tL_i)$,
extended by zero outside~$tK_i\cup t\circ\vp_i\bigl(
\Si_i\t(T,t^{-1}R')\bigr)$.

Under this identification between subsets of $tL_i$ and $N^t$, on
the support $tK_i\cup t\circ\vp_i\bigl(\Si_i\t(T,t^{a-1})\bigr)$
of $(1-F^t)v$ in $tL_i$, the metrics $g'$ on $tL_i$ and $h^t$ on
$N^t$ are close when $t$ is small. Applying Corollary \ref{cu6cor1}
to $(1-F^t)v$ on $tL_i$ gives $\lnm{(1-F^t)v}{2m/(m-2)}\le D_1
\lnm{\d((1-F^t)v)}{2}$, where norms are computed using $g'$ on
$tL_i$. As $g',h^t$ are close for small $t$, increasing $D_1$ to
$2D_1$ the same inequality holds with norms computed using $h^t$,
for small $t$. So making $\de>0$ smaller if necessary, for all
$t\in(0,\de)$ and $i=1,\ldots,n$ we have 
\e
\blnm{(1-F^t)v\vert_{P_i^t\cup Q_i^t}}{2m/(m-2)}\le
2D_1\blnm{\,\d((1-F^t)v)\vert_{P_i^t\cup Q_i^t}}{2}.
\label{cu6eq32}
\e

Now $(1-F^t)v$ is supported on $\bigcup_{i=1}^n(P_i^t\cup Q_i^t)$.
Therefore
\e
\begin{gathered}
\blnm{(1-F^t)v}{2m/(m-2)}\le\sum_{i=1}^n\blnm{(1-F^t)v
\vert_{P_i^t\cup Q_i^t}}{2m/(m-2)}
\quad\text{and}\\
\sum_{i=1}^n\blnm{\,\d((1-F^t)v)\vert_{P_i^t\cup Q_i^t}}{2}\le
\sqrt{n}\,\blnm{\,\d((1-F^t)v)}{2},
\end{gathered}
\label{cu6eq33}
\e
proving the second equation using $a_1+\cdots+a_n\le\sqrt{n}\,
(a_1^2+\cdots+a_n^2)^{1/2}$, which gives the factor $\sqrt{n}$.
Equations \eq{cu6eq32} and \eq{cu6eq33} give
\e
\begin{split}
\blnm{(1&-F^t)v}{2m/(m-2)}\le 2\sqrt{n}\,D_1\blnm{\,\d((1-F^t)v)}{2}\\
&\le 2\sqrt{n}\,D_1\bigl(\blnm{(1-F^t)\,\d v}{2}+
\lnm{v}{2m/(m-2)}\cdot\lnm{\d F^t}{m}\bigr).
\end{split}
\label{cu6eq34}
\e
Combining \eq{cu6eq31}, \eq{cu6eq34} and $\lnm{F^t\d v}{2},
\lnm{(1\!-\!F^t)\d v}{2}\le\lnm{\d v}{2}$ then proves
\e
\begin{split}
\bigl[1-(D_2+2\sqrt{n}\,D_1)&\lnm{\d F^t}{m}-D_2
\lnm{1-F^t}{2m/(m+2)}\bigr]\cdot\lnm{v}{2m/(m-2)}\\
&\le(D_2+2\sqrt{n}\,D_1)\lnm{\d v}{2}.
\end{split}
\label{cu6eq35}
\e

As the $L^m$ norm of $r^{-1}$ on $\Si_i\t(t^b,t^a)$ with the cone
metric is $O\bigl(\md{\log t}^{1/m}\bigr)$, using \eq{cu6eq30}
and $\vol\bigl(\supp(1-F^t)\bigr)=O(t^{ma}\bigr)$ we find that
\begin{equation*}
\lnm{\d F^t}{m}=O\bigl(\md{\log t}^{(1-m)/m}\bigr)
\quad\text{and}\quad
\lnm{1-F^t}{2m/(m+2)}=O(t^{a(m+2)/2})
\end{equation*}
for small $t$. Thus making $\de>0$ smaller if necessary we
can suppose that
\begin{equation*}
\bigl[1-(D_2+2\sqrt{n}\,D_1)\lnm{\d F^t}{m}-D_2
\lnm{1-F^t}{2m/(m+2)}\bigr]\ge\ha
\quad\text{for all $t\in(0,\de)$.}
\end{equation*}
Setting $A_7=2(D_2+2\sqrt{n}\,D_1)$, we see from \eq{cu6eq35}
that for all $v\in C^1(N^t)$ with $\int_{N^t}v\,\d V^t=0$ we
have $\lnm{v}{2m/(m-2)}\le A_7\lnm{\d v}{2}$. Since $C^1(N^t)$ is
dense in $L^2_1(N^t)$ and $L^2_1(N^t)\hookra L^{2m/(m-2)}(N^t)$
by the Sobolev Embedding Theorem, this completes the proof of
Theorem~\ref{cu6thm4}.
\end{proof}

\subsection{The main result, in the simplest case}
\label{cu65}

We can now prove our main result on desingularizations of
SL $m$-folds $X$ with conical singularities, in the simplest
case. To make the statement of the theorem short and easily
understood, we do not say much about the $\smash{\ti N^t}$.
But the construction gives a much more precise and detailed
description of the topology of $\smash{\ti N^t}$, and the topology
and geometry of its embedding in $(M,J,\om,\Om)$, which we can
read off from \S\ref{cu5}--\S\ref{cu6} above if we need to.

\begin{thm} Suppose $(M,J,\om,\Om)$ is an almost Calabi--Yau
$m$-fold and\/ $X$ a compact SL\/ $m$-fold in $M$ with conical
singularities at\/ $x_1,\ldots,x_n$ and cones $C_1,\ldots,C_n$.
Let\/ $L_1,\ldots,L_n$ be Asymptotically Conical SL\/ $m$-folds
in $\C^m$ with cones $C_1,\ldots,C_n$ and rates $\la_1,\ldots,
\la_n$. Suppose $\la_i<0$ for $i=1,\ldots,n$, and\/ $X'=X\sm
\{x_1,\ldots,x_n\}$ is connected.

Then there exists $\ep>0$ and a smooth family $\bigl\{
\smash{\ti N^t}:t\in(0,\ep]\bigr\}$ of compact, nonsingular
SL\/ $m$-folds in $(M,J,\om,\Om)$, such that\/ $\smash{\ti N^t}$
is constructed by gluing $tL_i$ into $X$ at\/ $x_i$ for
$i=1,\ldots,n$. In the sense of currents, $\smash{\ti N^t}\ra
X$ as~$t\ra 0$.
\label{cu6thm5}
\end{thm}

\begin{proof} By assumption $L_i$ is an AC SL $m$-fold with
rate $\la_i<0$. Now Theorem \ref{cu4thm2} shows that if $L$
is an AC SL $m$-fold in $\C^m$ with rate $\la<0$ then $L$ is
{\it also} Asymptotically Conical with rate $\la'$ for all
$\la'\in(2-m,0)$. It follows that if $\la_i\ge\ha(2-m)$ then
we can decrease the rate $\la_i$ of $L_i$ so that
$\la_i\in\bigl(2-m,\ha(2-m)\bigr)$. Therefore we can suppose
that $\la_i<\ha(2-m)$ for all $i=1,\ldots,n$, as in
Definition~\ref{cu6def1}.

Let $\de>0$ and $N^t$ for $t\in(0,\de)$ be as in Definition
\ref{cu6def2}, and make $\de$ smaller if necessary so that
Theorems \ref{cu6thm1}, \ref{cu6thm2} and \ref{cu6thm4} apply.
For each $t\in(0,\de)$, define a finite-dimensional vector
subspace $W^t\subset C^\iy(N^t)$ by $W^t=\an{1}$, the constant
functions. This will be $W\subset C^\iy(N)$ in Definition
\ref{cu5def1}. The natural projection $\pi_\sWt:L^2(N^t)\ra W^t$
is given by~$\pi_\sWt(v)=\vol(N^t)^{-1}\int_{N^t}v\,\d V^t$.

Let $\ve^t=\psi^m\sin\th^t\in C^\iy(N^t)$ for $t\in(0,\de)$
be as in Definition \ref{cu6def3}. Then $\Im\Om\vert_{N^t}=
\psi^m\sin\th^t\,\d V^t$. Thus
\begin{equation*}
\int_{N^t}\psi^m\sin\th^t\,\d V^t=\int_{N^t}\Im\Om=
[\Im\Om]\cdot[N^t]=[\Im\Om]\cdot[X]=0,
\end{equation*}
where $[\Im\Om]\in H^m(M,\R)$ and $[N^t],[X]\in H_m(M,\R)$, as
$N^t,X$ are homologous in $M$ and $\Im\Om\vert_{X'}\equiv 0$
as $X'$ is an SL $m$-fold. This implies that
\e
\pi_\sWt(\psi^m\sin\th^t)=\vol(N^t)^{-1}\int_{N^t}
\psi^m\sin\th^t\,\d V^t=0
\quad\text{for all $t\in(0,\de)$.}
\label{cu6eq36}
\e

Theorem \ref{cu6thm1} gives constants $\ka>1$ and $A_2>0$ such that
\e
\begin{gathered}
\lnm{\psi^m\sin\th^t}{2m/(m+2)}\le A_2t^{\ka+m/2},\quad
\cnm{\psi^m\sin\th^t}{0}\le A_2t^{\ka-1} \quad\text{and}\\
\lnm{\d(\psi^m\sin\th^t)}{2m}\le A_2t^{\ka-3/2}
\quad\text{for all $t\in(0,\de)$.}
\end{gathered}
\label{cu6eq37}
\e
Equations \eq{cu6eq36} and \eq{cu6eq37} imply that part (i) of
Theorem \ref{cu5thm2} holds for $N^t$ for all $t\in(0,\de)$,
replacing $N,W,\th$ by $N^t,W^t,\th^t$ respectively.

Let the Lagrangian neighbourhood $\Phi_\sNt:U_\sNt\ra M$ and the
$m$-form $\be^t$ on $U_\sNt$ be as in Definition \ref{cu6def4}.
Then Theorem \ref{cu6thm2} gives constants $A_1,A_3,\ldots,A_6>0$
such that parts (ii)--(v) of Theorem \ref{cu5thm2} hold for $N^t$
for all $t\in(0,\de)$, replacing $N,\be,h$ by $N^t,\be^t,h^t$
respectively.

As $\pi_\sWt(v)=0$ if and only if $\int_{N^t}v\,\d V^t=0$, Theorem
\ref{cu6thm4} gives $A_7>0$ such that if $v\in L^2_1(N^t)$ with
$\pi_\sWt(v)=0$ then $v\in L^{2m/(m-2)}(N^t)$ and $\lnm{v}{2m/(m-2)}
\le A_7\lnm{\d v}{2}$, for all $t\in(0,\de)$. Thus part (vi) of
Theorem \ref{cu5thm2} holds for $N^t$ for all $t\in(0,\de)$,
replacing $N,W$ by $N^t,W^t$ respectively. As $W^t=\an{1}$ part
(vii) of Theorem \ref{cu5thm2} is trivial for $N^t,W^t$, since
$\d^*\d w=\d w=0$ for all $w\in W^t$, and if $w\in W^t$ with
$\int_{N^t}w\,\d V^t=0$ then $w=0$. Thus part (vii) holds for
any $A_8>0$, and we take~$A_8=1$.

We have not yet showed that the inequality $\cos\th^t\ge\ha$
in Definition \ref{cu5def1} holds. From parts (i) and (ii) of
Theorem \ref{cu5thm2} we see that $\md{\sin\th^t}\le A_2A_3^{-m}
t^{\ka-1}$ on $N^t$. Thus for small $t\in(0,\de)$ we have
$\md{\sin\th^t}\le\frac{\sqrt{3}}{2}$ as $\ka>1$, so that
$\md{\cos\th^t}\ge\ha$. But $\cos\th^t$ is continuous, $N^t$
is connected, and $\cos\th^t\equiv 1$ on $K$ as $K$ is special
Lagrangian, so we must have $\cos\th^t\ge\ha$ on $N^t$ for
small~$t\in(0,\de)$.

We have constructed $\ka>1$ and $A_1,\ldots,A_8>0$ such that
parts (i)--(vii) of Theorem \ref{cu5thm2} hold for $N^t$ in
$M$ for all $t\in(0,\de)$. Let $\ep,K>0$ be as given in
Theorem \ref{cu5thm2} depending on $\ka,A_1,\ldots,A_8$ and
$m$, and make $\ep>0$ smaller if necessary to ensure that
$\ep<\de$ and $\cos\th^t\ge\ha$ on $N^t$ for $t\in(0,\ep]$.
Then Theorem \ref{cu5thm2} shows that for all $t\in(0,\ep]$
we can deform $N^t$ to a nearby compact, nonsingular SL
$m$-fold $\smash{\ti N^t}$, given by $\smash{\ti N^t}=
(\Phi_\sNt)_*\bigl(\Ga(\d f^t)\bigr)$ for some $f^t\in
C^\iy(N^t)$ with~$\cnm{\d f^t}{0}\le Kt^\ka<A_1t$.

Since $N^t$ and $\Phi_\sNt$ depend smoothly on $t$, we see that
$f^t$ is the locally unique solution of a nonlinear elliptic
p.d.e.\ on $N^t$ depending smoothly on $t$. It quickly follows
from general theory that $f^t$ depends smoothly on $t$, and so
$\smash{\ti N^t}$ does. One can easily show from Definition
\ref{cu6def2} that $N^t\ra X$ as currents as $t\ra 0$. But
the estimates on $\d f^t$ and $\na\d f^t$ in \S\ref{cu5}
imply that $\smash{\ti N^t}-N^t\ra 0$ as currents as $t\ra 0$,
so $\smash{\ti N^t}\ra X$ as $t\ra 0$. This completes the
proof of Theorem~\ref{cu6thm5}.
\end{proof}

\section{Desingularizing when $X'$ is not connected}
\label{cu7}

We now generalize the material of \S\ref{cu6} to the case
when $X'$ is not connected. Suppose $X'$ has $q>1$ connected
components $X'_1,\ldots,X'_q$. Then Theorem \ref{cu6thm4} no
longer holds. The basic reason for this is that the Laplacian
$\De^t$ on $N^t$ has $q$ small eigenvalues $0=\la_1^t,\ldots,
\la_q^t$, with~$\la_k^t=O(t^{m-2})$.

The corresponding eigenfunctions $1=v_1^t,\ldots,v_k^t$ are
approximately constant on the part of $N^t$ coming from $X_k'$
for each $k=1,\ldots,q$, and change on the part of $N^t$ coming
from $tL_i$ for $i=1,\ldots,n$. If $q>1$, calculation shows that
the $v_k^t$ for $1<k\le q$ satisfy $\int_{N^t}v_k^t\d V^t=0$ and
$\lnm{v_k^t}{2m/(m-2)}=O(t^{-(m-2)/2})\cdot\lnm{\d v_k^t}{2}$,
so that Theorem \ref{cu6thm4} cannot hold.

To get round this, we define in \S\ref{cu71} a vector subspace
$W^t\subset C^\iy(N^t)$ which approximates $\an{v_1^t,\ldots,v_q^t}$.
This will be $W$ in Theorem \ref{cu5thm2}. In \S\ref{cu73} we show
that if $v\in L^2_1(N^t)$ is $L^2$-orthogonal to $W^t$ then
$\lnm{v}{2m/(m-2)}\le A_7\lnm{\d v}{2}$ for some $A_7>0$ independent
of $t$, and this replaces Theorem~\ref{cu6thm4}.

There is also some more work to do. Part (i) of Theorem
\ref{cu5thm2} requires that $\lnm{\pi_\sWt(\psi^m\sin\th^t)}{1}
\le A_2t^{\ka+m-1}$. In \S\ref{cu6} this was trivial, as there
$W^t=\an{1}$ and $\pi_\sWt(\psi^m\sin\th^t)\equiv 0$ for
topological reasons. We shall show in \S\ref{cu72} that here
for $\pi_\sWt(\psi^m\sin\th^t)$ to be sufficiently small the
invariants $Z(L_i)$ for $i=1,\ldots,n$ must satisfy equation
\eq{cu7eq1} below. The final results are given in~\S\ref{cu74}.

\subsection{Setting up the problem}
\label{cu71}

We shall consider the following situation, the analogue of
Definitions~\ref{cu6def1}--\ref{cu6def3}.

\begin{dfn} We work in the situation of Definition \ref{cu6def1}.
Thus, $(M,J,\om,\Om)$ is an almost Calabi--Yau $m$-fold for
$m>2$ with metric $g$, $\psi:M\ra(0,\iy)$ satisfies \eq{cu2eq3},
$X$ is a compact SL $m$-fold in $M$ with conical singularities
at $x_1,\ldots,x_n$ with identifications $\up_i$, cones $C_i$
and rates $\mu_i$, and $L_1,\ldots,L_n$ are AC SL $m$-folds
in $\C^m$ with cones $C_i$ and rates $\la_i<0$. We write $X'=
X\sm\{x_1,\ldots,x_n\}$ and $\Si_i=C_i\cap{\cal S}^{2m-1}$,
and use the other notation of~\S\ref{cu6}.

However, we do {\it not\/} assume that $X'$ is connected, as we did
in \S\ref{cu64}. Set $q=b^0(X')$, so that $X'$ has $q$ connected
components, and number them $X_1',\ldots,X_q'$. For $i=1,\ldots,n$
let $l_i=b^0(\Si_i)$, so that $\Si_i$ has $l_i$ connected
components, and number them~$\Si_i^1,\ldots,\Si_i^{\smash{l_i}}$.

Now $\Up_i\circ\vp_i$ is a diffeomorphism $\Si_i\t(0,R')\ra S_i
\subset X'$. For each $j=1,\ldots,l_i$, $\Up_i\circ\vp_i\bigl(
\Si_i^j\t(0,R')\bigr)$ is a connected subset of $X'$, and so
lies in exactly one of the $X_k'$ for $k=1,\ldots,q$. Define
numbers $k(i,j)=1,\ldots,q$ for $i=1,\ldots,n$ and $j=1,\ldots,l_i$
by $\Up_i\circ\vp_i\bigl(\Si_i^j\t(0,R')\bigr)\subset
X'_{\smash{k(i,j)}}$. Suppose that
\e
\sum_{\substack{1\le i\le n, \; 1\le j\le l_i: \\
k(i,j)=k}}\psi(x_i)^mZ(L_i)\cdot[\Si_i^j\,]=0
\quad\text{for all $k=1,\ldots,q$.}
\label{cu7eq1}
\e
Here $Z(L_i)\in H^{m-1}(\Si_i,\R)$ is as in Definition
\ref{cu4def2}, and $[\Si_i^j]\in H_{m-1}(\Si_i,\Z)$,
and `$\,\cdot\,$' is the contraction $H^{m-1}(\Si_i,\R)
\t H_{m-1}(\Si_i,\Z)\ra\R$. The reason for \eq{cu7eq1}
will appear in Proposition \ref{cu7prop3} below.

Define Lagrangian $m$-folds $N^t$ for $t\in(0,\de)$ as in
Definition \ref{cu6def2}, but when we choose $\tau\in(0,1)$
to satisfy \eq{cu6eq1} we also require that
\e
\ts\frac{m}{m+1}<\tau<1.
\label{cu7eq2}
\e
Clearly this is possible. Suppose that the topology of $X$ and
$L_i$ is such that the $N^t$ are {\it connected}. This requires
that $X$ be connected, but not that $X'$ or the $L_i$ be
connected. A sufficient (but not necessary) condition for the
$N^t$ to be connected is that $X$ and the $L_i$ are connected.

Let $h^t=g\vert_{N^t}$, so that $(N^t,h^t)$ is a Riemannian
manifold. As in Definition \ref{cu6def3} the $N^t$ are oriented,
and $\Om\vert_{N^t}=\psi^m{\rm e}^{\smash{i\th^t}}\,\d V^t$ for some
phase function ${\rm e}^{\smash{i\th^t}}$ on $N^t$, where $\d V^t$
is the volume form on $N^t$. We write $\ve^t=\psi^m\sin\th^t$,
so that $\Im\Om\vert_{N^t}=\ve^t\,\d V^t$ for~$t\in(0,\de)$.
\label{cu7def1}
\end{dfn}

To apply Theorem \ref{cu5thm2} to $N^t$ we need a vector
subspace $W^t\subset C^\iy(N^t)$, to be $W$ in Definition
\ref{cu5def1}. In \S\ref{cu6}, where $X'$ was connected,
we took $W^t=\an{1}$. However, when $X'$ is not connected
we must introduce a nontrivial $W^t$ with $\dim W^t=b^0(X')$
to repair the proof of Theorem \ref{cu6thm4}. Here is
the definition.

\begin{dfn} We work in the situation of Definition \ref{cu7def1}.
For $i=1,\ldots,n$ apply Theorem \ref{cu4thm4} to the AC SL
$m$-fold $L_i$ in $\C^m$, using the numbering $\Si_i^j$ chosen
in Definition \ref{cu7def1} for the connected components of
$\Si_i$. This gives a vector space $V_i$ of bounded harmonic
functions on $L_i$ with $\dim V_i=l_i$. For each ${\bf c}_i=
(c_i^1,\ldots,c_i^{l_i})\in\R^{l_i}$ there exists a unique
$v^{{\bf c}_i}_{\smash{i}}\in V_i$ with
\e
\na^k\bigl(\vp_i^*(v^{{\bf c}_i}_{\smash{i}})-c^j_i\,\bigr)=
O\bigl(\md{{\bf c}_i}r^{\be-k}\bigr)
\quad\text{on $\Si_i^j\t(T,\iy)$ as $r\ra\iy$,}
\label{cu7eq3}
\e
for all $i=1,\ldots,n$, $j=1,\ldots,l_i$, $k\ge 0$ and~$\be\in(2-m,0)$.

We shall define a vector subspace $W^t\subset C^\iy(N^t)$
for $t\in(0,\de)$, with an isomorphism $W^t\cong\R^q$. Fix
${\bf d}=(d_1,\ldots,d_q)\in\R^q$, and set $c_i^j=d_{k(i,j)}$
for $i=1,\ldots,n$ and $j=1,\ldots,l_i$. Let ${\bf c}_i=
(c_i^1,\ldots,c_i^{l_i})$. This defines vectors ${\bf c}_i\in
\R^{l_i}$ for $i=1,\ldots,n$, which depend linearly on $\bf d$.
Hence we have harmonic functions $v^{{\bf c}_i}_{\smash{i}}\in
V_i\subset C^\iy(L_i)$, which also depend linearly on~$\bf d$.

Let $F:(0,\iy)\ra[0,1]$ and $\tau\in(0,1)$ be as in Definition
\ref{cu6def2}. Make $\de>0$ smaller if necessary so that
$tT<\ha t^\tau$ for all $t\in(0,\de)$. For $t\in(0,\de)$,
define a function $w_{\bf d}^t\in C^\iy(N^t)$ as follows:
\begin{itemize}
\item[(i)] The subset $K\subset N^t$ has $q$ connected
components $K\cap X_k'$ for $k=1,\ldots,q$. Define
$w_{\bf d}^t\equiv d_k$ on $K\cap X_k'$ for~$k=1,\ldots,q$.
\item[(ii)] Define $w_{\bf d}^t$ on $P_i^t\subset N^t$ by
$(\Up_i\circ t\circ\vp_i)^*(w_{\bf d}^t)\equiv
v^{{\bf c}_i}_{\smash{i}}$ on $K_i$ for~$i=1,\ldots,n$.
\item[(iii)] Define $w_{\bf d}^t$ on $Q_i^t\subset N^t$ by
\e
(\Xi_i^t)^*(w_{\bf d}^t)(\si,r)=\bigl(1-F(2t^{-\tau}r)\bigr)
\vp_i^*(v^{{\bf c}_i}_{\smash{i}})(\si,t^{-1}r)+F(2t^{-\tau}r)c_i^j
\label{cu7eq4}
\e
on $\Si_i^j\t(tT,R')$, for $i=1,\ldots,n$ and~$j=1,\ldots,l_i$.
\end{itemize}
It is easy to see that $w_{\bf d}^t$ is smooth over the joins
between $P_i^t,Q_i^t$ and $K$, so $w_{\bf d}^t\in C^\iy(N^t)$. Also
$w_{\bf d}^t$ is linear in ${\bf d}$, as $v^{{\bf c}_i}_{\smash{i}}$
is. Thus $W^t=\{w_{\bf d}^t:{\bf d}\in\R^q\}$ is a vector subspace of
$C^\iy(N^t)$ isomorphic to $\R^q$, for all~$t\in(0,\de)$.

If ${\bf d}=(1,\ldots,1)$ then $c_i^j\equiv 1$, so ${\bf c}_i=(1,
\ldots,1)$ for $i=1,\ldots,n$, and thus $v^{{\bf c}_i}_{\smash{i}}
\equiv 1$ for $i=1,\ldots,n$ by Theorem \ref{cu4thm4}. Therefore
$w_{\smash{(1,\ldots,1)}}^{\smash{t}}\equiv 1$ by (i)--(iii)
above, and $1\in W^t$ for all $t\in(0,\de)$. This corresponds
to the condition $1\in W$ in Definition \ref{cu5def1}. Define
$\pi_\sWt:L^2(N^t)\ra W^t$ to be the projection onto
$W^t$ using the $L^2$-inner product, as for $\pi_{\sst W}$
in Definition~\ref{cu5def1}.
\label{cu7def2}
\end{dfn}

In Definition \ref{cu6def2} we defined the $N^t$ by gluing
together $X'$ and $tL_i$ using $F(t^{-\tau}r)$ in \eq{cu6eq2}.
In contrast, equation \eq{cu7eq4} uses $F(2t^{-\tau}r)$. As
$F$ increases from 0 to 1 on $[1,2]$, this means that in
Definition \ref{cu6def2} we glued $X'$ and $tL_i$ together
on the annuli $\Xi_i^t\bigl(\Si_i\t(t^\tau,2t^\tau)\bigr)$,
but in \eq{cu7eq4} we glue $v^{{\bf c}_i}_{\smash{i}}$ and
$c_i^j$ together on the annuli $\Xi_i^t\bigl(\Si_i^j\t(\ha
t^\tau,t^\tau)\bigr)$. This will be important in \S\ref{cu72},
where the fact that $w_{\bf d}^t$ is {\it constant\/} on the
annuli $\Xi_i^t\bigl(\Si_i\t(t^\tau,2t^\tau)\bigr)$ will
simplify the calculations.

The point of Definition \ref{cu7def2} is that $w_{\bf d}^t$ is
made by gluing together the constant function $d_k$ on $X_k'$
with the harmonic function $t_*(v^{{\bf c}_i}_{\smash{i}})$ on
$tL_i$. This is asymptotic to $d_{\smash{k(i,j)}}$ on the
$j^{\rm th}$ end of $tL_i$, which is glued into $X_{\smash{k(i,j)}}'$.
Thus, for small $t$ the functions $d_{\smash{k(i,j)}}$ and
$t_*(v^{{\bf c}_i}_{\smash{i}})$ are nearly equal in the annulus
$\Xi_i^t\bigl(\Si_i^j\t(\ha t^\tau,t^\tau)\bigr)$ where they are
glued together.

Thus we expect $w_{\bf d}^t$ to be {\it nearly harmonic}
on $N^t$ for small $t$, that is, $\d^*\d w_{\bf d}^t$ is small
compared to $w_{\bf d}^t$. The next proposition estimates
$w_{\bf d}^t,\d w_{\bf d}^t$ and~$\d^*\d w_{\bf d}^t$.

\begin{prop} In the situation of Definition \ref{cu7def2}, for all\/
$t\in(0,\de)$, ${\bf d}=(d_1,\ldots,d_q)\in\R^q$, $\be\in(2-m,0)$
and\/ $i=1,\ldots,n$ we have
\ea
\md{w_{\bf d}^t}&\le\max\bigl(\md{d_1},\ldots,\md{d_q}\bigr)
\quad\text{on $N^t$,}
\label{cu7eq5}\\
\d w_{\bf d}^t=\d^*\d w_{\bf d}^t&=0
\quad\text{on $K$ and\/ $\Xi_i^t\bigl(\Si_i\t[t^\tau,R')\bigr)$,}
\label{cu7eq6}\\
\bmd{(\Xi_i^t)^*(\d w_{\bf d}^t)}(\si,r)
&=O\bigl(\md{{\bf d}}t^{-\be}r^{\be-1}\bigr)
\quad\text{on $\Si_i\t(tT,t^\tau)$,}
\label{cu7eq7}\\
\bmd{(\Xi_i^t)^*(\d^*\d w_{\bf d}^t)}(\si,r)
&=O\bigl(\md{{\bf d}}t^{-\be}r^{\be-1}\bigr)
\quad\text{on $\Si_i\t(tT,\ha t^\tau]$,}
\label{cu7eq8}\\
\bmd{(\Xi_i^t)^*(\d^*\d w_{\bf d}^t)}(\si,r)
&=O\bigl(\md{{\bf d}}t^{-\be+\tau(\be-2)}\bigr)
\quad\text{on $\Si_i\t(\ha t^\tau,t^\tau)$,}
\label{cu7eq9}\\
\text{and}\quad
\md{\d w_{\bf d}^t}&=O\bigl(\md{{\bf d}}t^{-1}\bigr),\quad
\md{\d^*\d w_{\bf d}^t}=O\bigl(\md{{\bf d}}t^{-1}\bigr)
\quad\text{on $P_i^t$.}
\label{cu7eq10}
\ea
Here $\d^*$ and\/ $\md{\,.\,}$ are computed using $h^t$
or~$(\Xi_i^t)^*(h^t)$.
\label{cu7prop1}
\end{prop}

\begin{proof} Since $v^{{\bf c}_i}_{\smash{i}}$ is harmonic it
has no strict maxima or minima in $L_i$ by the maximum principle,
and as it approaches $d_{k(i,j)}$ on the $j^{\rm th}$ end we
see that $\bmd{v^{{\bf c}_i}_{\smash{i}}}\le\max\bigl(\md{d_1},
\ldots,\md{d_q}\bigr)$ on $L_i$. Equation \eq{cu7eq5} then follows
from the definition of $w_{\bf d}^t$ above. Also $w_{\bf d}^t$ is
constant on $K$ and $\Xi_i^t\bigl(\Si_i\t[t^\tau,R')\bigr)$,
so \eq{cu7eq6} holds.

Let $h_i$ be the metric $g'\vert_{L_i}$ on $L_i$. Then we have
\ea
(\Up_i\circ t)^*(h^t)&=t^2h_i+O(t^3)&\;&\text{on $K_i$, and}
\label{cu7eq11}\\
(\Up_i\circ t\circ\vp_i)^*(h^t)&=t^2\vp_i^*(h_i)+O(t^2r)
&\;&\text{on $\Si_i\t(tT,t^\tau)$.}
\label{cu7eq12}
\ea
Therefore the metrics $h^t$ on $P_i^t$ and $\Xi_i^t\bigl(\Si_i
\t(tT,t^\tau)\bigr)$ and $t^2h_i$ on the corresponding regions
of $L_i$ are uniformly equivalent for small $t$. Using this,
\eq{cu7eq3} and \eq{cu7eq4} we deduce \eq{cu7eq7} and the
first equation of~\eq{cu7eq10}.

By the same method we also find that
\ea
\md{\na^2w_{\bf d}^t}&=O\bigl(\md{{\bf d}}t^{-2}\bigr) \
&\;&\text{on $P_i^t$, and}
\label{cu7eq13}\\
\bmd{(\Xi_i^t)^*(\na^2w_{\bf d}^t)}(\si,r)
&=O\bigl(\md{{\bf d}}t^{-\be}r^{\be-2}\bigr)
&\;&\text{on $\Si_i\t(tT,t^\tau)$,}
\label{cu7eq14}
\ea
computing $\na$ and $\md{\,.\,}$ using $h^t$ or $(\Xi_i^t)^*(h^t)$.
Equation \eq{cu7eq9} then follows from \eq{cu7eq14} with $r\in
(\ha t^\tau,t^\tau)$, as~$\md{\d^*\d w_{\bf d}^t}\le
\md{\na^2w_{\bf d}^t}$.

Now from \eq{cu7eq11} and a similar equation for derivatives
we find that
\e
\d^*_{h^t}\d w_{\bf d}^t=\d^*_{t^2h_i}\d w_{\bf d}^t+
O(t)\cdot\md{\na^2w_{\bf d}^t}+O(1)\cdot\md{\d w_{\bf d}^t}
\quad\text{on $P_i^t$,}
\label{cu7eq15}
\e
where $\d^*_{\smash{h^t}},\d^*_{\smash{t^2h_i}}$ are $\d^*$
computed using $h^t,t^2h_i$. But $w_{\bf d}^t=(\Up_i\circ t)_*
(v^{{\bf c}_i}_{\smash{i}})$ on $P_i^t$, and $v^{{\bf c}_i}_{
\smash{i}}$ is harmonic w.r.t.\ $h_i$, and so w.r.t.\ $t^2h_i$.
Hence $\d^*_{t^2h_i}\d w_{\bf d}^t=0$ on $P_i^t$, and combining
\eq{cu7eq13}, \eq{cu7eq15} and the first equation of \eq{cu7eq10}
gives the second equation of \eq{cu7eq10}. We prove \eq{cu7eq8} in
the same way, using \eq{cu7eq7}, \eq{cu7eq12} and~\eq{cu7eq14}.
\end{proof}

Another way to think about $W^t$ is that the Laplacian $\De^t=
\d^*\d$ of $h^t$ on $N^t$ has $q$ {\it small eigenvalues}
$\la_1^t,\ldots,\la_q^t$, counted with multiplicity and including
0, and $W^t$ approximates the sum of the corresponding eigenspaces
of $\De^t$. From Proposition \ref{cu7prop1} one can show that
\begin{equation*}
\lnm{w_{\bf d}^t}{2}=O(\md{{\bf d}})\quad\text{and}\quad
\lnm{\d w_{\bf d}^t}{2}=O(\md{{\bf d}}t^{(m-2)/2}),
\end{equation*}
so that $\lnm{\d w_{\bf d}^t}{2}^2=O(t^{m-2})\cdot\lnm{w_{\bf d}^t}{2}^2$.
This implies that the dominant eigenvectors of $\De^t$ in $w_{\bf d}^t$
have eigenvalues $O(t^{m-2})$, so $\la_k^t=O(t^{m-2})$ as~$t\ra 0$.

\subsection{Part (i) of Theorem \ref{cu5thm2}: estimating
$\lnm{\pi_\sWt(\psi^m\sin\th^t)}{1}$}
\label{cu72}

We need to bound $\lnm{\pi_\sWt(\psi^m\sin\th^t)}{1}$
to prove part (i) of Theorem \ref{cu5thm2} for $N^t,W^t$.
Writing $\ve^t=\psi^m\sin\th^t$ as in Definition \ref{cu6def3},
we shall do this by first estimating $\int_{N^t}w_{\bf d}^t\ve^t
\,\d V^t$ for all ${\bf d}\in\R^q$. As $w_{\bf d}^t\equiv d_k$
on $K\cap X_k'$ and on $\Xi_i^t\bigl(\Si_i^j\t[t^\tau,R')\bigr)$
when $k(i,j)=k$, we see that
\e
\begin{split}
\int_{N^t}w_{\bf d}^t\ve^t\,\d V^t=&
\sum_{i=1}^n\biggl(\,\int_{P_i^t}w_{\bf d}^t\ve^t\,\d V^t+\int_{
\Xi_i^t(\Si_i\t(tT,t^\tau))}w_{\bf d}^t\ve^t\,\d V^t\biggr)\\
+&\sum_{k=1}^qd_k\biggl(\,\int_{K\cap X_k'}
\!\!\!\!\!\!\ve^t\,\d V^t
+\!\!\sum_{\substack{1\le i\le n, \; 1\le j\le l_i: \\
k(i,j)=k}}\int_{\Xi_i^t(\Si_i^j\t
[t^\tau,R'))}\!\!\!\!\!\!\!\!\!\!\!\!\!\ve^t\,\d V^t\biggr).
\end{split}
\label{cu7eq16}
\e
The next two propositions bound the bracketed terms on each line.

\begin{prop} For all\/ $t\in(0,\de)$, ${\bf d}\in\R^q$ and\/
$i=1,\ldots,n$ we have
\e
\int_{P_i^t}w_{\bf d}^t\ve^t\,\d V^t+\int_{\Xi_i^t(\Si_i\t(tT,t^\tau))}
w_{\bf d}^t\ve^t\,\d V^t=O\bigl(\md{{\bf d}}t^{(m+1)\tau}\bigr).
\label{cu7eq17}
\e
\label{cu7prop2}
\end{prop}

\begin{proof} We have $\md{w_{\bf d}^t}\le\md{{\bf d}}$ on $N^t$ by
\eq{cu7eq5}. On $P_i^t$ we have $\md{\ve^t}\le Ct$ by \eq{cu6eq7}.
As $\vol(P_i^t)=O(t^m)$, this implies that the first integral in
\eq{cu7eq17} is $O\bigl(\md{{\bf d}}t^{m+1}\bigr)$. Similarly, on
$\Si_i\t(tT,t^\tau)$ we have $\bmd{(\Xi_i^t)^*(\ve^t)}\le Cr\le
Ct^\tau$ by \eq{cu6eq5}. As $\vol\bigl(\Xi_i^t(\Si_i\t(tT,t^\tau))
\bigr)=O(t^{m\tau})$, this implies that the second integral in
\eq{cu7eq17} is $O\bigl(\md{{\bf d}}t^{(m+1)\tau}\bigr)$.
The proposition follows.
\end{proof}

\begin{prop} For all\/ $t\in(0,\de)$, ${\bf d}\in\R^q$ and\/
$k=1,\ldots,q$ we have
\e
\begin{split}
\int_{K\cap X_k'}\!\!\!\!\!\!\ve^t\,\d V^t
&+\!\!\sum_{\substack{1\le i\le n, \; 1\le j\le l_i: \\
k(i,j)=k}}\int_{\Xi_i^t(\Si_i^j\t
[t^\tau,R'))}\!\!\!\!\!\!\!\!\!\!\!\!\!\ve^t\,\d V^t=\\
&-t^m\!\!\!\!\!\!\!\!\sum_{\substack{1\le i\le n, \; 1\le j\le l_i: \\
k(i,j)=k}}\!\!\!\!\!\!\!\! \psi(x_i)^mZ(L_i)\cdot[\Si_i^j\,]
+O\bigl(t^{(m+1)\tau}\bigr),
\end{split}
\label{cu7eq18}
\e
where $Z(L_i)\in H^{m-1}(\Si_i,\R)$ is as in \S\ref{cu41},
and\/~$[\Si_i^j\,]\in H_{m-1}(\Si_i,\Z)$.
\label{cu7prop3}
\end{prop}

\begin{proof} As $\ve^t\d V^t=\Im\Om\vert_{N^t}$, the left
hand side of \eq{cu7eq18} is the integral of $\Im\Om$ over
the $m$-chain 
\begin{equation*}
Z_k=\bigl(K\cap X_k'\bigr)+
\sum_{\substack{1\le i\le n, \; 1\le j\le l_i: \\
k(i,j)=k}}\Xi_i^t\bigl(\Si_i^j\t[t^\tau,R')\bigr)
\end{equation*}
for $k=1,\ldots,q$, which is a closed subset of $N^t$, with
boundary $(m\!-\!1)$-chain
\e
\pd Z_k=
-\sum_{\substack{1\le i\le n, \; 1\le j\le l_i: \\
k(i,j)=k}}\Xi_i^t\bigl(\Si_i^j\t\{t^\tau\}\bigr).
\label{cu7eq19}
\e

For each $i=1,\ldots,n$ and $j=1,\ldots,l_i$, define an $m$-chain
$A_i^j$ in $B_R$ to be the image of $\Si_i^j\t[0,1]$
under the map $\Si_i^j\t[0,1]\ra B_R$ given by
\begin{equation*}
(\si,r)\longmapsto r\Phi_\sCi\bigl(\si,t^\tau,t^2\chi_i^1
(\si,t^{\tau-1}),t^2\chi_i^2(\si,t^{\tau-1})\bigr).
\end{equation*}
As $\Up_i\circ\Phi_\sCi\bigl(\si,t^\tau,t^2\chi_i^1(\si,t^{\tau-1}),
t^2\chi_i^2(\si,t^{\tau-1})\bigr)\equiv\Xi_i^t(\si,t^\tau)$ for $\si
\in\Si_i$ by Definition \ref{cu6def2}, we see that
\e
\pd\bigl(\Up_i(A_i^j)\bigr)=\Xi_i^t\bigl(\Si_i^j\t\{t^\tau\}\bigr),
\label{cu7eq20}
\e
regarding $\Up_i(A_i^j)$ as an $m$-chain in~$M$.

Now define another $m$-chain $Z_k'$ for $k=1,\ldots,q$ to be
\begin{equation*}
Z_k'=\,\ov{\!X_k'\!}\,-\sum_{\substack{1\le i\le n, \; 1\le j\le
l_i: \\ k(i,j)=k}}\Up_i(A_i^j).
\end{equation*}
As $\,\ov{\!X_k'\!}\,$ is an $m$-chain without boundary, we see
from \eq{cu7eq19} and \eq{cu7eq20} that $\pd Z_k'=\pd Z_k$, and
in fact it is easy to see that $Z_k'$ and $Z_k$ are homologous
in $M$. Since $\Im\Om$ is a closed $m$-form on $M$, this implies
that $\int_{Z_k'}\Im\Om=\int_{Z_k}\Im\Om$. But $\Im\Om\vert_{X_k'}
\equiv 0$ as $X'_k$ is special Lagrangian. Hence we see that
\e
\int_{Z_k}\Im\Om=-\sum_{\substack{1\le i\le n, \; 1\le j\le l_i: \\
k(i,j)=k}}\int_{A_i^j}\Up_i^*(\Im\Om).
\label{cu7eq21}
\e

From Definition \ref{cu3def3} we have $\up_i^*(\Om)=\psi(x_i)^m\Om'$,
where $\Om'$ is as in \eq{cu2eq1}. Thus as $\Up_i^*(\Om)$ is smooth
on $B_R$, Taylor's theorem gives
\e
\Up_i^*(\Im\Om)=\psi(x_i)^m\Im\Om'+O(r)\quad\text{on $B_R$.}
\label{cu7eq22}
\e
Now $A_i^j$ is an $m$-chain in $B_R\subset\C^m$ with boundary in
the AC SL $m$-fold $tL_i$, and $[\pd A_i^j]\in H_{m-1}(tL_i,\R)$
is the image of $[\Si_i^j]\in H_{m-1}(\Si_i,\R)$ under the map
$H_{m-1}(\Si_i,\R)\ra H_{m-1}(L_i,\R)$ dual to the natural map
$H^{m-1}(L_i,\R)\ra H^{m-1}(\Si_i,\R)$. It then follows easily
from Definition \ref{cu4def2} and Lemma \ref{cu4lem} that
\e
\int_{A_i^j}\Im\Om'=Z(tL_i)\cdot[\Si_i^j]=t^mZ(L_i)\cdot[\Si_i^j].
\label{cu7eq23}
\e

But as $r=O(t^\tau)$ on $A_i^j$ and $\vol(A_i^j)=O(t^{m\tau})$
we see from \eq{cu7eq22} that
\e
\int_{A_i^j}\bigl(\Up_i^*(\Im\Om)-\psi(x_i)^m\Im\Om'\bigr)
=O(t^{(m+1)\tau}).
\label{cu7eq24}
\e
Equation \eq{cu7eq18} now follows from \eq{cu7eq21}, \eq{cu7eq23},
\eq{cu7eq24}, and the fact that the left hand side of \eq{cu7eq18}
is $\int_{Z_k}\Im\Om$. This completes the proof.
\end{proof}

We can now explain the reason for the condition \eq{cu7eq1} in
Definition \ref{cu7def1}. If \eq{cu7eq1} holds then the first
term on the right hand side of \eq{cu7eq18} is zero, and therefore
\eq{cu7eq16} and Propositions \ref{cu7prop2} and \ref{cu7prop3}
show that $\int_{N^t}w_{\bf d}^t\ve^t\,\d V^t=O\bigl(\md{{\bf d}}
t^{(m+1)\tau}\bigr)$. This in turn implies that~$\lnm{\pi_{\sst
W^t}(\ve^t)}{1}=O(t^{(m+1)\tau})$.

Therefore $\lnm{\pi_\sWt(\psi^m\sin\th^t)}{1}\le A_2
t^{\ka+m-1}$ for some $A_2>0$ and $\ka>1$ and all $t\in(0,\de)$,
as $\ve^t=\psi^m\sin\th^t$ and $\tau>\frac{m}{m+1}$ by \eq{cu7eq2}.
This is one of the conditions in part (i) of Theorem \ref{cu5thm2}
for $N^t,W^t$. However, if \eq{cu7eq1} does not hold then
$\lnm{\pi_\sWt(\psi^m\sin\th^t)}{1}=O(t^m)$, and part (i)
of Theorem \ref{cu5thm2} for $N^t,W^t$ does not hold for all
$t\in(0,\de)$ with $\ka>1$, so the construction fails.

Here is the analogue of Theorem~\ref{cu6thm1}.

\begin{thm} Making $\de>0$ smaller if necessary, there exist\/
$A_2>0$ and\/ $\ka>1$ such that\/ $\ve^t=\psi^m\sin\th^t$ on
$N^t$ satisfies $\lnm{\ve^t}{2m/(m+2)}\le A_2t^{\ka+m/2}$,
$\cnm{\ve^t}{0}\le A_2t^{\ka-1}$,
$\lnm{\d\ve^t}{2m}\le A_2t^{\ka-3/2}$ and\/
$\lnm{\pi_\sWt(\ve^t)}{1}\le A_2t^{\ka+m-1}$ for
all\/ $t\in(0,\de)$, as in part\/ {\rm(i)} of Theorem~\ref{cu5thm2}.
\label{cu7thm1}
\end{thm}

\begin{proof} Theorem \ref{cu6thm1} shows that there exist
$A_2>0$ and $\ka>1$ such that $\lnm{\ve^t}{2m/(m+2)}\le
A_2t^{\ka+m/2}$, $\cnm{\ve^t}{0}\le A_2t^{\ka-1}$, and
$\lnm{\d\ve^t}{2m}\le A_2t^{\ka-3/2}$ for all $t\in(0,\de)$.
It remains to consider the condition~$\lnm{\pi_\sWt(\ve^t)}{1}\le
A_2t^{\ka+m-1}$.

Combining equations \eq{cu7eq1} and \eq{cu7eq16} and
Propositions \ref{cu7prop2} and \ref{cu7prop3} gives
\e
\int_{N^t}w_{\bf d}^t\ve^t\,\d V^t=O\bigl(\md{{\bf d}}t^{(m+1)\tau}
\bigr)\quad\text{for all ${\bf d}\in\R^q$ and $t\in(0,\de)$.}
\label{cu7eq25}
\e
One can show from Definition \ref{cu7def2} that
$\lnm{w_{\bf d}^t}{2}\ge C\md{{\bf d}}$ for some $C>0$ and
all ${\bf d}\in\R^q$ and $t\in(0,\de)$. This and \eq{cu7eq25}
imply that $\lnm{\pi_\sWt(\ve^t)}{2}=O(t^{(m+1)\tau})$.
But $\lnm{\pi_\sWt(\ve^t)}{1}\le \vol(N^t)^{1/2}
\lnm{\pi_\sWt(\ve^t)}{2}$, and $\vol(N^t)=O(1)$. Therefore
\e
\lnm{\pi_\sWt(\ve^t)}{1}=O(t^{(m+1)\tau})
\quad\text{for all $t\in(0,\de)$.}
\label{cu7eq26}
\e

Make $\ka>1$ smaller if necessary so that $\ka+m-1\le(m+1)\tau$.
This is possible as $(m+1)\tau>m$ by \eq{cu7eq2}. Now make $A_2>0$
bigger and $\de>0$ smaller if necessary so that $\lnm{\pi_{\sst
W^t}(\ve^t)}{1}\le A_2t^{\ka+m-1}$ for all $t\in(0,\de)$. This is
possible by \eq{cu7eq26}, as $\ka+m-1\le(m+1)\tau$. The previous
inequalities $\lnm{\ve^t}{2m/(m+2)}\le A_2t^{\ka+m/2}$,
$\cnm{\ve^t}{0}\le A_2t^{\ka-1}$ and $\lnm{\d\ve^t}{2m}\le
A_2t^{\ka-3/2}$ for $t\in(0,\de)$ still hold with the new
$A_2,\ka$, as we have increased $A_2$ and decreased $\ka$, and
$t<1$. Thus there exist $A_2,\ka$ satisfying the conditions of
the theorem.
\end{proof}

\subsection{Parts (vi) and (vii) of Theorem \ref{cu5thm2}}
\label{cu73}

We now explain how to modify the material of \S\ref{cu64} to
the case when $X'$ is not connected. Thus we prove that parts
(vi) and (vii) of Theorem \ref{cu5thm2} hold for $N^t$ and
$W^t$. Here is the analogue of Proposition~\ref{cu6prop3}.

\begin{prop} In the situation of\/ \S\ref{cu71}, there exists
$D_2>0$ such that for all\/ $v\in C^1_{\rm cs}(X')$ we have
\e
\lnm{v}{2m/(m-2)}\le D_2\bigl(\lnm{\d v}{2}+
\ts\sum_{k=1}^q\bmd{\int_{X'_k}v\,\d V_g\,}\,\bigr).
\label{cu7eq27}
\e
\label{cu7prop4}
\end{prop}

\begin{proof} Applying Proposition \ref{cu6prop3} to
each connected component $X'_k$ of $X'$ gives constants
$D_{2,k}>0$ for $k=1,\ldots,q$ such that
\e
\blnm{v\vert_{X_k'}}{2m/(m-2)}\le D_{2,k}\bigl(
\blnm{\,\d v\vert_{X_k'}}{2}+\bmd{\ts\int_{X'_k}v\,\d V_g\,}\,\bigr).
\label{cu7eq28}
\e
Then summing \eq{cu7eq28} over $k=1,\ldots,q$ and using
\begin{equation*}
\ts\blnm{v}{2m/(m-2)}\le\sum_{k=1}^q\blnm{v\vert_{X_k'}}{2m/(m-2)}
\;\>\text{and}\;\>
\ts\sum_{k=1}^q\blnm{\,\d v\vert_{X_k'}}{2}\le q^{1/2}\blnm{\,\d v}{2}
\end{equation*}
proves \eq{cu7eq27}, with~$D_2=q^{1/2}\max
\bigl(D_{2,1},\ldots,D_{2,q})$.
\end{proof}

Here is the analogue of Theorem \ref{cu6thm4}. The condition
$\int_{N^t}vw\,\d V^t=0$ for all $w\in W^t$ is equivalent to
$\pi_\sWt(v)=0$, so the theorem proves part (vi) of Theorem
\ref{cu5thm2} for $N^t,W^t$, with $A_7$ independent of~$t$.

\begin{thm} Making $\de>0$ smaller if necessary, there exists
$A_7>0$ such that for all\/ $t\in(0,\de)$, if\/ $v\in L^2_1(N^t)$
with\/ $\int_{N^t}vw\,\d V^t=0$ for all\/ $w\in W^t$ then $v\in
L^{2m/(m-2)}(N^t)$ and\/~$\lnm{v}{2m/(m-2)}\le A_7\lnm{\d v}{2}$.
\label{cu7thm2}
\end{thm}

\begin{proof} Let $a,b,G$ and $F^t$ be as in the proof of Theorem
\ref{cu6thm4}. Then $F^t:N^t\ra[0,1]$ is smooth with $F^t\equiv 1$
on $K$, and $F^t\equiv 0$ on $P_i^t$ for $i=1,\ldots,n$. The support
of $F^t$ is a subset of $N^t\cap X'$, so we can also regard $F^t$ as
a compactly-supported function on $X'$. For $k=1,\ldots,q$, let
$F^t_k$ be the smooth, compactly-supported function on $X'$ equal
to $F^t$ on $X'_k$ and zero on $X'_{k'}$ for $k'\ne k$. Then
$F^t_k$ is supported on $X'_k$, and $F^t=\sum_{k=1}^qF^t_k$.
Moreover, extending $F^t_k$ by zero outside $N^t\cap X'$ we can
also regard $F^t_k$ as a smooth, compactly-supported function
on $N^t$, and $F^t=\sum_{k=1}^qF^t_k$ holds on $N^t$ as well.

Suppose now that $t\in(0,\de)$ and $v\in C^1(N^t)$ with
$\int_{N^t}vw\,\d V^t=0$ for all $w\in W^t$. Then $F^tv$
is supported in $N^t\cap X'$, so we can regard it as a
compactly-supported function on $X'$ and apply Proposition
\ref{cu7prop4} to it. This gives
\e
\begin{split}
\lnm{F^tv}{2m/(m-2)}
&\le D_2\bigl(\blnm{\,\d(F^tv)}{2}+
\ts\sum_{k=1}^q\bmd{\int_{X'_k}F^tv\,\d V_g\,}\,\bigr)\\
&=D_2\bigl(\blnm{F^t\d v+v\,\d F^t}{2}+
\ts\sum_{k=1}^q\bmd{\int_{N^t}F^t_kv\,\d V^t\,}\,\bigr).
\end{split}
\label{cu7eq29}
\e
Here in the second line we use the fact that $F^t$ is
supported in $N^t\cap X'$, so $h^t=g\vert_{X'}$ and
$\d V_g=\d V^t$ in the support of~$F^t$.

Let $e_1,\ldots,e_q$ be the usual basis of $\R^q$, so that
$e_k=(\de_{1k},\de_{2k},\ldots,\de_{qk})$ for $k=1,\ldots,q$.
As $\int_{N^t}vw_{e_k}^t\,\d V^t=0$ by choice of $v$ we have
\begin{equation*}
\ts\bmd{\int_{N^t}F^t_kv\,\d V^t\,}=
\ts\bmd{\int_{N^t}(F^t_k-w_{e_k}^t)v\,\d V^t\,}\le
\blnm{v}{2m/(m-2)}\cdot\blnm{F^t_k-w_{e_k}^t}{2m/(m+2)}.
\end{equation*}
Substituting this into \eq{cu7eq29} and using H\"older's
inequality we get
\e
\begin{split}
\lnm{F^tv}{2m/(m-2)}\le D_2\bigl(&\lnm{F^t\d v}{2}+
\lnm{v}{2m/(m-2)}\cdot\lnm{\d F^t}{m}\\
&+\blnm{v}{2m/(m-2)}\cdot
\ts\sum_{k=1}^q\blnm{F^t_k-w_{e_k}^t}{2m/(m+2)}\bigr).
\end{split}
\label{cu7eq30}
\e
This is the analogue of~\eq{cu6eq31}.

Now by the definitions of $F^t_k$ and $w_{e_k}^t$ we have
\begin{equation*}
F^t_k(x)=w_{e_k}^t(x)=\begin{cases}
1, & \text{$x\in K\cap X'_k$},\\
1, & \text{$x\in \Xi_i^t\bigl(\Si_i^j\t[t^a,R')\bigr)$
when $k(i,j)=k$,}\\
0, & \text{$x\in K\cap X'_{k'}$ for $k'\ne k$},\\
0, & \text{$x\in \Xi_i^t\bigl(\Si_i^j\t[t^a,R')\bigr)$
for $k(i,j)\ne k$.}
\end{cases}
\end{equation*}
Hence $F^t_k-w_{e_k}^t$ is {\it zero} on most of $N^t$.
The support of $F^t_k-w_{e_k}^t$ is contained in the union
of $P_i^t$ and $\Xi_i^t\bigl(\Si_i\t(tT,t^a)\bigr)$ over
$i=1,\ldots,n$, which has volume $O(t^{ma})$, and here
$\md{F^t_k-w_{e_k}^t}\le 1$ as $0\le F^t_k,w_{e_k}^t\le 1$. Hence
\e
\blnm{F^t_k-w_{e_k}^t}{2m/(m+2)}=O(t^{a(m+2)/2})
\quad\text{for $k=1,\ldots,q$.}
\label{cu7eq31}
\e

Following the proof of Theorem \ref{cu6thm4} without
change, we prove \eq{cu6eq34}. Using \eq{cu7eq30}
instead of \eq{cu6eq31}, in place of \eq{cu6eq35}
we obtain
\begin{align*}
\bigl[1-(D_2+2\sqrt{n}\,D_1)&\lnm{\d F^t}{m}
-D_2\ts\sum_{k=1}^q\blnm{F^t_k-w_{e_k}^t}{2m/(m+2)}
\bigr]\cdot\lnm{v}{2m/(m-2)}\\
&\le(D_2+2\sqrt{n}\,D_1)\lnm{\d v}{2}.
\end{align*}
Using \eq{cu7eq31} instead of $\lnm{1-F^t}{2m/(m+2)}=O(t^{a(m+2)/2})$,
the rest of the proof follows that of Theorem~\ref{cu6thm4}.
\end{proof}

We now prove part (vii) of Theorem \ref{cu5thm2} for $N^t,W^t$,
with $A_7$ as above.

\begin{thm} Making $\de>0$ smaller if necessary, for all\/
$t\in(0,\de)$ and\/ $w\in W^t$ we have $\lnm{\d^*\d w}{2m/(m+2)}
\le\ha A_7^{-1}\lnm{\d w}{2}$, where $A_7>0$ is as in Theorem
\ref{cu7thm2}. Also there exists $A_8>0$ such that for all\/
$t\in(0,\de)$ and\/ $w\in W^t$ with\/ $\int_{N^t}w\,\d V^t=0$
we have~$\cnm{w}{0}\le A_8t^{1-m/2}\lnm{\d w}{2}$.
\label{cu7thm3}
\end{thm}

\begin{proof} Following the method of Proposition \ref{cu6prop2},
using the estimates of Proposition \ref{cu7prop1} and taking
$2-m<\be<\ha(2-m)$ we find that for some $D_3>0$ and all
$t\in(0,\de)$ and $w_{\bf d}^t\in W^t$ we have
\ea
\lnm{\d w_{\bf d}^t}{2}&\le D_3\md{{\bf d}}t^{(m-2)/2}
\quad\text{and}
\label{cu7eq32}\\
\lnm{\d^*\d w_{\bf d}^t}{2m/(m+2)}&\le D_3\md{{\bf d}}t^{m/2}
+D_3\md{{\bf d}}t^{2\be(\tau-1)+\tau(m-2)}.
\label{cu7eq33}
\ea
Also, as $w_{\bf d}^t=d_k$ on $K\cap X_k'$, from \eq{cu7eq5} we
deduce that
\e
\cnm{w_{\bf d}^t}{0}=\max\bigl(\md{d_1},\ldots,\md{d_q})\le\md{{\bf d}}.
\label{cu7eq34}
\e

However, \eq{cu7eq32}--\eq{cu7eq34} are not enough to prove
the theorem, as \eq{cu7eq32} gives an {\it upper bound\/} for 
$\lnm{\d w_{\bf d}^t}{2}$, but we actually need a {\it lower bound}.
Now on $P_i^t$ and $\Xi_i^t\bigl(\Si_i\t(tT,\ha t^\tau]\bigr)$, by
definition $w_{\bf d}^t$ coincides with $v^{{\bf c}_i}_i$ on
the corresponding regions of $L_i$. Using \eq{cu7eq11} to
compare the volume forms $\d V^t$ of $h^t$ on $N^t$ and
$\d V_{\smash{h_i}}$ of $h_i=g'\vert_{L_i}$ on $L_i$, from
\eq{cu7eq6}, \eq{cu7eq7} and \eq{cu7eq10} we can show that
\e
\lnm{\d w_{\bf d}^t}{2}^2\!=\!t^{m-2}\sum_{i=1}^n
\lnm{\d v^{{\bf c}_i}_i}{2}^2\!+\!O\bigl(\ms{{\bf d}}t^{m-1}\bigr)
\!+\!O\bigl(\ms{{\bf d}}t^{\be(\tau-1)+\tau(m-2)/2}\bigr),
\label{cu7eq35}
\e
where $\lnm{\d w_{\bf d}^t}{2}$ is computed on $N^t$ using $h^t$,
and $\lnm{\d v^{{\bf c}_i}_{\smash{i}}}{2}$ on $L_i$ using~$h_i$.

Now $v^{{\bf c}_i}_{\smash{i}}$ depends only on $\bf d$ and not
on $t$, and in fact only on $d_{k(i,j)}$ for $j=1,\ldots,l_i$.
Also $\lnm{\d v^{{\bf c}_{\smash{i}}}_i}{2}$ depends only on
the differences $d_{k(i,j)}-d_{k(i,j')}$ for $1\le j<j'\le l_i$,
as adding an overall constant to $d_k$ adds a constant to 
$v^{{\bf c}_{\smash{i}}}_i$, and does not change $\d v^{{\bf
c}_{\smash{i}}}_i$.

For each $i=1,\ldots,n$, define an equivalence relation
$\sim$ on pairs $(i,j)$ for $j=1,\ldots,l_i$ by $(i,j)\sim
(i,j')$ if $\Si_i^j$ and $\Si_i^{\smash{j'}}$ are ends at
infinity of the same connected component of $L_i$. Then one
can show that
\e
\lnm{\d v^{{\bf c}_i}_i}{2}^2\ge
C_i\sum_{\substack{1\le j<j'\le l_i:\\(i,j)\sim(i,j')}}
\ms{d_{k(i,j)}-d_{k(i,j')}},
\label{cu7eq36}
\e
for all ${\bf d}\in\R^q$ and $C_i>0$ depending only on $L_i$
for~$i=1,\ldots,n$.

The point here is that if $d_{k(i,j)}=d_{k(i,j')}$ whenever
$(i,j)\sim(i,j')$ then $v^{{\bf c}_i}_i$ is constant on each
connected component of $L_i$, and $\d v^{{\bf c}_i}_i\equiv 0$,
so both sides of \eq{cu7eq36} are zero. But otherwise
$v^{{\bf c}_i}_i$ is not constant, and both sides of \eq{cu7eq36}
are positive. Summing \eq{cu7eq36} over $i=1,\ldots,n$, and
remembering that the $N^t$ are connected by Definition
\ref{cu7def1}, we can prove that
\e
\sum_{i=1}^n\lnm{\d v^{{\bf c}_i}_i}{2}^2\ge 2D_4^2
\bigl(\ms{{\bf d}}-\ts\frac{1}{q}(d_1+\ldots+d_q)^2\bigr),
\label{cu7eq37}
\e
for all ${\bf d}\in\R^q$, and some $D_4>0$ depending only on
$C_1,\ldots,C_n$, the $k(i,j)$, and the equivalence relations~$\sim$.

Combining \eq{cu7eq35} and \eq{cu7eq37} shows that for all
${\bf d}\in\R^q$ with $d_1+\cdots+d_q=0$ and $t\in(0,\de)$ we have
\begin{equation*}
\lnm{\d w_{\bf d}^t}{2}^2\ge 2D_4^2\ms{{\bf d}}t^{m-2}
+O\bigl(\ms{{\bf d}}t^{m-1}\bigr)
+O\bigl(\ms{{\bf d}}t^{\be(\tau-1)+\tau(m-2)/2}\bigr).
\end{equation*}
As $\be<\ha(2-m)$ and $\tau\in(0,1)$, both error terms are
smaller than $\ms{{\bf d}}t^{m-2}$ for small $t$. Hence, making
$\de>0$ smaller if necessary we see that when~$t\in(0,\de)$,
\e
\lnm{\d w_{\bf d}^t}{2}\ge D_4\md{{\bf d}}t^{(m-2)/2}
\quad\text{for ${\bf d}\in\R^q$ with $d_1+\cdots+d_q=0$.}
\label{cu7eq38}
\e

Now make $\de>0$ smaller if necessary so that for all
$t\in(0,\de)$, we have
\e
D_3\md{{\bf d}}t^{m/2}
+D_3\md{{\bf d}}t^{2\be(\tau-1)+\tau(m-2)}\le
\ha A_7^{-1}D_4\md{{\bf d}}t^{(m-2)/2}.
\label{cu7eq39}
\e
This is possible as both powers of $t$ on the left
are greater than the power on the right. Combining
\eq{cu7eq33}, \eq{cu7eq38} and \eq{cu7eq39} shows that
when~$t\in(0,\de)$,
\e
\lnm{\d^*\d w_{\bf d}^t}{2m/(m+2)}\le\ha
A_7^{-1}\lnm{\d w_{\bf d}^t}{2}
\;\>\text{for ${\bf d}\!\in\!\R^q$ with
$d_1\!+\!\cdots\!+\!d_q\!=\!0$.}
\label{cu7eq40}
\e

Now let $t\in(0,\de)$ and $w\in W^t$. Then $w=w_{{\bf d}'}^t$ for
some ${\bf d}'\in\R^q$. Let $c=\frac{1}{q}(d_1'+\ldots+d_q')$,
let $d_k=d_k'-c$, and ${\bf d}=(d_1,\ldots,d_q)$. Then
$d_1+\ldots+d_q=0$, and $w=w_{\bf d}^t+c$, since $w_{\bf d}^t$
depends linearly on $\bf d$ and $w_{\smash{(1,\ldots,1)}}^t\equiv
1$. Thus $\d w=\d w_{\bf d}^t$, and \eq{cu7eq40} holds, giving
$\lnm{\d^*\d w}{2m/(m+2)}\le\ha A_7^{-1}\lnm{\d w}{2}$, as
we have to prove.

Suppose also that $\int_{N^t}w\,\d V^t=0$. Then by
the proof of \eq{cu7eq34} we see that
\begin{equation*}
\min(d_1',\ldots,d_q')=\min_{N^t}w\le 0
\le\max_{N^t}w=\max(d_1',\ldots,d_q').
\end{equation*}
Using this it is easy to show that
\e
\cnm{w}{0}=\max\bigl(\md{d_1'},\ldots,\md{d_q'})\le 2
\max\bigl(\md{d_1},\ldots,\md{d_q})\le 2\md{{\bf d}}.
\label{cu7eq41}
\e
Define $A_8=2D_4^{-1}>0$. Then using \eq{cu7eq38}, 
\eq{cu7eq41} and $\d w=\d w_{\bf d}^t$ we find that
\begin{equation*}
\cnm{w}{0}\le 2\md{{\bf d}}=A_8t^{1-m/2}\cdot D_4\md{{\bf d}}t^{(m-2)/2}
\le A_8t^{1-m/2}\lnm{\d w}{2},
\end{equation*}
for all $t\in(0,\de)$ and $w\in W^t$ with $\int_{N^t}w\,\d V^t=0$.
This completes the proof.
\end{proof}

\subsection{The main results, when $X'$ is not connected}
\label{cu74}

We can now state our second main result on
desingularizations of SL $m$-folds $X$ with conical
singularities, this time allowing $X'$ not connected.

\begin{thm} Suppose $(M,J,\om,\Om)$ is an almost Calabi--Yau
$m$-fold and\/ $X$ a compact SL\/ $m$-fold in $M$ with conical
singularities at\/ $x_1,\ldots,x_n$ and cones $C_1,\ldots,C_n$.
Define $\psi:M\ra(0,\iy)$ as in \eq{cu2eq3}. Let\/ $L_1,\ldots,L_n$
be Asymptotically Conical SL\/ $m$-folds in $\C^m$ with cones
$C_1,\ldots,C_n$ and rates $\la_1,\ldots,\la_n$. Suppose
$\la_i<0$ for $i=1,\ldots,n$. Write $X'=X\sm\{x_1,\ldots,x_n\}$
and\/~$\Si_i=C_i\cap{\cal S}^{2m-1}$.

Set\/ $q=b^0(X')$, and let\/ $X_1',\ldots,X_q'$ be the connected
components of\/ $X'$. For $i=1,\ldots,n$ let\/ $l_i=b^0(\Si_i)$,
and let\/ $\Si_i^1,\ldots,\Si_i^{\smash{l_i}}$ be the connected
components of\/ $\Si_i$. Define $k(i,j)=1,\ldots,q$ by $\Up_i
\circ\vp_i\bigl(\Si_i^j\t(0,R')\bigr)\subset X'_{\smash{k(i,j)}}$
for $i=1,\ldots,n$ and $j=1,\ldots,l_i$. Suppose that
\e
\sum_{\substack{1\le i\le n, \; 1\le j\le l_i: \\
k(i,j)=k}}\psi(x_i)^mZ(L_i)\cdot[\Si_i^j\,]=0
\quad\text{for all\/ $k=1,\ldots,q$.}
\label{cu7eq42}
\e

Suppose also that the compact\/ $m$-manifold\/ $N$ obtained by
gluing $L_i$ into $X'$ at\/ $x_i$ for $i=1,\ldots,n$ is connected.
A sufficient condition for this to hold is that\/ $X$ and\/ $L_i$
for $i=1,\ldots,n$ are connected.

Then there exists $\ep>0$ and a smooth family $\smash{\bigl\{
\ti N^t:t\in(0,\ep]\bigr\}}$ of compact, nonsingular SL\/
$m$-folds in $(M,J,\om,\Om)$ diffeomorphic to $N$, such
that\/ $\smash{\ti N^t}$ is constructed by gluing $tL_i$
into $X$ at\/ $x_i$ for $i=1,\ldots,n$. In the sense of
currents in Geometric Measure Theory, $\smash{\ti N^t}\ra
X$ as~$t\ra 0$.
\label{cu7thm4}
\end{thm}

The proof follows that of Theorem \ref{cu6thm5}, but
using Definition \ref{cu7def1} instead of Definitions
\ref{cu6def1}--\ref{cu6def3}, and defining $W^t$ as in
Definition \ref{cu7def2} rather than $W^t=\an{1}$. We
use Theorem \ref{cu7thm1} instead of Theorem \ref{cu6thm1},
and Theorems \ref{cu7thm2} and \ref{cu7thm3} instead of Theorem
\ref{cu6thm4}. Note that Theorem \ref{cu6thm2} still holds in
this situation, as it does not assume that $X'$ is connected.
The extra hypotheses \eq{cu7eq42} and that $N$ is connected
come from Definition~\ref{cu7def1}.

If $X'$ is connected, so that $q=1$, then $k(i,j)\equiv 1$
and \eq{cu7eq42} becomes
\begin{equation*}
\sum_{i=1}^n\psi(x_i)^mZ(L_i)\cdot\sum_{j=1}^{l_i}[\Si_i^j\,]=0.
\end{equation*}
But $\sum_{j=1}^{l_i}[\Si_i^j\,]=[\Si_i]$, and $Z(L_i)\cdot[\Si_i]=0$
as $Z(L_i)$ is the image of a class in $H^{m-1}(L_i,\R)$ by Definition
\ref{cu4def2}, and $\Si_i$ is the boundary of $L_i$, so $[\Si_i]$ maps
to zero in $H_{m-1}(L_i,\R)$. Therefore \eq{cu7eq42} holds
automatically when $X'$ is connected, and Theorem \ref{cu7thm4}
reduces to Theorem \ref{cu6thm5} in this case.

\end{document}